\documentclass[a4paper]{article}
\usepackage{RR}
\usepackage{hyperref}
\usepackage[T1]{fontenc}
\RRdate{March 2009}

\RRauthor{
Vincent Acary\thanks{vincent.acary@inrialpes.fr}%
  \and
Bernard Brogliato\thanks{Bernard.brogliato@inrialpes.fr}%
}
\authorhead{Acary \& Brogliato}
\RRtitle{Simulations numériques par la méthode d'Euler implicite des systèmes à modes glissants}
\RRetitle{Implicit Euler numerical simulations of sliding mode systems}
\titlehead{Implicit Euler numerical simulations of sliding mode systems}
\RRresume{Dans ce rapport, on montre que la discrétisation en temps de type Euler implicite conduit à une très bonne stabilisation d'une classe de systèmes commutés avec des modes glissants, et de leurs dérivées sur la surface de glissement. Les oscillations artificielles qui sont généralement mentionnées pour l'implémentation discréte de ce type de systèmes sont évitées. De plus, la méthode (de type ``event-capturing'' ou ``time--stepping'') permet de traiter des accumulations d'événements (Phénomène de Zenon) et des surfaces de commutations multiples (\textit{i.e.} des surfaces de glissement de codimension $\geq 2$). Dans ce rapport,  les détails de l'implémentation  sont donnés et des exemples numériques illustrent  ses propriétés. Cette méthode peut être une alternative aux méthodes complexes de suppression des oscillations, en gardant la nature intrinsèquement discontinue de la dynamique sur les surfaces de glissement. Le lien avec les commandes à modes glissants en temps discret est étudié.
}
\RRabstract{In this report it is shown that the implicit Euler time-discretization of some classes of switching systems with sliding modes, yields a very good stabilization of the trajectory and of its derivative on the sliding surface. Therefore the spurious oscillations which are pointed out elsewhere when an explicit method is used, are avoided. Moreover the method (an {\em event-capturing}, or {\em time-stepping} algorithm) allows for accumulation of events (Zeno phenomena) and for multiple switching surfaces (i.e., a sliding surface of codimension $\geq 2$). The details of the implementation are given, and numerical examples illustrate the developments. This method may be an alternative method for chattering suppression, keeping the intrinsic discontinuous nature of the dynamics on the sliding surfaces. Links with discrete-time sliding mode controllers are studied. }
\RRmotcle{Systèmes commutés, Inclusion Différentielles de Filippov, Problèmes de complémentarité, Méthode d'Euler implicite, modes glissants, opérateurs, maximaux monotones, Problème linéaire de complémentarité mixte, Discrétisation Bloqueur d'Ordre Zéro (BOZ)
}
\RRkeyword{Switching systems, Filippov's differential inclusions, complementarity problems, backward Euler algorithm, sliding modes, maximal monotone mappings, mixed linear complementarity problem, ZOH discretization. 
}
\RRprojet{Bipop}
\RRtheme{\THNum} 
 \RCGrenoble 

\usepackage{a4wide}
\usepackage{amsmath}
\usepackage{amssymb}
\usepackage{color}
\usepackage{graphicx,epsfig}
\usepackage[latin1]{inputenc}
\usepackage{subfigure}
\usepackage{psfrag}

\usepackage{url}

\DeclareMathOperator{\sgn}{sgn}
\usepackage{algorithmic}
\usepackage{algorithm}

\newcommand{\scontract}{{\,{\Bar\otimes}\,}}

\usepackage{fancyhdr}\usepackage{lastpage}

\begingroup
\def\2#1{\ifnum#1<10 0\fi\the#1}
\xdef\isodayandtime{\the\year-\2\month-\2\day\space\2{\count0}:%
\2{\count2}}
\endgroup

{\newtheorem{ndr1va}{\textbf{\textsc{Redaction note V.A.}}}[section]}

{%
\noindent\begin{ndr1va}\hrule\vspace{1em}%
\ttfamily\small \  \\
\indent}%
{%
\begin{flushright}%
\  \\
\end{flushright}%
\vspace{-1.5em}\hrule
\end{ndr1va}%
}

\usepackage{pifont}
\def\cqfd{\ifmmode\sqw\else{\ifhmode\unskip\fi\nobreak\hfil
\penalty50\hskip1em\null\nobreak\hfil\ding{111}
\parfillskip=0pt\finalhyphendemerits=0\endgraf}\fi}

\makeatother

\def\geq{\geqslant}
\def\leq{\leqslant}

\newtheorem{definition}{Definition}
\newtheorem{proposition}{Proposition}
\newtheorem{lemma}{Lemma}

\newtheorem{remark}{Remark}
\newtheorem{assumption}{Assumption}
\newtheorem{example}{Example}

\newtheorem{corollary}{Corollary}

\newcommand{\ZZ}{\mbox{\rm \lower0.3pt\hbox{$\angle\!\!\!$}Z}}
\newcommand{\RR}{\mbox{\rm $I\!\!R$}}
\newcommand{\NN}{\mbox{\rm $I\!\!N$}}
\newcommand{\Qend}{\hfill $\Box$}

\begin{document}
\RRNo{6886}
\makeRR   

\thispagestyle{plain}

\section{Introduction}

Sliding mode controllers are widely used because of their intrinsic robustness properties \cite{Orlov2009,LNCIS334,LNCIS375,tsai}. Some important fields of application are induction motors \cite{proca2003,wang2005,barambones}, aircraft control \cite{promtun2009,hess,wells,koren}, hard disk drives \cite{hu2009,herrmann2005}, solar energy systems \cite{garcia2009}. However they are known to generate chattering which renders their application delicate. Solutions to cope with chattering or reduce its effects have been proposed, see e.g. \cite{alexander1998,alexander1999,boiko,brogliato1995,tsai,young1999}, which also have their own limitations \cite{young1999}. One drawback of these solutions is that they usually destroy the intrinsic discontinuous nature of sliding mode control. Fundamentally, these control schemes are of the switching discontinuous type and they yield closed-loop systems that can be recast into Filippov's differential inclusions. The numerical simulation of such nonsmooth dynamical systems is non trivial and it has received a lot of attention, see e.g. \cite{stewart1996,stewart1990,camlibel2002,kastner,dontchev,leine2004,danca}, to cite a few.  Both event-driven methods and time-stepping methods have been developed, see e.g. \cite{acary2008} for a survey. In this paper we focus on time-stepping methods, which have an interest not only for the sake of numerical simulation, but also for the real implementations of sliding mode controllers on discrete-time systems \cite{wang}. Recently it has been shown that the {\em explicit} Euler method generates unwanted effects like spurious oscillations (also called chattering effects) around the switching surface \cite{galias2006,galias2007,wang,yu2008}. In parallel, the digital implementation of sliding mode controllers has been studied thoroughly in \cite{Galias2008a,koshkouei2000}, where the Zero--Order Holder (ZOH) discretization is used.

The purpose of this paper is to analyze the {\em implicit} ({\em backward}) Euler method for some particular classes of differential inclusions, that include sliding mode controllers. It is shown that, besides convergence and order results, the advantage of the implicit method is that it allows one to get a very accurate and smooth stabilization on the switching surface (of codimension one or larger than one). Roughly speaking, this is due to the fact that the switches are no longer monitored by the state at step $k$, but by a {\em multiplier} (a slack variable in a nonlinear programming language). The multivalued part of the sgn$(\cdot)$ function, i.e. a multifunction, is then correctly taken into account, avoiding stiff problems. The advantage of such ``dual'' methods in terms of their accuracy on the sliding surface has already been noticed in \cite{stewart1990,stewart1996} in an event-driven context, where the motivation was the simulation of mechanical systems with Coulomb friction. From a numerical point of view, our study shows that convergence and order results may not be sufficient to guarantee that the derivative of the state is correctly approximated on the switching surface. The implicit method adapts naturally to an arbitrary large number of switching surfaces, that is not the case of most of the other methods which  become quite cumbersome as soon as more than two switching surfaces are considered. A further advantage of the proposed method is that contrary to other methods that have been studied and which destroy the intrinsic discontinuous nature of sliding mode systems \footnote{see \cite{young1999} for a discussion on this point.} (like the so-called {\em boundary layer control}, or various filtering techniques), our method keeps the multivalued discontinuity and consequently the fundamental aspects and properties of sliding mode control from a Filippov's systems point of view.  Moreover, sampling rates need not be high to reduce chattering, contrary to other discrete sliding mode controllers. A second contribution of this paper is to show that the results that hold for the backward Euler scheme, extend to ZOH discretizations of sliding mode systems.

The paper is organized as follows: Section~\ref{section1} presents a motivating example for using an implicit Euler implementation of the simplest sliding mode system. In Section~\ref{section2}, a class of differential inclusions is introduced and existence and uniqueness results are given under the  maximal monotonicity assumption. Through several examples,  the Equivalent--Control--Based Sliding--Mode--Control (ECB-SMC)  and the Lyapunov--based discontinuous robust control are shown to fit well within this class of differential inclusion. In Section~\ref{section3}, some convergence and chattering free finite--time stabilization results are given. These central results of the paper show that the implicit Euler implementation of the differential inclusion yields a chattering free convergence in finite time on the sliding surface. Section~\ref{sec:DTSMC} is devoted the study of Discrete--time Sliding Mode Control and the extension to ZOH discretization. Some hints on the numerical implementation of the implicit Euler scheme are given in Section \ref{Sec:TimeDiscretization} and the paper ends with some numerical experiments in Section~\ref{secSimus}.

{\bf Notations and definitions:} Let $A\in \RR^{n\times m}$, then  $A_{\bullet i}$ is the $i$th column and  $A_{i \bullet}$ is the $i$th row. The open ball of radius $r > 0$ centered at a point $x\in \RR^n$ is denoted by $B_r(x)$. For a set of indices $\alpha \subset \{1,\ldots, n\}$ and a column vector $x\in \RR^n$, the column vector $x_\alpha$ will denoted the sub-vector of corresponding indices in $\alpha$, that is $x_{\alpha} = [x_i,i \in \alpha]^T$.

\section{A simple example}
\label{section1}

To start with we consider the simplest case: 
\begin{equation}
\dot{x}(t) \in -\mbox{sgn}(x(t))=\left\{\begin{array}{lll} 1 & \mbox{if} & x(t) <0 \\ -1 & \mbox{if} & x(t) >0 \\ \mbox{[-1,1]} & \mbox{if} & x(t)=0 \end{array}\right.,\;\;x(0)=x_{0}
\label{eq:sept}
\end{equation}
with $x(t) \in \RR$. This system possesses a unique Lipschitz continuous solution for any $x_{0}$. The backward Euler discretization of (\ref{eq:sept}) reads as:
\begin{equation}
 \begin{cases}
   x_{k+1}-x_{k} = -h s_{k+1}\\[2mm]
   s_{k+1} \in \mbox{sgn}(x_{k+1})
 \end{cases}
\label{eq:huit}
\end{equation}
This method converges with at least order $\frac{1}{2}$ (see Proposition \ref{proposition2} below). Let us now state a result which shows that once the iterate $x_{k}$ has reached a value inside some threshold around zero for some $k$, then the dual variable $s_{k+1}$ keeps its value and so does $x_{k+n}$ for all $n \geq 1$. 

\begin{lemma}
For all $h>0$ and $x_0 \in \RR$, there exists $k_0$ such that  $x_{k_0+n}=0$ and  $\displaystyle\frac{x_{k_0+n+1}-x_{k_0+n}}{h}=0$  for all $n \geq 1$. 
\label{lemmayahoo}
\end{lemma}

{\bf Proof:} The value $k_0$ is defined as the first time step such that $x_{k_0} \in [-h , h] $. If $x_0 \in [-h , h]$, then $k_0=0$. Otherwise, the solution of the time-discretization~(\ref{eq:huit}) is given by $x_k = x_0-\sgn(x_0) k h, s_k = \sgn(x_o)$ while $x_k\notin [-h , h]  $ for $ k < k_0$, and $k_0 = \lceil \frac{|x(0)|}{h}\rceil -1 $. The symbol      $\lceil x  \rceil$ is the  ceiling function which  gives the smallest integer greater than or equal to $x$.

Let us now consider that $x_{k_0}\in [-h, h]$. The only  possible solution for  
\begin{equation}
 \begin{cases}
   x_{k_0+1}-x_{k_0} = -h s_{k_0+1}\\[2mm]
   s_{k_0+1} \in \mbox{sgn}(x_{k_0+1})
 \end{cases}
\label{eq:huit-bis-hypo}
\end{equation}
is $x_{k_0+1} = 0$ and $s_{k_0+1} = \displaystyle\frac{x_{k_0}}{h}$. For the next iteration, we have to solve
\begin{equation}
 \begin{cases}
   x_{k_0+2} = -h s_{k_0+2}\\[2mm]
   s_{k_0+2} \in \mbox{sgn}(x_{k_0+2})
 \end{cases}
\label{eq:huit-ter-hypo}
\end{equation}
and we obtain $x_{k_0+2} = 0 $ and $s_{k_0+2}=0$.  The same holds for all $x_{k_0+n}$,$s_{k_0+n}$, $n \geq 3$, redoing the same reasoning. Clearly then the terms  $\displaystyle\frac{x_{k_0+n+1}-x_{k_0+n}}{h}$ approximating the derivative, are zero for any $h >0$.  \Qend

This result is robust with respect to the numerical threshold that can be encountered in floating point operations. Indeed, let us assume that  $x_{k_0}-h=\varepsilon \ll 1$, that is, $\varepsilon >0 $ is zero at the machine precision. We obtain $s_{k_0+1} = -1 $ and $x_{k_0+1} = \varepsilon $ that is zero at the machine precision. For $n=2$, we obtain $x_{k_0+2} = 0$ and $s_{k_0}= \displaystyle\frac {\varepsilon}{h}$. This robustness stems from the fact  that the dynamics is not only monitored by the sign of $x_k$ but also by the fact that the ``dual'' variable $s_{k+1}$ belongs to $[-1,1]$. 


\begin{figure}[hbtp]
\begin{center}
\input{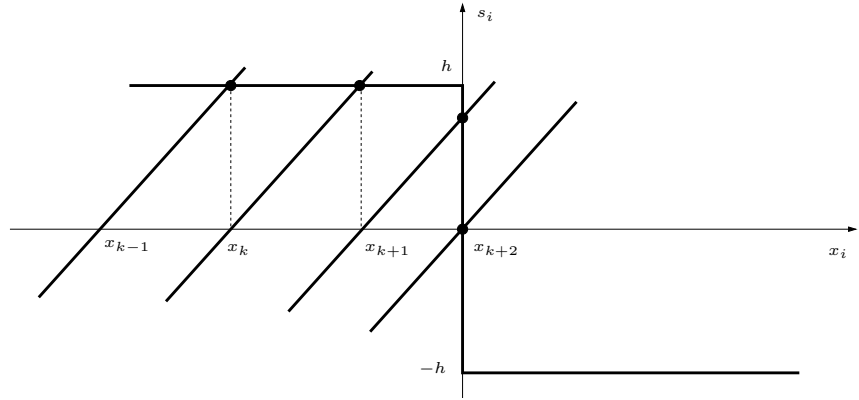}
\caption{Iterations of the backward Euler method.}
\label{figgraph.fig}
\end{center}
\end{figure}

Consequently this result shows that there are no spurious oscillations around the switching surface, contrary to other time-stepping schemes like the explicit Euler method \cite{galias2006,galias2007}. Remarkably Lemma \ref{lemmayahoo} holds for any $h >0$, which means that even a large time step assures a smooth stabilization on the sliding surface. It is noteworthy that solving the system~(\ref{eq:huit})  with unknown $x_{k+1}$ and $s_{k+1}$ is equivalent to calculate the intersection between the graph of the multivalued mapping $x_{k+1} \mapsto -h \mbox{sgn}(x_{k+1})$ and the straight line $x_{k+1} \mapsto x_{k+1}-x_{k}$. This is illustrated on Figure \ref{figgraph.fig}, where  few iterations are depicted until the state reaches zero.

From a control perspective the input is implemented on $[t_{k},t_{k+1})$ as $u_{k}=-\mbox{sgn}(x_{k+1})$ as a function of $x_{k}$ and $h$, where $h$ is the sampling time. There is no problem of causality in such an implementation. It is noteworthy that in the implicit method there is absolutely no issue related to calculating sgn$(0)$, or more exactly sgn$(\epsilon)$ where $\epsilon$ is a very small quantity whose sign is uncertain. The implicit method automatically computes a value inside the multivalued part of the sign multifunction and may be considered as {\em the} time-discretization of the multifunction sgn$(\cdot)$. It is easy to show that the explicit method yields an oscillation around $x=0$, as shown in more general situations in \cite{galias2006,galias2007}. Other time-stepping methods exist, like the so-called {\em switched model} \cite{acary2008,leine2004}, however it fails to correctly solve the integration problem when the number of switched surfaces is too large (see also \cite{alexander1998} for similar issues when the so-called {\em sigmoid blending} mechanism is implemented). Moreover this method may yield a stiff system, and from a control point of view it introduces a high-gain feedback that may not be desirable in practical applications.

\begin{figure}[htbp]
 \centering  
 \psfrag{time}[][]{\Huge $\mathsf t$}
 \psfrag{h}[][]{\Huge $\mathsf h$}
 \psfrag{error}[][]{\Huge $\mathsf{errors}$}
 \psfrag{xv}[][]{}
 \psfrag{Ex1}[r][l]{\huge $\mathsf{x(t)}$}
 \psfrag{Ev1}[r][l]{\huge $\mathsf{-s(t)}$}
 \psfrag{einf}[r][l]{\Huge$\mathsf{(i)}$}
 \psfrag{e1}[r][l]{\Huge$\mathsf{(ii)}\quad\quad$}
 \psfrag{e2}[r][l]{\Huge$\mathsf{(iii)}$}
 \subfigure[ State and Control vs. Time $h = 0.2$]
   {\label{Fig:Trivial2a}\includegraphics[angle=-90,width=0.50\textwidth]{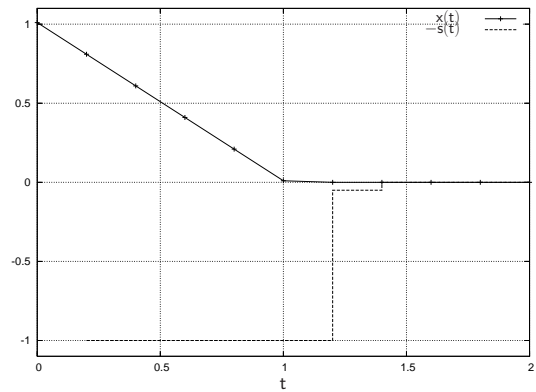}} 
   \subfigure[ State and Control vs. Time $h = 0.02$]
   {\label{Fig:Trivial2b}\includegraphics[angle=-90,width=0.50\textwidth]{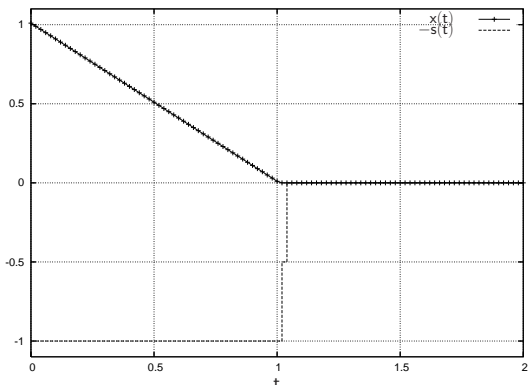}}
  \subfigure[ State and Control vs. Time $h = 0.01$]
  {\label{Fig:Trivial2c}\includegraphics[angle=-90,width=0.49\textwidth]{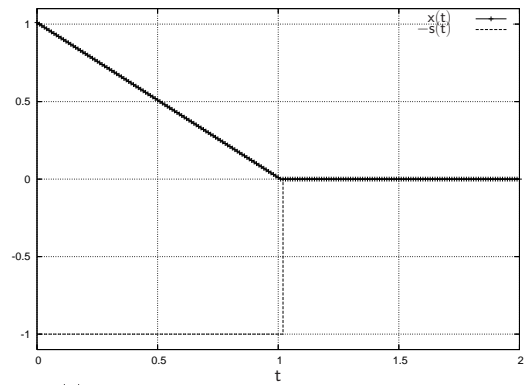}} 
 \subfigure[Numerical error $\|e_s\|_\infty$ (solid line (i)), $\|e_s\|_1$ (dashed line (ii)), $\|e_s\|_2$ (dotted line (iii)),  with respect to $h$ in logscale ]
  {
 \includegraphics[angle=-90,width=0.80\textwidth]{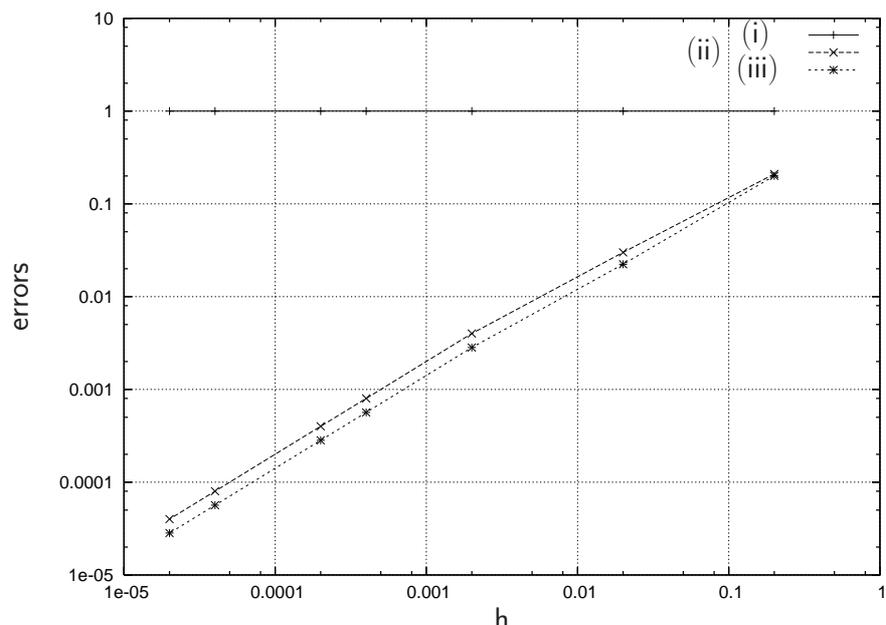}\label{Fig:ConvergenceError}}
 \caption{A simple example for $x_0=1.01$ at $t_0=0$.}\label{Fig:Trivial2}
\end{figure}

On Figure \ref{Fig:Trivial2}(a)-(c), the discrete state $x_k$ and the control $s_k$ are displayed for $x_0=1.01$ at $t_0=0$ and for various values of the time--step $h$ that are sufficiently large to illustrate the behavior of the time--stepping scheme and its convergence.

Let us define two discrete function norms to measure the convergence:
\begin{equation}
 \label{eq:trivial-conv1}
 \begin{array}{l}
 \|e_f\|_{\infty} = \sum_{i=0}^N | f_k - f(t_k)| \\ \\
 \|e_f\|_{p} = ( h \sum_{i=0}^N | f_k - f(t_k)|^p )^{1/p}.
\end{array}
\end{equation}
We can compute that
\begin{equation}
 \label{eq:trivial-conv}
 \|e_s\|_{\infty} = 1 \text{ for all } h>0
\end{equation}
and therefore there is no convergence in infinite norm  $\|.\|_{\infty}$ for $s=\mbox{sgn}(x)$. In   $\|.\|_{1}$ and   $\|.\|_{2}$, we can respectively observe the  convergence with order $1$ on Figure~\ref{Fig:ConvergenceError}.

\paragraph{Complementarity framework} Let us end this section by restating the systems~(\ref{eq:sept}) and (\ref{eq:huit}) into the complementarity framework. Let us introduce equivalent formulations of the inclusion $s(t) \in \sgn(x(t))$ such that
\begin{equation}
  \label{eq:sign}
    s(t) \in \sgn(x(t)) \Leftrightarrow x(t) \in N_{[-1,1]}(  s(t)   ) \Leftrightarrow s(t)\in [-1,1] \text{ and }
   \begin{cases}
     x(t) = 0 \text{ if } s(t) \in ]-1,1[ \\
     x(t) \leq 0 \text{ if } x(t) = - 1 \\
     x(t) \geq 0 \text{ if } x(t) =  1 \\  
   \end{cases}
\end{equation}
where $N_{[-1,1]}$ is the normal cone in the sense of Convex Analysis to the interval $[-1,1]$.  The definition of the normal cone in the present case,
\begin{equation}
  \label{eq:normalcone}
  N_{[-1,1]}(s) = \{ -v_1+v_2 , 0 \leq v_1 \perp s+1 \geq 0, 0 \leq v_2 \perp 1-s \geq 0 \}
\end{equation}
yields the following complementarity representation of the sign multi-valued function
\begin{equation}
  \label{eq:normalcone1}
  x(t) \in N_{[-1,1]}(s(t))  \Leftrightarrow
  \begin{cases}
    x(t) = -v_1(t) + v_2(t) \\
    0 \leq v_1(t) \perp s(t)+1 \geq 0 \\
    0 \leq v_2(t) \perp 1-s(t) \geq 0
  \end{cases}
\end{equation}
In order to directly substitute the value of $s(t)$ into the dynamics $\dot x(t) = -s(t)$, a other complementarity formulation can be defined. By setting $\lambda_1(t)=\frac{1}{2}  (1-s(t))$ and $\lambda_2(t) = v_1(t)$, one gets
\begin{equation}
  \label{eq:normalcone2}
  x(t) \in N_{[-1,1]}(s(t))  \Leftrightarrow
  \begin{cases}
      s(t)= 1 -2 \lambda_1(t) \\
      0 \leq \lambda_1(t) \perp x(t)+\lambda_2(t) \geq 0 \\
      0 \leq \lambda_2(t) \perp 1-\lambda_1(t) \geq 0
  \end{cases}
\end{equation}

\section{A class of differential inclusions}
\label{section2}

Let us now introduce the following  class of differential inclusions, where $x(t) \in \RR^{n}$: 

\begin{equation}
\left\{\begin{array}{l}
\dot{x}(t)\in -A(x(t))+ f(t,x(t)),\,\,\,\mbox{a.e. on}\,\,(0,T) \\  \\
x(0)=x_{0}
\end{array}\right.
\label{eq:un}
\end{equation}
The following assumption is made:

\begin{assumption} The following hold:

\begin{itemize}

\item {\bf (i)} $A(\cdot)$ is a multivalued maximal monotone operator from $\RR^{n}$ into $\RR^{n}$, with domain $D(A)$, \textit{i.e.}, for all $x\in D(A), y \in D(A)$ and all $x'\in A(x), y'\in A(y)$, one has
  \begin{equation}
    \label{eq:mono}
    (x'-y')^T(x-y) \geq 0
  \end{equation}
\item {\bf (ii)} There exists $L \geq 0$ such that for all $t \in [0,T]$, for all $x_{1}, x_{2} \in \RR^{n}$, one has $||f(t,x_{1})-f(t,x_{2})|| \leq L ||x_{1}-x_{2}||$.
\item {\bf (iii)} There exists a function $\Phi(\cdot)$ such that for all $R \geq 0$:
$$\Phi(R)=\sup \left\{ \parallel \frac{\partial f}{\partial t}(\cdot,v) \parallel_{{\mathcal L}^{2}((0,T);\RR^{n})}\,\mid\,\parallel v \parallel_{{\mathcal L}^{2}((0,T);\RR^{n})} \leq R\right\} < +\infty$$. 

\end{itemize}
\label{ass1}
\end{assumption}
The following is proved in \cite{bastien2002,bastien2007}. 
\begin{proposition}
Let Assumption \ref{ass1} hold, and let $x_{0} \in D(A)$. Then the differential inclusion (\ref{eq:un}) has a unique solution $x: (0,T) \rightarrow \RR^{n}$ that is Lipschitz continuous. 
\label{proposition1}
\end{proposition}

In this paper we shall focus on inclusions of the form:
\begin{equation}
\left\{\begin{array}{l}
\dot{x}(t)\in f(t,x(t))-B \mbox{Sgn}(Cx(t)+D)
,\,\,\,\mbox{a.e. on}\,\,(0,T) \\  \\
x(0)=x_{0}
\end{array}\right.
\label{eq:deux}
\end{equation}
with $B \in \RR^{n \times m}$,  and $\mbox{Sgn}(Cx(t)+D)\stackrel{\Delta}{=}(\mbox{sgn}(C_{1}x+D_{1}),...,\mbox{sgn}(C_{m}x+D_{m}))^{T} \in \RR^{m}$. It will be shown how to recast (\ref{eq:deux}) into (\ref{eq:un}).
\begin{example}[Equivalent-control-based sliding-mode-control (ECB-SMC)]
  Consider a system $\dot{x}(t)=Fx(t)+Gu$, with an
  equivalent-control-based sliding-mode-control (ECB-SMC) of the form
  $u(x)=-(HG)^{-1}HFx-\alpha (HG)^{-1}\mbox{Sgn}(Hx)$, $\alpha >0$
  (see e.g. \cite{yu2008}). Then the closed-loop system
  $\dot{x}(t)=(F-G(HG)^{-1}HF)x(t)-\alpha G(HG)^{-1}\mbox{Sgn}(Hx(t))$
  fits within (\ref{eq:deux}).
\end{example}

Let us now state a well-posedness result which is a  consequence of Proposition \ref{proposition1}. 
\begin{corollary}
Consider the differential inclusion in (\ref{eq:deux}). Suppose that {\bf (ii)} and {\bf (iii)} of Assumption \ref{ass1}) hold. If there exists an $n \times n$  matrix $P=P^{T}>0$ such that 
\begin{equation}
PB_{\bullet i}=C_{i\bullet}^{T}
\label{eq:cinq}
\end{equation}
for all $1 \leq i \leq m$, then for any initial data the differential inclusion (\ref{eq:deux}) has a unique solution $x: (0,T) \rightarrow \RR^{n}$ that is Lipschitz continuous. 
\label{corollary1}
\end{corollary}

{\bf Proof:} The proof uses a state variable change introduced in \cite{brogliato2004}.  Let $R$ be the symmetric square root of $P$, i.e. $R^{2}=P$. Let us perform the  state transformation $z=Rx$. Then we get
\begin{equation}
\dot{z}(t) \in Rf(t,R^{-1}z(t))-RB \mbox{Sgn}(CR^{-1}z(t)+D)
\end{equation}
Notice that $B \mbox{Sgn}(CR^{-1}z(t)+D)=\sum_{i=1}^{m} B_{\bullet i} \mbox{sgn}(C_{i\bullet}R^{-1}z+D_{i})$. Therefore $RB \mbox{Sgn}(CR^{-1}z(t)+D)=\sum_{i=1}^{m} RB_{\bullet i} \mbox{sgn}(C_{i\bullet}R^{-1}z+D_{i})=\sum_{i=1}^{m} R^{-1}C_{i \bullet}^{T} \mbox{sgn}(C_{\bullet i}R^{-1}z+D_{i})$. We can rewrite the system as
\begin{equation}
\dot{z}(t) \in Rf(t,R^{-1}z(t))-\sum_{i=1}^{m} R^{-1}C_{i\bullet}^{T} \mbox{sgn}(C_{i\bullet}R^{-1}z(t)+D_{i})
\label{eq:six}
\end{equation}
The multivalued mapping $\xi \mapsto \mbox{sgn}(\xi)$ is monotone. By \cite[Exercise 12.4]{wets} it follows that each multivalued mapping $z \mapsto R^{-1}C_{i\bullet}^{T} \mbox{sgn}(C_{i\bullet}R^{-1}z(t)+D_{i})$ is monotone. From \cite[Proposition 1.3.11]{goeleven}  it follows that $R^{-1}C_{i\bullet}^{T} \mbox{sgn}(C_{i\bullet}R^{-1}z(t)+D_{i})=\partial f_{i}(z)$ with $f_{i}(z)=|C_{i\bullet}R^{-1}z(t)+D_{i}|$. By \cite[Theorem 5.7]{rockafellar} it follows that $f_{i}(\cdot)$ is convex. Being the subdifferential of a convex function, the multivalued mapping $z \mapsto \partial f_{i}(z)$ is maximal (monotone) \cite[Corollary 31.5.2]{rockafellar}. Therefore by Proposition \ref{proposition1} the inclusion in  (\ref{eq:six}) possesses a unique Lipschitz solution on $(0,T)$  for any $T >0$ and since $R$ is full--rank so does~(\ref{eq:deux}). \Qend

\begin{example}
Consider the sliding mode system in \cite[Equ.(1)--(4)]{galias2006}. One has  $B=(0\;\;1)^{T}$, $C=(c_{1}\;\;1)$, $D=0$. Then the condition in (\ref{eq:cinq}) holds with $P=\left(\begin{array}{cc} p_{11} & c_{1} \\ c_{1} & 1 \end{array}\right)$ and $p_{11} > (c_{1})^{2}$ assures that $P >0$. 
\label{example1}
\end{example}

\begin{example}
Consider $B=\left(\begin{array}{cc} 1 & 2 \\ 2 & -1 \end{array}\right)$, $\mbox{Sgn}(Cx+D)=(\mbox{sgn}(x_{1}+2x_{2}),\mbox{sgn}(2x_{1}-x_{2}))^{T}$. Trajectories may slide on codimension one surfaces $x_{1}+2x_{2}=0$ or $2x_{1}-x_{2}=0$ and on the codimension 2 surface ($x_{1}+2x_{2}=0$ and $2x_{1}-x_{2}=0$). 
\label{Ex:Multi}
\end{example}

\begin{example}

One solution to reduce chattering is the observer based SMC. Let us consider the following example taken from \cite{young1999}, whose closed-loop dynamics is given by:

\begin{equation}
\left(\begin{array}{c} \dot{x}(t) \\ \dot{e}(t) \\ \dot{x}_{s}(t) \\ \ddot{x}_{s}(t)    \end{array}\right)=\left(\begin{array}{cccc}  0 & 0 & 0 & 0 \\ k & -k & -k & 0 \\ 0 & 0 & 0 & 1 \\ \frac{1}{\tau^{2}} & 0 &  -\frac{1}{\tau^{2}} &  -\frac{2}{\tau}  \end{array}\right) \left(\begin{array}{c} x(t) \\ e(t) \\ x_{s}(t) \\ \dot{x}_{s}(t)   \end{array}\right)-\left(\begin{array}{c} 1 \\ 0 \\ 0 \\ 0   \end{array}\right) \mbox{sgn}(Cx(t))
\label{observer}
\end{equation}
with $C=(1\;-1\;0\;0)$. For the notations see \cite[\S II.C]{young1999}. This system satisfies the condition (\ref{eq:cinq}) with $P=\left(\begin{array}{cccc} 1 & -1 & 0 & 0 \\ -1 & p_{22} & 0 & 0 \\  0 & 0  & p_{33} & 0 \\ 0 & 0 & 0 & p_{44}  \end{array}\right)$, $p_{22} >1$, $p_{33} >0$, $p_{22}>0$.  

\label{Ex:observer}
\end{example}

Notice that the condition (\ref{eq:cinq}) implies that $B_{\bullet i}^T P B_{\bullet i} = B_{\bullet i}^T C_{i \bullet}^T = B_{i\bullet} C_{\bullet i} >0 $. When $m=1$ this is a relative degree one condition.  It is noteworthy that (\ref{eq:cinq}) does not imply that $B$ has full column rank. In particular it does not preclude $m >n$. Dissipative systems with no feedthrough matrix satisfy an input-output constraint similar to (\ref{eq:cinq}) \cite{brogliatoBook}.

\begin{example}[Lyapunov-based discontinuous robust control]  Let us show how the above material adapts to this type of feedback controller. The class of dynamical systems is
\begin{equation}
\dot{x}(t)=f(x(t))+Bu(t)+B\gamma(t),\;\;x(0)=x_{0}
\label{cor1}
\end{equation}
where $x(t) \in \RR^{n}$, $B \in \RR^{n \times m}$, $f(\cdot)$ satisfies assumption \ref{ass1}, and $\gamma(\cdot) \in \RR^{m}$ is a bounded disturbance satisfying $|\gamma_{i}(t)| < \rho_{i}$ for all $1 \leq i \leq m$, all $t \geq 0$ and some finite $\rho_{i}$. The problem is the stabilization of the system at the origin $x=0$, knowing that there exists a function $V(\cdot)$ such that the uncontrolled undisturbed system $\dot{x}(t)=f(x(t))$ admits $V(\cdot)$ as a Lyapunov function. In particular, one has $\dot{V}(x(t)) =\nabla V(x(t))^{T} f(x(t)) \leq 0$ along the trajectories of the free system. Let us rewrite the system in (\ref{cor1}) as
\begin{equation}
\dot{x}(t)=f(x(t))+\sum_{i=1}^{m}B_{\bullet i}u_{i}+\sum_{i=1}^{m}B_{\bullet i}\gamma_{i}(t)
\label{cor2}
\end{equation}
Let us propose the control input $u_{i}(x)=-\rho_{i} \mbox{sgn}(\nabla V^{T}(x)B_{\bullet i})$. We obtain:
\begin{equation}
\dot{x}(t) \in f(x(t))-\sum_{i=1}^{m}\rho_{i} B_{\bullet i}  \mbox{sgn}(\nabla V(x)^{T}B_{\bullet i}) +\sum_{i=1}^{m}B_{\bullet i}\gamma_{i}(t)
\label{cor3}
\end{equation}

We can state the following result.
\begin{corollary}
Suppose that $V(x)=\frac{1}{2}x^{T}Px$, $P=P^{T}> 0$. The system in (\ref{cor3}) has a unique Lipschitz continuous solution on $[0,+\infty)$ for any $x_{0}$.
\label{corollaryCor} 
\end{corollary}

{\bf Proof:} We have $\nabla V(x)^{T}B_{\bullet i}=B_{\bullet,i}^{T}Px$. Let $z=Rx$, where $R >0$ is the symmetric square root of $P$. We may rewrite (\ref{cor3}) as
$$\dot{z}(t) \in Rf(R^{-1}z(t))-\sum_{i=1}^{m}\rho_{i} RB_{\bullet i}  \mbox{sgn}(B_{\bullet i}^{T}Rz) +\sum_{i=1}^{m}RB_{\bullet i}\gamma_{i}(t)$$
Then following the same steps as for the proof of Corollary \ref{corollary1} we conclude that Proposition \ref{proposition1} applies to this system, hence to (\ref{cor3}).   \Qend

Such a controller assures the global asymptotic stability of the equilibrium $x=0$. This is made possible because of the multivalued characteristic of the discontinuous input. The closed-loop system possesses the origin as its unique equilibrium, because of the multivaluedness property. The restriction to quadratic Lyapunov functions stems from monotonicity preserving conditions, and is not straightforwardly avoided. 

\label{example6}
\end{example}

\section{Convergence results and Chattering Free Finite--time Stabilization}
\label{section3}

The differential inclusion (\ref{eq:un}) is time-discretized on $[0,T]$ with a backward Euler scheme as follows:

\begin{equation}
\left\{\begin{array}{l}
\displaystyle\frac{x_{k+1}-x_{k}}{h}+A(x_{k+1}) \ni f(t_{k},x_{k}),\,\,\mbox{for all}\,\,k \in \{0,...,N-1\} \\ \\
x_{0}=x(0)
\end{array}\right.
\label{app2}
\end{equation}
where $h=\frac{T}{N}$. The fully implicit method uses $f(t_{k+1},x_{k+1})$ instead of $f(t_{k},x_{k})$. The convergence and order results stated in Proposition \ref{proposition2} below have been derived for the semi-implicit scheme (\ref{app2}) in \cite{bastien2002}. So the analysis in this section is based on such a discretization. However this is only a particular case of a more general $\theta-$method which is used in practical implementations. The next result is proved in \cite{bastien2002}. 

\begin{proposition}
Under assumption \ref{ass1}, there exists $\eta$ such that for all $h>0$ one has

\begin{equation}
\mbox{For all}\,\,t \in [0,T],\,\,||x(t)-x^{N}(t)|| \leq \eta\,\sqrt{h}
\label{app3}
\end{equation}

Moreover $\lim_{h \rightarrow 0^{+}} \max_{t \in [0,T]} ||x(t)-x^{N}(t)||^{2}+\int_{0}^{t} ||x(s)-x^{N}(s)||^{2}ds=0$.

\label{proposition2}
\end{proposition}

Thus the numerical scheme in (\ref{app2}) has at least order $\frac{1}{2}$, and convergence holds. The conditions of Assumption \ref{ass1} for the existence and uniqueness results of Proposition \ref{proposition1} are sufficient only. Other criteria exist, like Filippov's criterion for uniqueness of solutions \cite[Proposition 5]{cortes}. Similarly it is possible that time-stepping methods converge for systems that satisfy such a criterion, despite no result seems to be available in the literature. As seen in Lemma \ref{lemmayahoo}, the precision of the method may be much better than what is to be expected from (\ref{app3}) on large portions of the trajectories.

The differential inclusion in (\ref{eq:deux}) is therefore discretized as follows:
\begin{equation}
\left\{\begin{array}{l}
\displaystyle\frac{x_{k+1}-x_{k}}{h} \in f(t_{k},x_{k})-B \mbox{Sgn}(Cx_{k+1}+D)
,\,\,\,\mbox{a.e. on}\,\,(0,T) \\  \\
x(0)=x_{0}
\end{array}\right.
\label{eq:discdeux}
\end{equation}
One sees that advancing the implicit method from step $k$ to step $k+1$ involves solving generalized equations with unknown $x_{k+1}$, of the form $0 \in F_{s}(x_{k+1})+F_{m}(x_{k+1})$ where $F_{s}(\cdot)$ is singlevalued while $F_{m}(\cdot)$ is multivalued. $h$, $t_{k}$ and $x_{k}$ appear as parameters of the generalized equations. Solving such generalized equations thus boils down to computing the intersection between the graph of $F_{s}(\cdot)$ and  the graph of $F_{m}(\cdot)$ as illustrated in section \ref{section1}. The result of Proposition \ref{proposition2} applies to (\ref{eq:discdeux}). As we shall see next, such an implicit method also assures a good estimate of the derivative $\dot{x}$ and a smooth stabilization of the discrete-time solution on the sliding surface.  

 Before stating the smooth stabilization result, let us consider a  preliminary result. Let us denote the output of the dynamical as:
 \begin{equation}
   \label{eq:outout}
   y(t) \stackrel{\Delta}{=} Cx(t) + D
 \end{equation}
\begin{lemma}
\label{Lemma:bound1}
  Let us assume that a sliding mode exists for some indices $ i \in \alpha \subset \{1\ldots m\}$ such that
  \begin{equation}
    \label{eq:AssSlide}
    \exists t_* >0, \quad y_\alpha(t) = C_{\alpha\bullet} x(t) +D_\alpha = 0, \quad \text{ for all } t > t_*.
  \end{equation}
Then there exists $\rho>0$ such that for all $t > t_*$ and  for all  $x(t)$  such that $ C_{\alpha\bullet} x(t) +D_\alpha = 0$, one has
\begin{equation}
  \label{eq:AssSlide2-1}
 \| (C_{\alpha\bullet} f( x(t),t)) \| \leq \rho 
\end{equation}
Furthermore,  let Assumption 1.(ii) holds, then the following bound is satisfied in the neighborhood of the sliding subspace,
\begin{equation}
  \label{eq:AssSlide2}
\exists r > 0, \exists \kappa >0, \exists \rho > 0 \text{ such that } \forall t > t_* , \forall \bar x \in B_r(x),   \| (C_{\alpha\bullet} f(\bar x,t)) \| \leq \kappa r + \rho
\end{equation}
for all $x(t)$  such that $C_{\alpha\bullet} x(t) +D_\alpha = 0$.\end{lemma}

{\bf Proof:} From~(\ref{eq:AssSlide}), we have $\dot y_{\alpha}(t) \in C_{\alpha\bullet} f(x(t),t)  -h C_{\alpha\bullet}B\mbox{Sgn}(y(t)) $.  For $t > t_*$,  the sliding mode $y_\alpha(t) = 0$ implies that $\dot y_\alpha(t) = C_{\alpha\bullet} \dot x(t) = 0  $  for all $t > t_*$ and therefore
\begin{equation}
 C_{\alpha\bullet} f(x(t),t) \in  C_{\alpha\bullet}B\mbox{Sgn(y(t))} \label{eq:AssSlide3}
\end{equation}
The inclusion (\ref{eq:AssSlide3}) yields 
\begin{equation}
  \label{eq:AssSlide4bis}
 \exists \rho > 0 , \| (C_{\alpha\bullet} f(x(t),t)) \| \leq  \rho
\end{equation}
for all $x(t)$  such that $C_{\alpha\bullet} x(t) +D_\alpha = 0$.
By the assumption 1.(ii), the Lipschitz continuity of $f(\cdot,\cdot)$ allows us to write for some $\kappa >0$
\begin{equation}
  \label{eq:AssSlide4}
  \forall \bar x(t) \in B_r(x(t)),\quad  \|   C_{\alpha\bullet} (f(\bar x(t),t)- f(x(t),t)) \| \leq \|C_{\alpha\bullet} \| L r \stackrel{\Delta}{=} \kappa r.
\end{equation}
Combining (\ref{eq:AssSlide4bis}) and (\ref{eq:AssSlide4}) ends the proof. \Qend

Lemma \ref{lemmayahoo} extends to (\ref{eq:discdeux}) as follows when the sliding surface of codimension $| \alpha |$ is attained.

\begin{lemma}
 Let us assume that a sliding mode occurs for the index $\alpha \subset \{1\ldots m\}$, that is $y_\alpha(t) = 0, t > t_*$. Let $C$ and $B$ be such that (\ref{eq:cinq}) holds and $C_{\alpha\bullet}B_{\bullet\alpha}>0$. Then there exists $h_c >0$ such that $\forall h<h_c$, there  exists $k_0 \in \NN$ such that $ y_{k_0+n} =  C x_{k_0+n+1}+D=0$ for all integers $n \geq 1$.  
\label{lemmaoups}
\end{lemma}

{\bf Proof:} At each time--step, we have to solve for $y_{k+1} =C x_{k+1}+D$ and $s_{k+1}$ the generalized equation
\begin{equation}
\left\{\begin{array}{l}
y_{k+1} = y_k + hCf(t_k,x_k) -hCB s_{k+1} \\[2mm]
s_{k+1}\in  \mbox{Sgn}(y_{k+1})
\end{array}\right.
\label{eq:yy}
\end{equation}
Under condition (\ref{eq:cinq}), the convergence of the time--stepping scheme is ensured by Proposition~\ref{proposition2}. The convergence and  the existence of the sliding mode ensure that
\begin{equation}
  \label{eq:AssSlide5}
  \exists k_0, \exists K_1 >0, \exists K_2 >0, \exists t_1 >t_* \text{ such that}\quad \|y_{\alpha,k_0}\| \leq K_1 \sqrt h \text{ and }    \|x_{k_0} - x(t_1)\| \leq K_2 \sqrt h 
\end{equation}
for $ C_{\alpha\bullet} x(t_1) +D_{\alpha} = 0$.
Using (\ref{eq:AssSlide2}) for $x(t_1)$ and a sufficiently small $h$ such that $r = K_2 \sqrt h  $, we have the following bound
\begin{equation}
  \label{eq:AssSlide6}
  \| y_{\alpha,k_0} + hC_{\alpha,\bullet}f(t_{k_0},x_{k_0}) \| \leq  \sqrt h( K_1  +  h \kappa K_2    + \sqrt h  \rho)
\end{equation}
Introducing the complementary index set $\beta  = \{ i, y_i(t) = C_{i\bullet} x(t) +D_i \neq  0 \}$, for $ t >t_*$ almost everywhere  and using (\ref{eq:AssSlide6}) we obtain that there exists $\rho_1>0$ such that 
\begin{equation}
  \label{eq:AssSlide8}
 \| y_{\alpha,k_0} + hC_{\alpha,\bullet}f(t_{k_0},x_{k_0}) - h C_{\alpha\bullet}B_{\bullet\beta} \mbox{Sgn}(y_{\beta,k_0+1})  \| \leq  \sqrt h( K_1  +  h \kappa K_2    + \sqrt h  (\rho+\rho_1))
\end{equation}
and therefore it is possible to choose $h_1$ such that for all $h <h_1$
\begin{equation}
  \label{eq:AssSlide7}
   \left|\left[-h(C_{\alpha\bullet}B_{\bullet\alpha})^{-1}\left[y_{\alpha,k_0} + hC_{\alpha,\bullet}f(t_{k_0},x_{k_0}) - h C_{\alpha\bullet}B_{\bullet\beta} \mbox{Sgn}(y_{\beta,k_0+1})  \right] \right]_i \right| \leq 1,\text{ for all } i \in \alpha.
\end{equation}
If~(\ref{eq:AssSlide7}) is satisfied, the unique solution of  (\ref{eq:yy}) at the iteration $k_0+1$ is given by
\begin{equation}
  \label{eq:AssSlide10}
  y_{\alpha,k_0+1} = 0 ; s_{\alpha,k_0+1} = -h(C_{\alpha\bullet}B_{\bullet\alpha})^{-1}\left[y_{\alpha,k_0} + hC_{\alpha,\bullet}f(t_{k_0},x_{k_0}) - h C_{\alpha\bullet}B_{\bullet\beta}  \mbox{Sgn}(y_{\beta,k_0+1})\right]
\end{equation}
The next iterate will by given by the solution of the generalized equation,
\begin{equation}
  \label{eq:AssSlide11}
  \left\{\begin{array}{l}
      y_{k_0+2} = hCf(t_{k_0+1},x_{k_0+1}) -hCB s_{k_0+2} \\[2mm]
      s_{k_0+2}\in  \mbox{Sgn}(y_{k_0+2})
\end{array}\right. .
\end{equation}
Using the fact that $ y_{\alpha,k_0+1} = C_{\alpha\bullet} x_{k_0+1} +D_{\alpha}=0 $, we can use~(\ref{eq:AssSlide4bis}) to conclude that there exists $h_2$ such that for  all $h <h_2$
\begin{equation}
  \label{eq:AssSlide7bis}
   \left|\left[-h(C_{\alpha\bullet}B_{\bullet\alpha})^{-1}\left[ hC_{\alpha,\bullet}f(t_{k_0+1},x_{k_0+1}) - h C_{\alpha\bullet}B_{\bullet\beta} \mbox{Sgn}(y_{\beta,k_0+2})  \right] \right]_i \right| \leq 1,\text{ for all } i \in \alpha, 
\end{equation}
and therefore the solution of (\ref{eq:AssSlide11}) is
\begin{equation}
  \label{eq:AssSlide10bis}
  y_{\alpha,k_0+2} = 0 ; s_{\alpha,k_0+2} = -h(C_{\alpha\bullet}B_{\bullet\alpha})^{-1}\left[hC_{\alpha,\bullet}f(t_{k_0+1},x_{k_0+1}) - h C_{\alpha\bullet}B_{\bullet\beta}  \mbox{Sgn}(y_{\beta,k_0+2})  \right]
\end{equation}
The bound (\ref{eq:AssSlide4bis}) is uniform and can be applied for the next steps. Choosing $h_c$ as the minimum of the considered time steps $h_1, h_2, \ldots$, the proof is obtained for $y_{\alpha,k_0+n}, n \geq 1$.
\Qend

The finite-time convergence of the time-discretization of similar nonsmooth dynamical systems (essentially mechanical systems with dry friction) is proved in \cite{Baji2006}. Our results may therefore be considered as the continuation of studies on the finite-time convergence for algorithms of the proximal type.

\section{Discrete--time Sliding Mode Control (SMC)}
\label{sec:DTSMC}
This section is devoted to show how the above discretizations may be used in a digital control framework.

\subsection{Example of an implicit Euler controller (IEC)}

Let us come back to the inclusion in (\ref{eq:sept}). For this simple system, the ZOH and the Euler discretization yield the same--discrete system. Assume the integrator $\dot{x}(t)=u(t)$ is sampled with sampling period $h >0$. On the time interval $[t_{k},t_{k+1})$ one has $x(t)=x_{k}+h_{t}u_{k}$, where $h_{t}=t-t_{k}$. The controller $u(x)=-\mbox{sgn}(x)$ is known as the equivalent control-based SMC \cite{yu2008}. Let us implement a ``backward'' controller $u_{k}=-\mbox{sgn}(x_{k+1})$ at time $t_{k}$, following the above lines.  Suppose that $x_{k} \in [-h,h]$. Then following the same calculations as in the proof of Lemma \ref{lemmayahoo}, we obtain that $s_{k+1} = \frac{x_k}{h}$. Therefore on $[t_{k},t_{k+1})$:
\begin{equation}
x(t)=x_{k}-\frac{h_{t}}{h}x_{k}
\end{equation}
and it follows that $x(t_{k+1})=x_{k+1}=0$. On the next sampling interval $[t_{k+1},t_{k+2})$ one obtains $s_{k+2}=0$ 
\begin{equation}
x(t)=x_{k+1}-\frac{h_{t}}{h}x_{k+1}=0-\frac{h_{t}}{h}0=0
\end{equation}
and so on on the next intervals, where the zero value is obviously some small value at the machine accuracy. if we suppose that  $x_{k} \notin [-h,h]$, the value of $s_{k+1}$ is $1$ or $-1$ according to the sign of $x_k$. To summarize the control is given explicitly in terms of $x_k$ and $h$ by
\begin{equation}
  \label{eq:proj1}
  u_k= - \mbox{proj}_{[-1,1]} (\frac{x_k}{h})
\end{equation}
where $\mbox{proj}_C$ denotes the Euclidean projection operator onto the set $C$.

 As alluded to above, such an ``implicit'' input is causal and can be computed at $t_{k}$ with the values of the state at $t_{k}$ by~(\ref{eq:proj1}). It requires at each step to solve a rather simple multivalued problem which a Mixed Linear Complementarity Problem (MLCP, see Section \ref{Sec:TimeDiscretization} below). It is not of the high gain type. 

\begin{remark}
The fact that the function sgn$(\cdot)$ generates only binary values ($+1$ or $-1$) does not hamper the above method to work. Indeed the implicit Euler method allows us to compute values of the sign multifunction inside its multivalued part at $x_{k}=0$. 
\end{remark}



\subsection{Extension to ZOH discretized systems}
\label{subsecZOH}

The ZOH discretization of linear time invariant systems $\dot{x}(t)=Fx(t)+Gu(t)$ with an ECB-SMC controller, $u(x) = - (CG)^{-1}(C Fx + \alpha \mbox{Sgn}(C x)), \alpha >0 $ results in a discrete-time system of the form:
\begin{equation}
x_{k+1}=\Phi x_{k} - \Gamma s_{k}\;\;\mbox{for all}\;\;t \in [kh,(k+1)h)
\label{ZOH1}
\end{equation}
where $h >0$ is the sampling period, and
\begin{eqnarray}
  \label{eq:ZOHGAMMA}
  \Phi&=&\exp(Fh)-\int_{0}^{h}\exp(F\tau)d\tau G   (CG)^{-1}   C F \\
  \Gamma&=&\int_{0}^{h} \exp(F\tau) G (CG)^{-1}d\tau
\end{eqnarray}
 with $G \in \RR^{n \times m}$, $C \in \RR^{m \times n}$, when a explicit Euler implementation of the control is performed \cite{wang,yu2003}. For an implicit Euler implementation, let us set
\begin{equation}
  \left\{\begin{array}{l}
      u_k =  - (CG)^{-1}(C Fx_k  + s_{k+1})\\[2mm]
      s_{k+1}= \; \mbox{Sgn}(C x_{k+1}),
\end{array}\right.
\label{ZOH2}
\end{equation}
which corresponds to the implicit discrete time version of the ECB-SMC controller. We therefore get on each sampling period:
\begin{equation}
x_{k+1}=\Phi x_{k}- \Gamma s_{k+1} \;\;\mbox{for all}\;\;t \in [kh,(k+1)h)
\label{ZOH3}
\end{equation}
At each time--step, one has to solve
\begin{equation}
\left\{\begin{array}{l}
x_{k+1} = \Phi x_k  -\Gamma s_{k+1} \\[2mm]
y_{k+1} = C x_{k+1} + D \\[2mm]
s_{k+1} \in \mbox{Sgn} (y_{k+1}) 
\end{array}\right.
\label{eq:ZOHZOH1}.
\end{equation}
Inserting the first line of (\ref{eq:ZOHZOH1}) into the second line we obtain the following one--step system 
\begin{equation}
\left\{\begin{array}{l}
y_{k+1} = C \Phi x_k   + D  -  C \Gamma s_{k+1}    \\[2mm]
s_{k+1} \in \mbox{Sgn} (y_{k+1}) 
\end{array}\right.
\label{eq:ZOHZOH2}.
\end{equation}
Comparing with the time--discretized systems in (\ref{eq:discdeux}) and (\ref{eq:yy})  one sees that the term $h CB$ is replaced in case of a ZOH method by the term $ C \Gamma$. Provided the problem has a unique solution one can compute the controller in (\ref{ZOH2}) with the knowledge of $x_{k}$, $h$, $F$, $G$ and $C$. We will see in the next Section how the computation can be carried out in  practice.


\section{Implementation of  Discrete--Time Systems}
\label{Sec:TimeDiscretization}
Let us consider in this section the following discrete--time system:
\begin{equation}
\left\{\begin{array}{l}
x_{k+1} = R x_k + p - S s_{k+1} \\[2mm]
y_{k+1} = C x_{k+1} + D \\[2mm]
s_{k+1} \in \mbox{Sgn} (y_{k+1}) 
\end{array}\right.
\label{eq:TDS1}
\end{equation}
where $k\geq 0$ is an integer, $x_k$ the discrete state, $y_k$ the discrete output and $s_{k}$ the discrete input. The discrete system~(\ref{eq:TDS1}) is a common representative for the discretization given by (\ref{eq:discdeux}), (\ref{app2}) or  (\ref{eq:ZOHZOH1}) and  the matrices $R\in \RR^{n\times n}$, $S\in \RR^{n\times m}$ and the vector $p\in \RR^n$ are determined by the chosen time--discretization method and detailed in Section~\ref{Sec:TimeDiscretizationmethods}. The matrices $C$ and $D$ are given by their definition in~(\ref{eq:deux}). 

\subsection{Mixed Linear Complementarity Problem (MLCP)}
The time--discretized system~(\ref{eq:TDS1})  appears to be a Mixed Linear Complementarity Problem (MLCP) that we have to solve at each time--step. Let us define what is a MLCP in its general form with bounds constraints as it has been proposed in \cite{Dirkse.Ferris1995}:

\begin{definition}[MLCP]
  Given a matrix $M \in \RR^{m \times m}$, a vector $q\in  \RR^{m}$ and lower and upper bounds $l,u\in \overline\RR^{m}$, find $z\in \RR^m$, $w,v \in \RR^m_+$ such that
  \begin{equation}
    \label{eq:MLCP1}
    \left\{\begin{array}[]{c}
      M z + q  = w-v \\[2mm]
      l \leq z \leq u \\[2mm]
      (z-l)^Tw=0 \\[2mm]
      (u-z)^Tv = 0
    \end{array}\right.
  \end{equation}
where  $\overline\RR = \RR \cup \{+\infty,-\infty\}$. 
\end{definition}
Note that the problem~(\ref{eq:MLCP1})  implies  that
\begin{equation}
  \label{eq:MLCP10}
  - (M z +q) \in N_{[l,u]}(z).
\end{equation}
where the notation $N_{C}(x)$ is used for the normal cone in the Convex Analysis sense to a convex set $C$ at the point $x$. The box $[l,u]\subset \RR^m$ is defined by the Cartesian product of the intervals $[l_i,u_i], i \in \{1,\ldots,m\}$. The normal cone to a convex set  is a standard instance of a multi--valued mapping~\cite{rockafellar}. The relation (\ref{eq:MLCP10}) is equivalent to the MCP~(\ref{eq:MLCP1}) if we assume that  $w$ is the positive part of $Mz +q$, that is $w = (Mz+q)^+=max(0, M z+q))$ and $v$ is the negative part of $M z+ q $, that is $v = (Mz+q)^-= max(0, -(Mz+q))$.

In order to state the problem~(\ref{eq:TDS1}) as a MLCP, the variable $x_{k+1}$ is condensed into the second line such that
\begin{equation}
  \label{eq:TDS2}
  \left\{\begin{array}{l}
  y_{k+1} = C R x_k + C p - C S s_{k+1}  + D \\[2mm]
  s_{k+1} \in \mbox{Sgn} (y_{k+1})
\end{array}\right.
\end{equation}
and the following variable and parameters are defined as follows
\begin{equation}
   \label{eq:TDS3}
  \left\{ \begin{array}{l}
   z  =s_{k+1} ; \quad  y_{k+1} = w- v  \\ [2mm]
   M = CS, \quad q =  -(C R x_k + C p + D) \\[2mm]
   l_i = -1, u_i = 1, i = 1\ldots m. 
 \end{array}\right.
\end{equation}
Finally, the problem~(\ref{eq:TDS1}) can be recast into a MLCP  by observing that
\begin{equation}
  \label{eq:MLCP2}
  \begin{array}{c}
    s_{k+1} \in \mbox{Sgn}{( y_{k+1} )}\\
\Updownarrow \\   y_{k+1} \in N_{[-1,1]^m}(  s_{k+1}    ) \\
\Updownarrow \\
s_{k+1}\in [-1,1]^m \text{ and }
   \begin{cases}
     y_{j,k+1} = 0 \text{ if } s_{j,k+1} \in ]-1,1[ \\
     y_{j,k+1} \leq 0 \text{ if } s_{j,k+1} = - 1 \\
     y_{j,k+1} \geq 0 \text{ if } s_{j,k+1} =  1 \\  
   \end{cases}, j\in \{1,\ldots, m\}
 \end{array}
\end{equation}

The MLCP~(\ref{eq:MLCP1}) is a well-known problem in the mathematical programming theory arising for instance from the Karush--Kuhn--Tucker optimality conditions of a quadratic program or from the primal/dual optimality conditions of a linear program. The MCLP enjoys a large number of numerical algorithms and  several reliable solvers have been implemented. Several families of solvers may be cited: a) extensions of Lemke and principal pivotal techniques for LCP to MLCP \cite{Sargent1978,Rutherford1993,Dirkse.Ferris1995,Cao.Ferris1996}  b) extensions of projection/splitting techniques for MLCP \cite{Facchinei,Cottle.Pang.ea1992} and  c) semi--smooth Newton methods \cite{Munson.ea2001}. In this paper, the computations are done with the help the \textsc{Siconos/Numerics} open source Library \cite{Acary.Perignon2007} and/or the PATH solver \cite{Dirkse.Ferris1995}. The results of existence and uniqueness of solutions of (\ref{eq:MLCP1}) are related to the properties of $M$ (P-properties or coherent orientations of the associated affine map (normal map) for particular cases of bounds constraints). Without entering into further details, we refer to \cite{Harker.Pang1990,Facchinei} for the main results. The assumptions on the matrix $M$ drives the choice of particular solvers that can be in polynomial time rather than standard exponential time for brute force enumerative solvers.

\subsection{Some Time--Discretization Methods}
\label{Sec:TimeDiscretizationmethods}

In this Section, the formulation of the discrete--time system~(\ref{eq:TDS1}) is related to the continuous time system~(\ref{eq:deux}) through a given discretization method. 
\paragraph{Explicit Euler discretization of $f(\cdot,\cdot)$} Let us start with the explicit Euler discretization method of the term $f(t,x(t))$ as it has been given in (\ref{app2}). At each time step, the matrices in (\ref{eq:TDS1}) and in the MLCP~(\ref{eq:MLCP1}) can be identified as
\begin{equation}
  \label{eq:TDS3-a}
  R = I , p = h f(t_k,x_k), S = hB, \quad M = hCB, q =  -  (h C f(t_k,x_k)  + C x_k +  D) 
\end{equation}
Let the assumptions of Corollary~\ref{corollary1} be satisfied with $B$ full--column rank ($CB= B^T PB >0$). This result ensures the existence and uniqueness of a solution of the MCLP. Furthermore, standard pivotal techniques such as Lemke's method or projection/splitting such as Projected Successive Over-Relation (PSOR) compute the solution.

\paragraph{Implicit Euler and $\theta$- method} In a more general way, we can choose to time--discretize the term $f(t,x(t))$ by a implicit Euler scheme or a $\theta$-method. The main motivation for doing in this way is the higher accuracy and stability that we can obtain for such a numerical integration scheme (see \cite{Acary.BrogliatoRR2009} for an example of instability with the Explicit Euler method). Let us consider first that the mapping $f(\cdot,\cdot)$ is affine, that is $f(t,x(t)) = F x(t) +g$. The matrices in (\ref{eq:TDS1}) and in the MLCP~(\ref{eq:MLCP1}) can be identified as
\begin{equation}
  \label{eq:TDS4}
 \left\{ \begin{array}{l}
  R = (I-h\theta F)^{-1}(I+ h (1-\theta)F) , p = (I-h\theta F)^{-1} g , S = h(I-h\theta F)^{-1} B, \\[2mm]
 M = h C(I-h\theta F)^{-1}B, q =  -( (I-h\theta F)^{-1}(I+ h (1-\theta)F)x_k + (I-h\theta F)^{-1} g +D)
\end{array}\right.
\end{equation}
for $\theta \in [0,1]$.  For $\theta= 0$, the explicit Euler case is retrieved. For $\theta =1$, the implicit Euler scheme is used to discretize $f$.
If the mapping $f(\cdot,\cdot)$ is nonlinear, a newton linearization can be invoked. In this case, the solution at each time step is sought as a limit of solutions of successive MLCPs. We refer to \cite{Acary.BrogliatoRR2009} for a detailed presentation of these developments.

\begin{figure}
  \centering
\begin{picture}(0,0)%
\epsfig{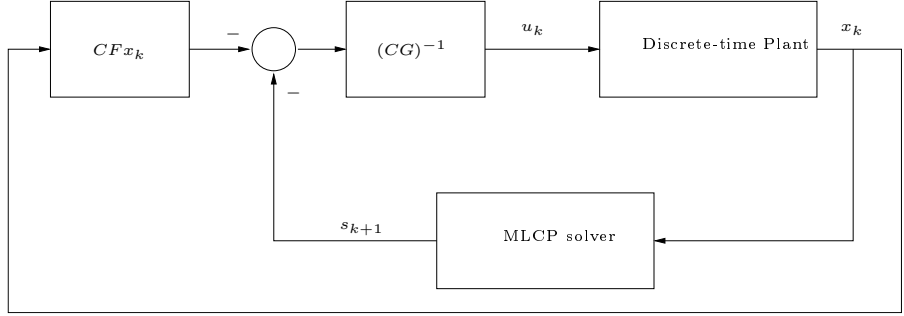}%
\end{picture}%
\setlength{\unitlength}{1973sp}%
\begingroup\makeatletter\ifx\SetFigFont\undefined%
\gdef\SetFigFont#1#2#3#4#5{%
  \reset@font\fontsize{#1}{#2pt}%
  \fontfamily{#3}\fontseries{#4}\fontshape{#5}%
  \selectfont}%
\fi\endgroup%
\begin{picture}(11124,3924)(1939,-5473)
\put(8101,-4561){\makebox(0,0)[lb]{\smash{{\SetFigFont{6}{7.2}{\familydefault}{\mddefault}{\updefault}{\color[rgb]{0,0,0}MLCP solver}%
}}}}
\put(5401,-2761){\makebox(0,0)[lb]{\smash{{\SetFigFont{6}{7.2}{\rmdefault}{\mddefault}{\updefault}{\color[rgb]{0,0,0}$-$}%
}}}}
\put(4651,-2011){\makebox(0,0)[lb]{\smash{{\SetFigFont{6}{7.2}{\rmdefault}{\mddefault}{\updefault}{\color[rgb]{0,0,0}$-$}%
}}}}
\put(3001,-2236){\makebox(0,0)[lb]{\smash{{\SetFigFont{6}{7.2}{\rmdefault}{\mddefault}{\updefault}{\color[rgb]{0,0,0}$CFx_k$}%
}}}}
\put(6526,-2236){\makebox(0,0)[lb]{\smash{{\SetFigFont{6}{7.2}{\rmdefault}{\mddefault}{\updefault}{\color[rgb]{0,0,0}$(CG)^{-1}$}%
}}}}
\put(8326,-1936){\makebox(0,0)[lb]{\smash{{\SetFigFont{6}{7.2}{\rmdefault}{\mddefault}{\updefault}{\color[rgb]{0,0,0}$u_k$}%
}}}}
\put(12301,-1936){\makebox(0,0)[lb]{\smash{{\SetFigFont{6}{7.2}{\rmdefault}{\mddefault}{\updefault}{\color[rgb]{0,0,0}$x_k$}%
}}}}
\put(6076,-4411){\makebox(0,0)[lb]{\smash{{\SetFigFont{6}{7.2}{\rmdefault}{\mddefault}{\updefault}{\color[rgb]{0,0,0}$s_{k+1}$}%
}}}}
\put(9826,-2161){\makebox(0,0)[lb]{\smash{{\SetFigFont{6}{7.2}{\rmdefault}{\mddefault}{\updefault}{\color[rgb]{0,0,0}Discrete-time Plant}%
}}}}
\end{picture}%
 \caption{Control system schema with implicit Euler implementation.}
\label{fig:ControlSchema}
\end{figure}

\paragraph{Zero--Order Holder (ZOH) method} The ZOH discretization presented in Section~\ref{subsecZOH} can be also formalized into the form (\ref{eq:TDS1}) and then (\ref{eq:MLCP1}) with 
\begin{equation}
  \label{eq:TDS4-a}
 \left. \begin{array}{l}
  R = \Phi , p = 0 , S = \Gamma, \quad M =   C\Gamma, q =  -(C \Phi x_k ) 
\end{array}\right.
\end{equation}
In practice, numerous methods are available to compute the ZOH discretization, \textit{i.e.}, $\Phi$ and $\Gamma$ which amounts to compute the matrix exponential and its time integral \cite{Moler.VanLoan2003}. In this work, the numerical computation is performed using an explicit Runge--Kutta method with high order of accuracy and a numerical tolerance near the machine precision threshold. On the Figure~\ref{fig:ControlSchema}, the control scheme is depicted showing that the controller is causal and computed form $x_k$.

\section{Two other classes of differential inclusions}
In this section, we introduce other classes of differential inclusions which extend (\ref{eq:deux}). The second class of differential inclusions is:
\begin{equation}
\left\{\begin{array}{l}
\dot{x}(t)\in f(t,x(t))-\sum_{i=1}^{m}(A_{i}x(t)+B_{i}) \mbox{sgn}(C_{i}x(t)+D_{i})
,\,\,\,\mbox{a.e. on}\,\,(0,T) \\  \\
x(0)=x_{0}
\end{array}\right.
\label{eq:trois}
\end{equation}
with $A_{i} \in \RR^{n \times n}$, $B_{i} \in \RR^{n \times 1}$, $C_{i} \in \RR^{1 \times n}$, $D_{i}$ is a scalar.

The third class that we shall analyze is:
\begin{equation}
\left\{\begin{array}{l}
\dot{x}(t)\in f(t,x(t))-g(x(t)) \mbox{Sgn}(h(x(t))
,\,\,\,\mbox{a.e. on}\,\,(0,T) \\  \\
x(0)=x_{0}
\end{array}\right.
\label{eq:quatre}
\end{equation}
where $g: \RR^{n} \rightarrow \RR^{n \times m}$ and $h: \RR^{n} \rightarrow \RR^{m}$ are smooth functions, $\mbox{Sgn}(h(x))=[\mbox{sgn}(h_{1}(x),...,\mbox{sgn}(h_{m}(x)]^{T} \in \RR^{m}$.

\begin{corollary}
Consider the differential inclusion in (\ref{eq:trois}). Suppose that {\bf (ii)} and {\bf (iii)} of Assumption \ref{ass1}) hold. Suppose that the multivalued mappings $x \mapsto (A_{i}x+B_{i}) \mbox{sgn}(C_{i}x+D_{i})$, $1 \leq i \leq m$, are hypomonotone. Then for any initial data  the differential inclusion (\ref{eq:trois}) has a unique solution $x: (0,T) \rightarrow \RR^{n}$ that is Lipschitz continuous. 
\label{corollary2}
\end{corollary}

The proof is straightforward and is omitted.

\begin{example}
The mapping $\RR \rightarrow \RR$, $x \mapsto (x+1) \mbox{sgn}(x)$ is hypomonotone. Indeed $(x+1) \mbox{sgn}(x)=|x|+\mbox{sgn}(x)$ and $x \mapsto |x|+k x$ is monotone for any $k \geq 1$. 
\label{example2}
\end{example}

\begin{example}
Let $k_{1}$, $k_{2}$ be reals. The mapping $F: \RR \rightarrow \RR$, $x \mapsto \left\{\begin{array}{lll} -k_{1}x + 1  & \mbox{if} & x \geq 0 \\ -k_{2}x-1 & \mbox{if} & x \leq 0 \\ \mbox{$[-1,1]$} & \mbox{if} & x=0 \end{array} \right.$, is hypomonotone with constant $k$ for any $k \geq \max(|k_{1}|,|k_{2}|)$. Let $k_{1}=-k_{2}=k$. Then $F(x)=(kx+1)\mbox{sgn}(x)$. The linearized Stribeck friction model (with multivalued part at zero tangential velocity) \cite{marton} is hypomonotone. 
\label{example3}
\end{example}

Let us state other cases where (\ref{eq:trois}) fits within Proposition \ref{proposition1}. 
\begin{lemma}
Let $B_{i}=\alpha C_{i}^{T}$ for some $\alpha >0$, $D_{i}=0$ and Ker$(C_{i}) \subseteq$ Ker$(A_{i})$. Then for any initial data  the differential inclusion (\ref{eq:trois}) has a unique solution $x: (0,T) \rightarrow \RR^{n}$ that is Lipschitz continuous.. 
\end{lemma}

{\bf Proof:} First notice that the set-valued mapping $x \mapsto \alpha C_{i}^{T}  \mbox{sgn}(C_{i}x+D_{i})$ is maximal monotone \cite[Exercise 12.4]{wets}.   Under the lemma's conditions, one sees that $x \mapsto A_{i}x \mbox{sgn}(C_{i}x)$  is continuous on the surface $\Sigma_{i}=\{x \in \RR^{n} \mid  C_{i}x=0 \}$. Indeed the jump of the vector field is equal to $2A_{i}x=0$ on $\Sigma_{i}$. Moreover it is Lipschitz continuous as it is piececewise linear. Hence Proposition \ref{proposition1} applies.  \Qend

\begin{corollary}
Consider the differential inclusion in (\ref{eq:quatre}). Suppose that {\bf (ii)} and {\bf (iii)} of Assumption \ref{ass1}) hold. Suppose that the multivalued mappings $x \mapsto g_{\bullet i}(x) \mbox{sgn}(h_{i}(x))$, $1 \leq i \leq m$, are hypomonotone. Then for any initial data  the differential inclusion (\ref{eq:quatre}) has a unique solution $x: (0,T) \rightarrow \RR^{n}$ that is Lipschitz continuous. 
\label{corollary3}
\end{corollary}

The proof is straightforward and is omitted.

\begin{example}
The mapping $F: \RR \rightarrow \RR$, $x \mapsto \frac{1}{1+x^{2}} \mbox{sgn}(\arctan(x))$, is hypomonotone with any $k \geq \frac{9}{8\sqrt{3}}$. 
\label{example4}
\end{example}

\begin{remark}
As noted in \cite{young1999} chattering may be due in sliding mode control appplications to the presence of parasitic dynamics. Simple modelling of these yield when inserted in (\ref{eq:sept}) the differential inclusion (see \cite[(7) (8)]{young1999})
\begin{equation}
\left(\begin{array}{c} \dot{x}(t) \\ \dot{x}_{s}(t) \\ \ddot{x}_{s}(t)   \end{array}\right)=\left(\begin{array}{ccc} 0 & 0 & 0 \\ 0 & 0 & 1 \\ \frac{1}{\tau^{2}} & -\frac{1}{\tau^{2}} & -\frac{2}{\tau}  \end{array}\right) \left(\begin{array}{c} x(t) \\ x_{s}(t) \\ \dot{x}_{s}(t) \end{array}\right)-\left(\begin{array}{c} 1 \\ 0 \\ 0  \end{array}\right)\mbox{sgn}(Cx(t))
\label{parasitics}
\end{equation}
with $C=(0\;1\;0)$. The relative degree of the triplet $(A,B,C)$ of this system is $r=3$, where $B=(1\;0\;0)^{T}$ and $A=\left(\begin{array}{ccc} 0 & 0 & 0 \\ 0 & 0 & 1 \\ \frac{1}{\tau^{2}} & -\frac{1}{\tau^{2}} & -\frac{2}{\tau}  \end{array}\right)$. This system does not fit within the above classes of inclusions. Similar conclusions hold for the other form of parasitics in \cite[(3) (4)]{young1999}. Such parasitics may be seen as a {\em non collocation} issue, that is known to greatly influence the stability of systems and usually may yield instability. The mere existence and uniqueness of solutions for such relative degree 3 systems is not trivial. In \cite{pogromski} an example is given that possesses an infinity of absolutely continuous Filippov's solutions, but a unique so-called forward solution. One interesting question is to determine what kind of solution is approximated by the backward Euler method applied to (\ref{parasitics}) which, according to \cite[Theorem 1]{pogromski} has a unique forward solution  since its leading Markov parameter is $CA^{2}B=\frac{1}{\tau^{2}} >0$. A possible solution for this non collocation issue is the observer design of Example~\ref{Ex:observer}.

\end{remark}

The differential inclusions in (\ref{eq:trois})--(\ref{eq:quatre}) are therefore discretized as follows:
\begin{equation}
\left\{\begin{array}{l}
\frac{x_{k+1}-x_{k}}{h} \in f(t_{k},x_{k})-\rho x_{k}-\sum_{i=1}^{m}(A_{i}x_{k+1}+B_{i}) \mbox{sgn}(C_{i}x_{k+1}+D_{i})+\rho_{i}x_{k+1}
,\,\,\,\mbox{a.e. on}\,\,(0,T) \\  \\
x(0)=x_{0}
\end{array}\right.
\label{eq:disctrois}
\end{equation}
and
\begin{equation}
\left\{\begin{array}{l}
\frac{x_{k+1}-x_{k}}{h}  \in f(t_{k},x_{k})-\rho x_{k}-g(x_{k+1}) \mbox{Sgn}(h(x_{k+1})+\rho x_{k+1}
,\,\,\,\mbox{a.e. on}\,\,(0,T) \\  \\
x(0)=x_{0}
\end{array}\right.
\label{eq:discquatre}
\end{equation}
where the $\rho_{i}$ are the hypomonotonicity constants and $\sum_{i=1}^{m}\rho_{i}=\rho$.  The result of Proposition \ref{proposition2} applies to  (\ref{eq:disctrois}) and (\ref{eq:discquatre}).

\subsection{A simple hypomonotone case}

As shown in Section \ref{section1} on a simple monotone example, in practice the intersection may be computed as follows.   Let us now illustrate this on the following system with hypomonotone multivalued part:
\begin{equation}
\dot{x}(t) \in -(x(t)+1)\mbox{sgn}(x(t))+u(t)
\label{sixe}
\end{equation}
with $x(t) \in \RR$. We may discretize it as:
\begin{equation}
-(x_{k+1}-x_{k}-hu_{k}-h\rho x_{k}) \in h(x_{k+1}+1)\mbox{sgn}(x_{k+1})+h \rho x_{k+1},\,\,k \geq 0,\,\,x_{0}=x(0)
\label{neuf}
\end{equation}
Notice that we may rewrite (\ref{neuf}) as
\begin{equation}
0 \in x_{k+1}-(x_{k}+hu_{k}+h\rho x_{k})+ h(x_{k+1}+1)\mbox{sgn}(x_{k+1})+h \rho x_{k+1}
\label{dix}
\end{equation}

Let us denote the mapping in the right-hand-side of (\ref{dix}) as $F(x_{k+1})$. The set-valued mapping $F(\cdot)$ is {\em strongly monotone} \cite[Definition 2.3.1]{Facchinei} for all $\rho \geq 1$. It follows from \cite[Theorem 2.3.3]{Facchinei} that the {\em generalized equation} $0 \in F(x_{k+1})$ has a unique solution. 

For $u(t)=0$, the following Lemma extends the Lemma~\ref{lemmayahoo}.
\begin{lemma}
For all $ 1>h>0$ and $x_0\in \RR$, there exists $k_0$ such that  $x_{k_0+n}=0$ and  $\displaystyle\frac{x_{k_0+n+1}-x_{k_0+n}}{h}=0$ for all $n \geq 1$.  
\label{lemmayoupi}
\end{lemma}

{\bf Proof:}If $x_0 \in [-h , h]$, then $k_0=0$. Otherwise, for $k<k_0$ and $x_k \notin [-h,h]$, the solution is given by :
\begin{equation}
  \label{eq:solu-hypo1}
  x_{k+1} =\frac{x_k-h\mbox{sgn}(x_k)}{1+h\mbox{sgn}(x_k)}; s_{k+1} = \mbox{sgn}(x_k)
\end{equation}
From the solution~(\ref{eq:solu-hypo1}), the step $k_0$ for which $x_{k_0} \in [-h,h]$ can be easily found.  Let us now consider that $x_{k_0}\in [-h, h]$. The only  possible solution for  
\begin{equation}
 \begin{cases}
   x_{k_0+1}-x_{k_0} = -h(x_{k_0+1}+1) s_{k_0+1}\\[2mm]
   s_{k_0+1} \in \mbox{sgn}(x_{k_0+1})
 \end{cases}
\label{eq:huit-bis}
\end{equation}
is $x_{k_0+1} = 0$ and $s_{k_0+1} = \displaystyle\frac{x_{k_0}}{h}$. For the next iteration, we have to solve
\begin{equation}
 \begin{cases}
   x_{k_0+2} = -h s_{k_0+2}\\[2mm]
   s_{k_0+2} \in \mbox{sgn}(x_{k_0+2})
 \end{cases}
\label{eq:huit-ter}
\end{equation}
and we obtain $x_{k_0+2} = 0 $ and $s_{k_0+2}=0$.  The same holds for all $x_{k_0+n}$,$s_{k_0+n}$, $n \geq 3$, redoing the same reasoning. Clearly then the terms  $\displaystyle\frac{x_{k_0+n+1}-x_{k_0+n}}{h}$ approximating the derivative, are zero for any $h >0$.  \Qend

We conclude that in this case also the system and its derivative are correctly approximated at the zero value on the sliding surface. There is no spurious oscillation around the switching surface.

\section{Detailed Implementation of Implicit Euler Discretization of the general case~(\ref{eq:quatre}) }

\label{sectionImplementation}
 This section is devoted to the implementation and the study  of the numerical algorithm. The interval of integration is $[0,T]$, $T>0$, and a grid $t_{0}=0$, $t_{k+1}=t_{k}+h$, $k \geq 0$, $t_{N}=T$ is constructed. The approximation of a function $f(\cdot)$ on $[0,T]$ is denoted as $f^{N}(\cdot)$, and is a piecewise constant function, constant on the intervals $[t_{k},t_{k+1})$. We denote $f^{N}(t_{k})$ as $f_{k}$. The time-step is $h>0$.

\subsection{Time--discretization }

Starting from  (\ref{eq:quatre}), let us introduce a new notation, 
\begin{equation}
\begin{array}{l}
\dot{x}(t) = f(x(t),t) - g(x(t)) s(t) \\[2mm]
y(t) = h(x(t)) \\[2mm]
s(t) \in \mbox{Sgn}(y(t))
\end{array}
\label{sept-bis-bis}
\end{equation}
where $s(t) \in \RR^m$  and $y(t) \in \RR^m$ are  complementary variables related through the $\mbox{Sgn}(\cdot)$ multi--valued mapping.   According to the class of systems (\ref{eq:deux}), (\ref{eq:trois}) or (\ref{eq:quatre}) that we are studying the functions $f(\cdot)$ and $g(\cdot)$ are defined either in a fully nonlinear framework or by affine functions. We present the time-discretization in its full generality and specialize the algorithms for each case in Section~\ref{Sec:Spec}.

Let us now proceed with the time discretization of (\ref{sept-bis-bis}) by a fully implicit scheme : 
\begin{equation}
  \begin{array}{l}
    \label{eq:toto1}
     x_{k+1} = x_{k} +h f(x_{k+\theta},t_{k+\theta}) - h g(x_{k+\gamma})  s_{k+1} \\[2mm]
     y_{k+1} =  h(x_{k+1}) \\[2mm] 
     s_{k+1} \in \mbox{Sgn}{(y_{k+1} )}
  \end{array}
\end{equation}
where $x_{k+\theta} = \theta x_{k+1} + (1-\theta) x_{k}$, $x_{k+\gamma} = \gamma x_{k+1} + (1-\gamma) x_{k}$,  and $t_{k+\theta} = \theta t_{k+1} + (1-\theta)  t_{k+1}$, with $\theta = [0,1]$ and $\gamma \in [0,1]$. As in \cite{acary2008}, we call the problem \eqref{eq:toto1} the ``one--step nonsmooth problem''.

 This time-discretization is slightly more general than a standard implicit Euler scheme. The main discrepancy lies in the choice of a $\theta$-method to integrate the nonlinear term. For $\theta=0$, we retrieve the explicit integration of the smooth and  single valued term $f(\cdot)$. Moreover for $\gamma =0$, the term $g(\cdot)$ is explicitly evaluated. The flexibility in the choice of $\theta$ and $\gamma$ allows the user to improve and control the accuracy, the stability and the numerical damping of the proposed method. For instance, if the smooth dynamics given by $f(\cdot)$ is stiff, or if we have to use large step sizes for practical reasons, the choice of $\theta > 1/2$ offers better stability properties with  respect to $h$.

\subsection{Mixed Complementarity Problem}
 The so-called "one--step nonsmooth problem''~(\ref{eq:toto1})  appears to be a Mixed Complementarity Problem (MCP) that we have to solve at each time--step. Let us define what is a MCP :

\begin{definition}[MCP]
  Given a function $f : \RR^q \rightarrow \RR^q $ and lower and upper bounds $l,u\in \bar\RR^{q}$, find $z\in \RR^q$, $w,v \in \RR^q_+$ such that
  \begin{equation}
    \label{eq:MCP1}
    \left\{\begin{array}[]{c}
      F(z) = w-v \\[2mm]
      l \leq z \leq u \\[2mm]
      (z-l)^Tw=0 \\[2mm]
      (u-z)^Tv = 0
    \end{array}\right.
  \end{equation}
where  $\bar\RR = \RR \cup \{+\infty,-\infty\}$. 
\end{definition}

Note that the problem~(\ref{eq:MCP1})  implies  that
\begin{equation}
  \label{eq:MCP10}
  -F(z) \in N_{[l,u]}(z).
\end{equation}
 The relation (\ref{eq:MCP10}) is equivalent to the MCP~(\ref{eq:MCP1}) if we assume that  $w$ is the positive part of $F(z)$, that is $w = F^+(z)=max(0, F(z))$ and $v$ is the negative part of $F(z)$, that is $v = F^-(z)= max(0, -F(z))$.

\paragraph{The One--step nonsmooth problem as a MCP}
Let us  define the MCP by
\begin{equation}
  \label{eq:MCP3}
  \begin{array}[l]{c}
    z = \left[
      \begin{array}{l}
        x_{k+1} \\ s_{k+1}
\end{array}
\right] \\ \\
    F(z) =
    \left[\begin{array}{c}
      x_{k+1}-x_k-h f(x_{k+\theta}, t_{k+\theta}) + hg(x_{k+\gamma}) s_{k+1} \\
      - h(x_{k+1})
    \end{array}\right]
    \\ \\
    l_i  =
    \begin{cases}
      -\infty, i = 1\ldots n \\
      -1, i = n+1\ldots m\\
    \end{cases},\quad
     u_i  =
    \begin{cases}
      +\infty, i = 1\ldots n \\
      +1, i = n+1\ldots m \\
    \end{cases} 
  \end{array}
\end{equation}
If $z$ solves the MCP~(\ref{eq:MCP3}), the bounds $u$ and $l$ and the condition $ (z-l)^Tw=0, (u-z)^Tv = 0$ imply that 
      \begin{equation}
        \label{eq:MCP4}     
        w_i  =
    \begin{cases}
     0, i = 1\ldots n \\
     y_i^- \geq 0, i = n+1\ldots m \\
    \end{cases},\quad
    v_i  =
    \begin{cases}
     0, i = 1\ldots n \\
     y_i^+\geq 0, i = n+1\ldots m \\
    \end{cases} 
\end{equation}
The MCP is then given by
\begin{equation}
  \label{eq:MCP5}
  \begin{array}{c}
    \left[\begin{array}{c}
        x_{k+1}-x_k-hf(x_{k+\theta}) + hg(x_{k+\gamma}) s_{k+1} \\
        - h(x_{k+1})
      \end{array}\right]
    = 
    \left[\begin{array}{c}
        0 \\
        y^-(x_{k+1})- y^+(x_{k+1}) = -  y(x_{k+1}) 
      \end{array}\right] \\[2mm] 
    -1 < s_{k+1} < 1 \\[2mm]
    (s_{k+1}+1)^T y^-(x_{k+1}) =0 \\[2mm]
    (1-s_{k+1})^T y^+(x_{k+1}) =0
 \end{array}
\end{equation}
It is clear that the problem~(\ref{eq:toto1}) is equivalent to the MCP defined by  (\ref{eq:MCP3}). The results of existence and uniqueness of solution of (\ref{eq:MCP1}) or equivalently~(\ref{eq:MCP10}) are related to the monotonicity properties of $F$ (P-properties or coherent orientations of the associated affine map for particular cases of bounds and affine function $F(\cdot)$. Without entering into further details, we refer to \cite{Harker.Pang1990,Facchinei} for the main results.

\paragraph{Numerical  Solvers}

The MCP (\ref{eq:MCP1}) can be solved by a large family  of solvers based on Newton--type Methods and interior-points techniques.  We refer to  \cite{Billups.ea1997} for a comparison of several solvers based on Newton's method. The numerical implementation of the MCP solvers are often based on the computation of the Jacobian matrix of the function $F(\cdot)$ with respect to $z$. The Jacobian matrix is explicitly given in the case defined in (\ref{eq:MCP3}) by
\begin{equation}
  \label{eq:JacMCLP}
  \nabla_z F(z) =
  \left[\begin{array}{cc}
    I - h  \theta \nabla_{x} f(x_{k+\theta}, t_{k+\theta})  + h \gamma \nabla_{x} g(x_{k+\gamma}) \scontract s_{k+1} & h g(x_{k+\gamma})  \\ \\
    \nabla_{x} h(x_{k+1} ) & 0
  \end{array}\right],
\end{equation}
where  $\scontract$ denotes the simple contracted tensor product and the third--order tensor $\nabla_{x} g(x)$ is the Jacobian of $g$ with the respect $x$  given by the following component:
\begin{equation}
  \label{eq:NL00}
   (\nabla_{x} g(x))_{klp} = \frac{ \partial g_{kl}(x)}{\partial x_p}.
\end{equation}

\subsection{Newton's linearization and Mixed Linear Complementarity Problems}

Due to the fact that  two of the systems that are studied in this paper involve affine functions $f(\cdot)$ and $g(\cdot)$, we propose to solve the "one--step nonsmooth problem'' (\ref{eq:toto1}) by performing an external Newton linearization, which yields a Mixed Linear Complementarity Problems (MLCP).

 \paragraph{Newton's linearization} The first line of the  problem~(\ref{eq:toto1}) can be written under the form of a residue $\mathcal R$ depending only on $x_{k+1}$ and $s_{k+1}$ such that 
\begin{equation}
  \label{eq:NL3}
  \mathcal R (x_{k+1},s_{k+1}) =0
\end{equation}
with $\mathcal R(x,s) = x - x_{k} -h f( \theta x + (1-\theta) x_{k}, t_{k+\theta}) + hg(\gamma x + (1- \gamma) x_{k})  s$.
The solution of this system of nonlinear equations is sought as a limit of the sequence $\{ x^{\alpha}_{k+1},s^{\alpha}_{k+1} \}_{\alpha \in \NN}$ such that
 \begin{equation}
   \label{eq:NL7}
   \begin{cases}
     x^{0}_{k+1} = x_k \\ \\
     \mathcal R_L( x^{\alpha+1}_{k+1},s^{\alpha+1}_{k+1} ) =x^{\alpha}_{k+1}-x_k-h f (x^{\alpha}_{k+\theta})  + \nabla_{x} \mathcal R (x^{\alpha}_{k+1},s^{\alpha}_{k+1} )(x^{\alpha+1}_{k+1}-x^{\alpha}_{k+1} ) +  h g(x^{\alpha}_{k+\gamma})s^{\alpha+1}_{k+1} =0
 \end{cases}
\end{equation}
The computation of the Jacobian of $\mathcal R$ with respect to $x$, denoted by $M(x,s)$ leads to 
\begin{equation}
   \label{eq:NL9}
   \begin{array}{l}
    M(x,s)= \nabla_{x} \mathcal R (x,s)= I - h  \theta \nabla_{x} f( \theta x + (1-\theta) x_{k}, t_{k+\theta} )  + h \gamma \nabla_{x} g(\gamma x + (1- \gamma) x_{k}) \scontract s .
 \end{array}
\end{equation}
At each time--step, we have to solve the following linearized problem,
\begin{equation}
   \label{eq:NL10}
     x^{\alpha}_{k+1}-x_k- h f (x^{\alpha}_{k+\theta},t_{k+\theta})  + M (x^{\alpha}_{k+1}, s^{\alpha}_{k+1})   (x^{\alpha+1}_{k+1}-x^{\alpha}_{k+1} ) +  h g(x^{\alpha}_{k+\gamma})s^{\alpha+1}_{k+1} =0 ,
\end{equation}
that is
\begin{equation}
 x^{\alpha+1}_{k+1}= x^{\alpha}_{k+1} + M^{-1}(x^{\alpha}_{k+1},s^{\alpha}_{k+1})\left[x_k- x^{\alpha}_{k+1} +h f (x^{\alpha}_{k+\theta},t_{k+\theta}) -  hg(x^{\alpha}_{k+\gamma}) s^{\alpha+1}_{k+1}  \right].\label{eq:NL11} 
\end{equation}
The matrix $M$ is clearly non singular for small $h$.
The same operation is performed with the second equation of (\ref{eq:toto1}) leading to the following linearized equation
\begin{equation}
  \label{eq:NL12}
  y^{\alpha+1}_{k+1} =  y^{\alpha}_{k+1} +\nabla_x h(x^{\alpha}_{k+1})\left[x^{\alpha+1}_{k+1}-x^{\alpha}_{k+1}\right]
\end{equation}
Inserting (\ref{eq:NL11}), we get the following linear relation between $y^{\alpha+1}_{k+1}$ and   $s^{\alpha+1}_{k+1}$, 
\begin{equation}
  \begin{array}[l]{l}
 y^{\alpha+1}_{k+1} =  y^{\alpha}_{k+1} +\nabla_x h(x^{\alpha}_{k+1})
 \left[M^{-1}(x^{\alpha}_{k+1},s^{\alpha}_{k+1}) (x_k- x^{\alpha}_{k+1} +h f (x^{\alpha}_{k+\theta},t_{k+\theta}) -  h g(x^{\alpha}_{k+\gamma})  s^{\alpha+1}_{k+1}  )\right]
\end{array}
\label{eq:NL13} 
\end{equation}

\paragraph{Mixed linear complementarity problem (MLCP)}To summarize, the problem to be solved in each Newton iteration is:\\{
  \begin{minipage}[l]{1.0\linewidth}
    \begin{equation}
      \left\{\begin{array}[l]{l}
        y^{\alpha+1}_{k+1} =   - W^{\alpha+1}_{k+1}  s^{\alpha+1}_{k+1} + b^{\alpha+1}_{k+1} \\ \\
        s^{\alpha+1}_{k+1} \in \mbox{Sgn}{( y^{\alpha+1}_{k+1} )}
      \end{array}\right.
      \label{eq:NL14} 
    \end{equation}
  \end{minipage}
}
with $W\in \RR^{m\times m}$ and $b\in\RR^{m}$ defined by
\begin{equation}
  \label{eq:NL15}
 \begin{array}[l]{l}
   W^{\alpha+1}_{k+1} = h \nabla_x h(x^{\alpha}_{k+1}) M^{-1}(x^{\alpha}_{k+1},s^{\alpha}_{k+1})  g(x^{\alpha}_{k+\gamma}) \\ \\
   b^{\alpha+1}_{k+1} = y^{\alpha}_{k+1} + \nabla_x h(x^{\alpha}_{k+1})
 \left[M^{-1}(x^{\alpha}_{k+1},s^{\alpha}_{k+1}) (x_k- x^{\alpha}_{k+1} +h f (x^{\alpha}_{k+\theta},t_{k+\theta}) )\right]
\end{array}
\end{equation}

The problem~(\ref{eq:NL14}) is equivalent to a MLCP which can be solved under suitable assumptions by many linear complementarity solvers such as pivoting techniques, interior point techniques and splitting/projection strategies. The  reformulation into a standard MLCP follows the same line as for the MCP in the previous section. One obtains,
    \begin{equation}
      \left\{\begin{array}[l]{l}
        y^{\alpha+1,+}_{k+1}-y^{\alpha+1}_{k+1}  =   - W^{\alpha+1}_{k+1}  s^{\alpha+1}_{k+1} + b^{\alpha+1}_{k+1} \\ \\
        0 \leq (s^{\alpha+1}_{k+1}+1)\perp  y^{\alpha+1,-}_{k+1}  \geq 0 \\[2mm]
        0 \leq  (1-s^{\alpha+1}_{k+1})\perp y^{\alpha+1,+}_{k+1} \geq 0 
      \end{array}\right.
      \label{eq:MLCP1-bis} 
    \end{equation}
As for the  MCP, there exists numerous methods to numerically solve MLCP. In the worst case when the matrix $W^{\alpha+1}_{k+1}$ has no special properties, the MCLP can be always solved by enumerative solvers for which various implementations can be found. With some positiveness properties \cite{Facchinei}, standard methods for LCP\cite{Cottle.Pang.ea1992} can be straightforwardly extended. Among these methods, we can cite the family of projection/splitting methods, interior point methods and semi-smooth Newton methods (see \cite{acary2008} for an overview).

\subsection{The special cases of the  affine  systems}
\label{Sec:Spec}

In this section, we specify the time--discretization  to the two other classes of systems~(\ref{eq:deux}) and ~(\ref{eq:trois}) and for particular value of $\theta$ and $\gamma$. 

\subsubsection{Time--discretization of the   system (\ref{eq:deux})}

For the system (\ref{eq:deux}), the function $g(x)$ is reduced to the matrix $B$ and the function $h(x)$ is affine, that is $h(x) = C x + D $.  The matrix $W^{\alpha+1}_{k+1}$ and $b^{\alpha+1}_{k+1}$ are then given by 
\begin{equation}
  \label{eq:NL15-deux}
 \left\{\begin{array}[l]{l}
   W^{\alpha+1}_{k+1} = h C M^{-1}(x^{\alpha}_{k+1},s^{\alpha}_{k+1},t_{k+1})  B \\ \\
   b^{\alpha+1}_{k+1} = y^{\alpha}_{k+1} + C
 \left[M^{-1}(x^{\alpha}_{k+1},t_{k+1}) (x_k- x^{\alpha}_{k+1} +h f (x^{\alpha}_{k+\theta},t_{k+\theta}) )\right]
\end{array}\right.
\end{equation}
with
\begin{equation}
  \label{eq:NL16-deux}
   M(x,t)= I - h  \theta \nabla_{x} f( \theta x + (1-\theta) x_{k},\theta t + (1-\theta) t_{k} )
\end{equation}

If  $CB  >0$ then the matrix $W^{\alpha+1}_{k+1}$ is also positive definite for sufficiently small $h$. This result ensures the existence and uniqueness of a solution of the MCLP. Furthermore, standard pivoting techniques such as Lemke's method or projection/splitting such as Projected Successive Over-Relation (PSOR) compute the solution.

\paragraph{Semi-implicit discretization with $\theta =0$} The matrix $M$ is then  reduced to
\begin{equation}
  \label{eq:NL17-deux}
   M(x,t)= I
\end{equation}
and 
\begin{equation}
  \label{eq:NL18-deux}
 \begin{array}[l]{l}
   W =  h C B \\ \\
   b_{k+1} =  y_k  +h C f(x_{k},t_{k}) 
\end{array}
\end{equation}
In this particular case, there is no need to perform some Newton iterations because the system to be solved at each time--step is linear. Furthermore, the MLCP has a solution for any $h>0$ under the assumptions that   $CB  >0$.

\paragraph{Fully implicit  discretization with an affine function}
 The same conclusion can be drawn if  $\nabla_{x} f( \theta x + (1-\theta) x_{k},\theta t + (1-\theta) t_{k} )$ is equal to a constant matrix $E$ that is when $f(\cdot)$ is linear time-invariant and given by $f(x,t) = E x(t) + F $. In this case, the matrix $M$ reduces to 
 \begin{equation}
  \label{eq:NL19-deux}
   M(x,t)= I -h \theta E
\end{equation}
and 
\begin{equation}
  \label{eq:NL20-deux}
 \begin{array}[l]{l}
   W = h C (I-h\theta E)^{-1} B \\ \\
   b_{k+1} = y_{k} +  h C (I-h\theta E)^{-1}
 \left[  E x_{k} + a\right]
\end{array}
\end{equation}

\subsubsection{Time--discretization of the  system  (\ref{eq:trois})}

 For the system (\ref{eq:trois}), we recall that $g(\cdot)$ and $h(\cdot)$ are given by
\begin{equation}
  \label{eq:NL-trois}
  \begin{array}{lcl}
    g(x) &=& [g_{\bullet j}(x) = A_j x + B_j, j = 1 \ldots m ] \in \RR^{n\times m}, \\
    h(x) &=& [h_{i}(x) = C_{i} x + D_{i}, i = 1 \ldots m ] \in \RR^{m}.
  \end{array}
\end{equation}
The components of $g(\cdot)$ can be explicitly expressed by
\begin{equation}
  \label{eq:NL-trois3}
   g_{kl}(x) = \sum_{p=1}^n A_{l,kp} x_p + B_l,
\end{equation}
and therefore, the Jacobian of $g(\cdot)$ is given by
\begin{equation}
  \label{eq:NL-trois2}
  (\nabla_{x} g(x))_{klp} = \frac{ \partial g_{kl}(x)}{\partial x_p} = A_{l,kp}.
\end{equation}
 The Jacobian of  $h$ takes the following simple form:
\begin{equation}
  \label{eq:NL-trois1}
  \begin{array}{lcl}
    \nabla h (x) &=& C = [ C_{i}, i = 1 \ldots m] \in \RR^{m\times n}.
\end{array}
\end{equation}

After the newton linearization, we have to solve at each Newton's loop the MLCP (\ref{eq:NL14}) with
\begin{equation}
  \label{eq:NL-trois4}
  \begin{array}[l]{l}
    W^{\alpha+1}_{k+1} = h C M^{-1}(x^{\alpha}_{k+1},s^{\alpha}_{k+1},t_{k+1})  g(x^{\alpha}_{k+\gamma}) \\ \\
    b^{\alpha+1}_{k+1} = y^{\alpha}_{k+1} + C
    \left[M^{-1}(x^{\alpha}_{k+1},t_{k+1}) (x_k- x^{\alpha}_{k+1} +h f (x^{\alpha}_{k+\theta},t_{k+\theta}) )\right]
  \end{array}
\end{equation}

\paragraph{Semi--implicit discretization with $\theta=\gamma=0$} If $\gamma$ and $\theta$ vanish, the residue $\mathcal R$ given by
\begin{equation}
\mathcal R(x,s) = x - x_{k} -h f( x_{k}, t_{k}) + hg(x_{k})  s\label{eq:Residue-semi}
\end{equation}
is linear in $x$ and $s$. In this particular case, there is no need to perform Newton's iterations. The MCLP defined by (\ref{eq:NL-trois4}) can be simplified to
\begin{equation}
  \label{eq:NL-trois5}
    W_{k+1} = h C g(x_{k}) \\ \\
    b_{k+1} = y_{k} + C
    \left[h f (x_{k},t_{k}) )\right]
\end{equation}

\subsection{Algorithms} 

We propose in this section two algorithms to sum-up the numerical implementation of the implicit Euler time--stepping scheme.  The Algorithm~\ref{Algo:EulerSliding-MCP} describes the implementation with a generic MCP solver and the Algorithm~\ref{Algo:EulerSliding-MCLP} describes the  numerical implementation of the algorithm with an external Newton linearization and a MCLP solver.
\begin{algorithm}[htbp]
   \begin{algorithmic}
 { \sf 
    \REQUIRE System definition: $\sf f, g, h$ 
    \REQUIRE $\sf x(0)$ the initial condition
    \REQUIRE $\sf t_0, T$ time--integration interval
    \REQUIRE $\sf h$ time--step 
    \REQUIRE $\sf \theta, \gamma$ numerical integration parameters
    \ENSURE  $\sf (\{ x_k\}, \{ s_k\},\{y_k \}), k \in\{1,2, \ldots \}$ 
    \STATE $ \sf $
    \STATE $\sf k \leftarrow 0;\quad  x_{0} \leftarrow x(0);\quad \sf y_{0} \leftarrow y(0) =h(x(0));\quad \sf tau_{0} \leftarrow 0 $  
    \STATE //\textit{ Time integration loop}
     \WHILE {$\sf t_k < T$} 
         \STATE Solve the MCP~(\ref{eq:MCP5}) for $\sf x_{k+1}, s_{k+1}, y_{k+1}$ with $\sf F, l$ and $\sf u$ given by ~(\ref{eq:MCP3})  and the Jacobian $\nabla_z F(z)$ given by (\ref{eq:JacMCLP}).
         \STATE //\textit{Update}
         \STATE $\sf x_{k} \leftarrow x_{k+1};\quad \sf s_{k} \leftarrow s_{k+1} ;\quad \sf y_{k} \leftarrow y_{k+1} $
         \STATE //\textit{time iteration}
         \STATE $\sf t_{k} \leftarrow t_{k+1};\quad \sf k \leftarrow k+1$
      \ENDWHILE
     \STATE $ \sf$}
   \end{algorithmic}
   \caption{Implicit Euler time-discretization with a generic MCP solver}  
 \label{Algo:EulerSliding-MCP}
\end{algorithm}
\begin{algorithm}[htbp]
   \begin{algorithmic}
 { \sf 
    \REQUIRE System definition: $\sf f, g, h$ 
    \REQUIRE $\sf x(0)$ the initial condition
    \REQUIRE $\sf t_0, T$ time--integration interval
    \REQUIRE $\sf h$ time--step
    \REQUIRE $\sf \theta, \gamma$ numerical integration parameters 
    \REQUIRE $\sf \varepsilon $ Newton's method tolerance
    \ENSURE  $\sf (\{ x_k\}, \{ s_k\},\{y_k \}), k \in\{1,2, \ldots \}$ 
    \STATE $ \sf $
    \STATE $\sf k \leftarrow 0;\quad  x_{0} \leftarrow x(0);\quad \sf y_{0} \leftarrow y(0) =h(x(0));\quad \sf tau_{0} \leftarrow 0 $  
    \STATE $ $
    \STATE //\textit{ Time integration loop}
     \WHILE {$\sf t_k < T$}
         \STATE $\sf \alpha \leftarrow 0;\quad\sf x^0_{k+1} \leftarrow x^0_{k};\quad s^0_{k+1} \leftarrow s^0_{k};\quad \sf y^0_{k+1} \leftarrow y^0_{k} $
        \STATE //\textit{Newton's loop}
          \WHILE {$\sf  \| \mathcal R(x^{\alpha}_{k+1},s^{\alpha}_{k+1} )\| > \varepsilon$}  
         \STATE  $\sf M^{-1}(x^{\alpha}_{k+1},s^{\alpha}_{k+1} ) \leftarrow (I - h \theta\nabla_x f(x^{\alpha}_{k+\theta},t_{k+1}) - h \gamma \nabla_x g(x^{\alpha}_{k+\theta}) \scontract  s^{\alpha}_{k+1}        )^{-1} $.
         \STATE  $\sf W^{\alpha+1}_{k+1} \leftarrow  h \nabla_x h(x^{\alpha}_{k+1}) M^{-1}(s^{\alpha}_{k+1},s^{\alpha}_{k+1})  g(x^{\alpha}_{k+1})$
           \STATE  $\sf b^{\alpha+1}_{k+1} \leftarrow y^{\alpha}_{k+1} + \nabla_x h(x^{\alpha}_{k+1}) M^{-1}(s^{\alpha}_{k+1})\,\left[x_k- x^{\alpha}_{k+1} + h f (x^{\alpha}_{k+1}) \right]$      
           \STATE $ $
           \STATE Solve the MLCP (\ref{eq:MLCP1-bis}) for $\sf y^{\alpha+1}_{k+1}, s^{\alpha+1}_{k+1}$
           \STATE $ $
          \STATE $\sf s^{\alpha}_{k+1} \leftarrow s^{\alpha+1}_{k+1};\quad  y^{\alpha}_{k+1} \leftarrow x^{\alpha+1}_{k+1} $
         \STATE $\sf x^{\alpha}_{k+1}\leftarrow x^{\alpha}_{k+1} + M^{-1}(s^{\alpha}_{k+1})\left[x_k- x^{\alpha}_{k+1} +h f (x^{\alpha}_{k+1}) +  hg(x^{\alpha}_{k+1})s^{\alpha+1}_{k+1}  \right]$
         \STATE $\sf \alpha \leftarrow \alpha+1$         
         \ENDWHILE
         \STATE //\textit{Update}
         \STATE $\sf x_{k+1}\leftarrow x^{\alpha}_{k+1};\quad \sf s_{k+1} \leftarrow s^{\alpha+1}_{k+1} ;\quad \sf y_{k+1} \leftarrow y^{\alpha+1}_{k+1} $
         \STATE //\textit{time iteration}
         \STATE $\sf t_{k} \leftarrow t_{k+1};\quad \sf k \leftarrow k+1$
      \ENDWHILE
     \STATE $ \sf$}
   \end{algorithmic}
   \caption{Implicit Euler time-discretization with an external Newton loop and a MLCP solver }  
 \label{Algo:EulerSliding-MCLP}
\end{algorithm}

In the case of the system~(\ref{eq:deux}) with an affine function $f(\cdot)$ or $\theta =0$, the MLCP matrix $W$ can be computed before the beginning of the time loop, saving a lot of computing effort.  In the case of the system (\ref{eq:trois}) with $\theta=\gamma=0$, the MLCP matrix $W$ can be computed before the beginning of the Newton loop.
\clearpage


\section{Numerical experiments}
\label{secSimus}

Let us illustrate the above developments with numerical integrations performed with the {\sc siconos} software platform of the INRIA\footnote{\url{http://siconos.gforge.inria.fr/}} \cite{acary2008,Acary.Perignon2007} which is designed for the simulation of multivalued nonsmooth systems. 

\subsection{Chattering free stabilization}

\begin{figure}[htbp]
 \psfrag{time}[][]{\Huge time $ \sf t$}
 \psfrag{x1}[][]{\Huge$\sf x_1(t)$}
 \psfrag{x2}[][]{\Huge$\sf x_2(t)$}
 \psfrag{s1}[][]{\Huge$\sf s_1(t)$}
 \psfrag{s2}[][]{\Huge$\sf s_2(t)$}
 \psfrag{xv}[][]{\Huge state}
 \centering
 \subfigure[$h=0.3$. Explicit Euler]
   {\label{Fig:exGaliasYu2007a}\includegraphics[angle=-90,width=0.40\textwidth]{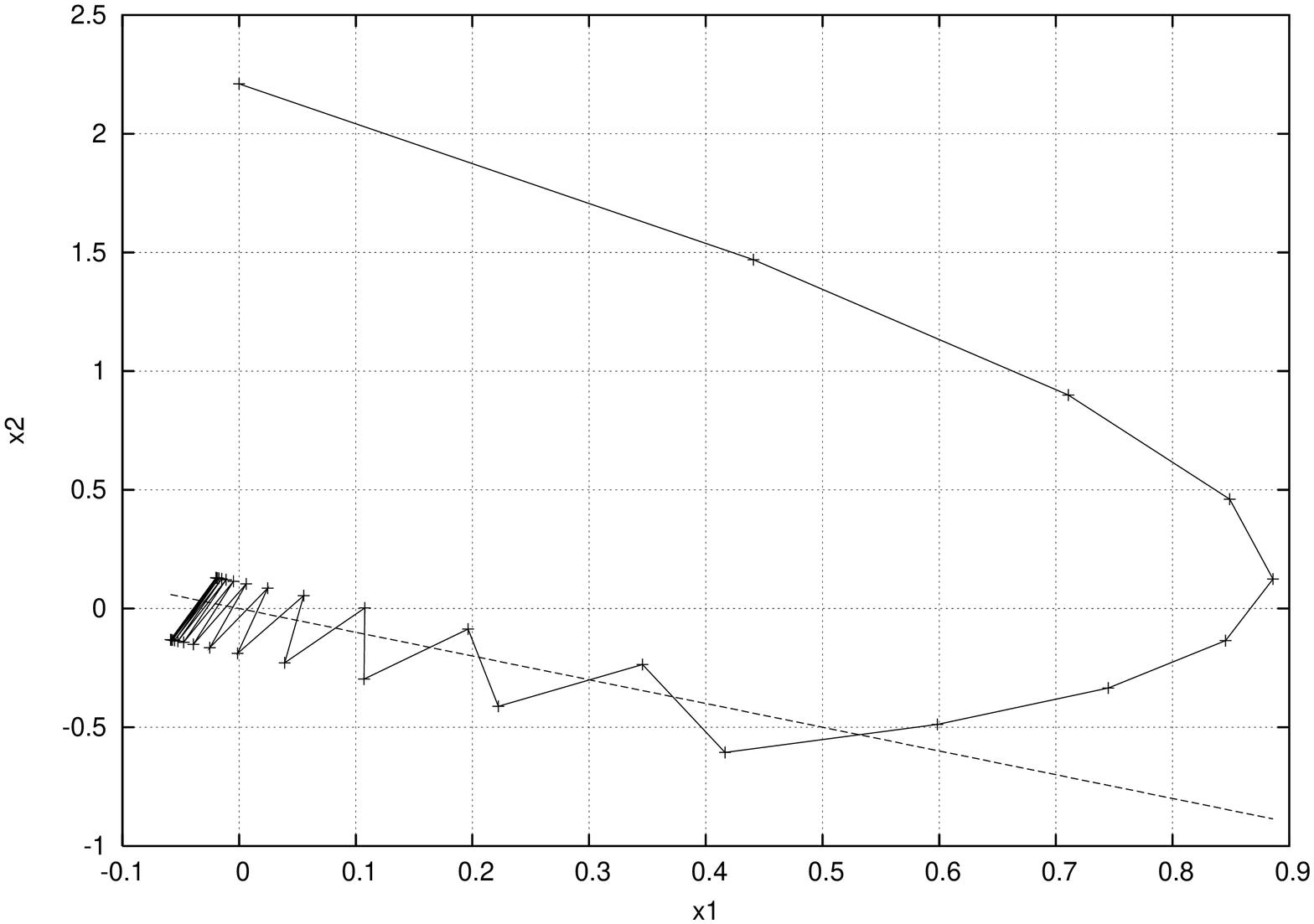}}
 \subfigure[$h=0.1$. Explicit Euler]
   {\label{Fig:exGaliasYu2007b}\includegraphics[angle=-90,width=0.40\textwidth]{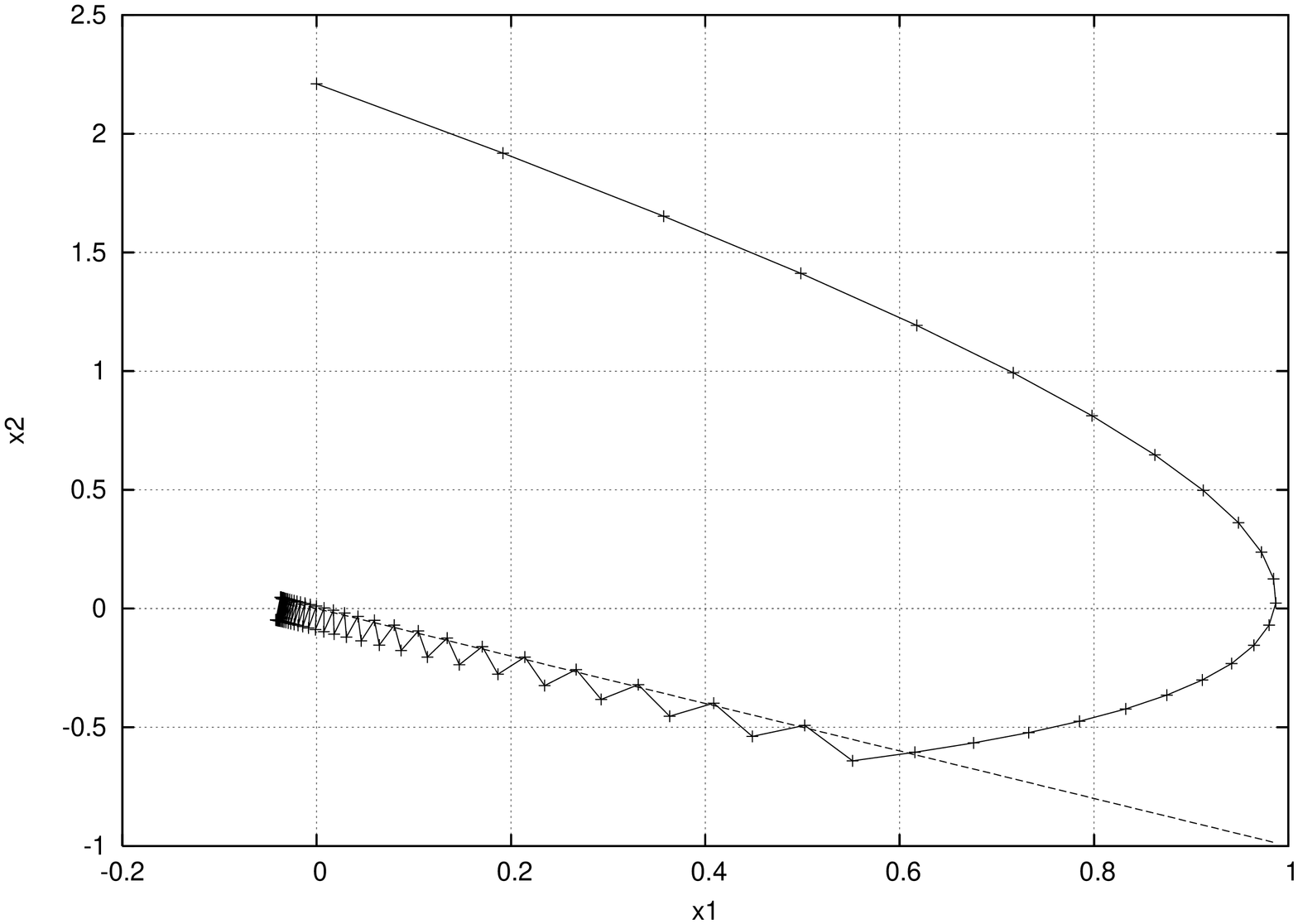}}
 \subfigure[$h=1$. Implicit Euler]
   {\label{Fig:GaliasYu2007a}\includegraphics[angle=-90,width=0.40\textwidth]{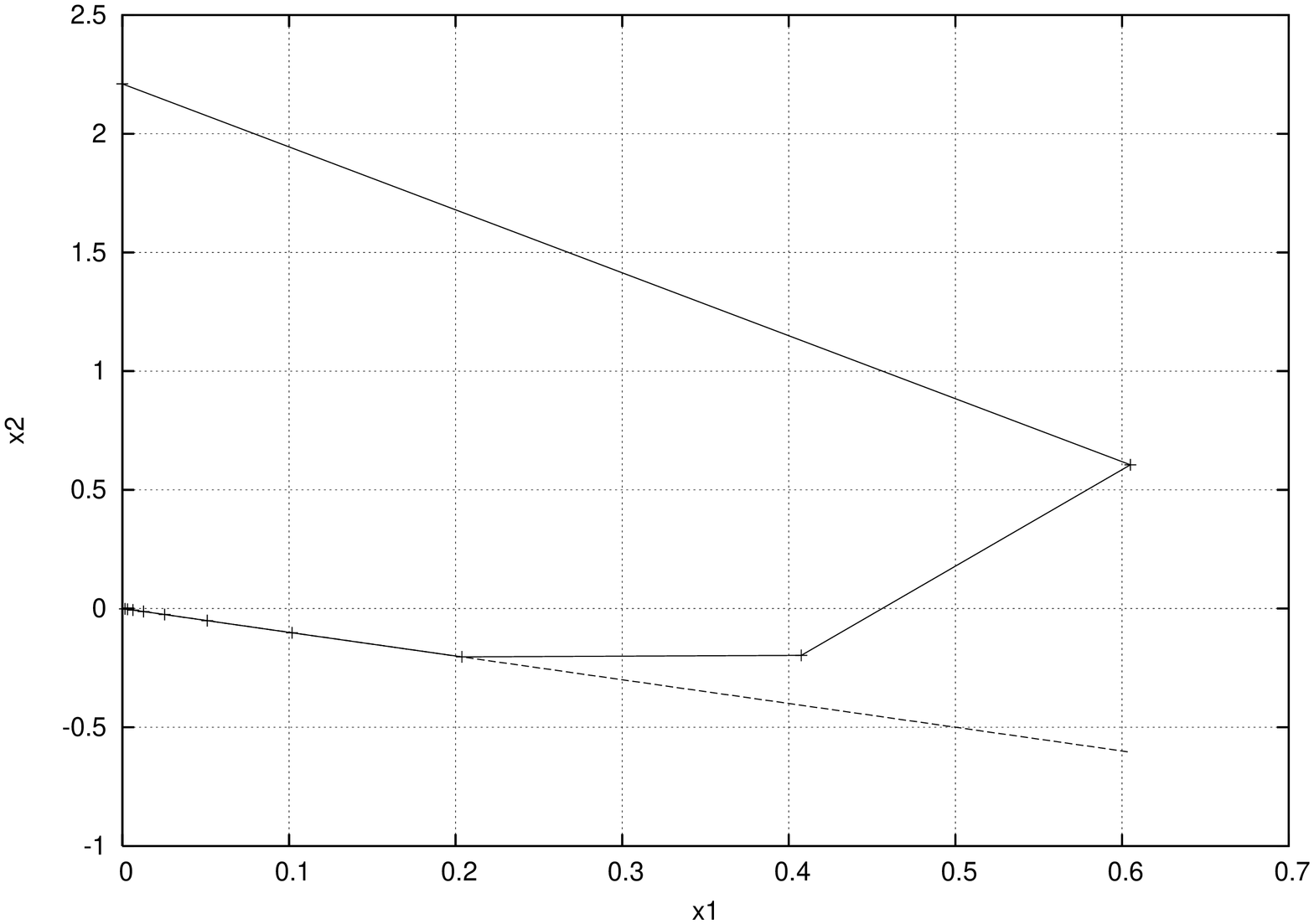}}
 \subfigure[$h=0.3$. Implicit Euler]
   {\label{Fig:GaliasYu2007b}\includegraphics[angle=-90,width=0.40\textwidth]{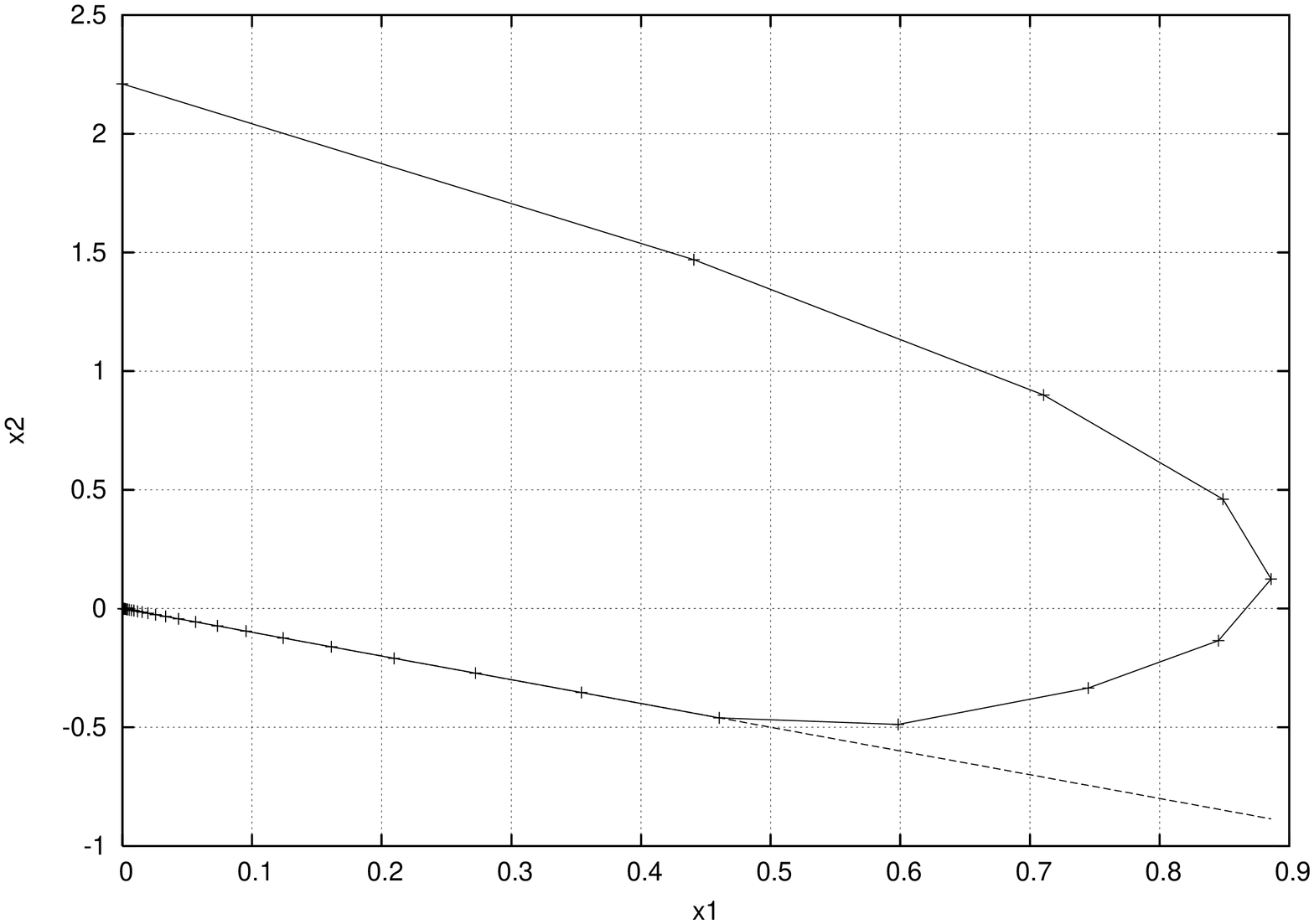}}
   \subfigure[$h=0.1$. Implicit Euler]
   {\label{Fig:GaliasYu2007c}\includegraphics[angle=-90,width=0.40\textwidth]{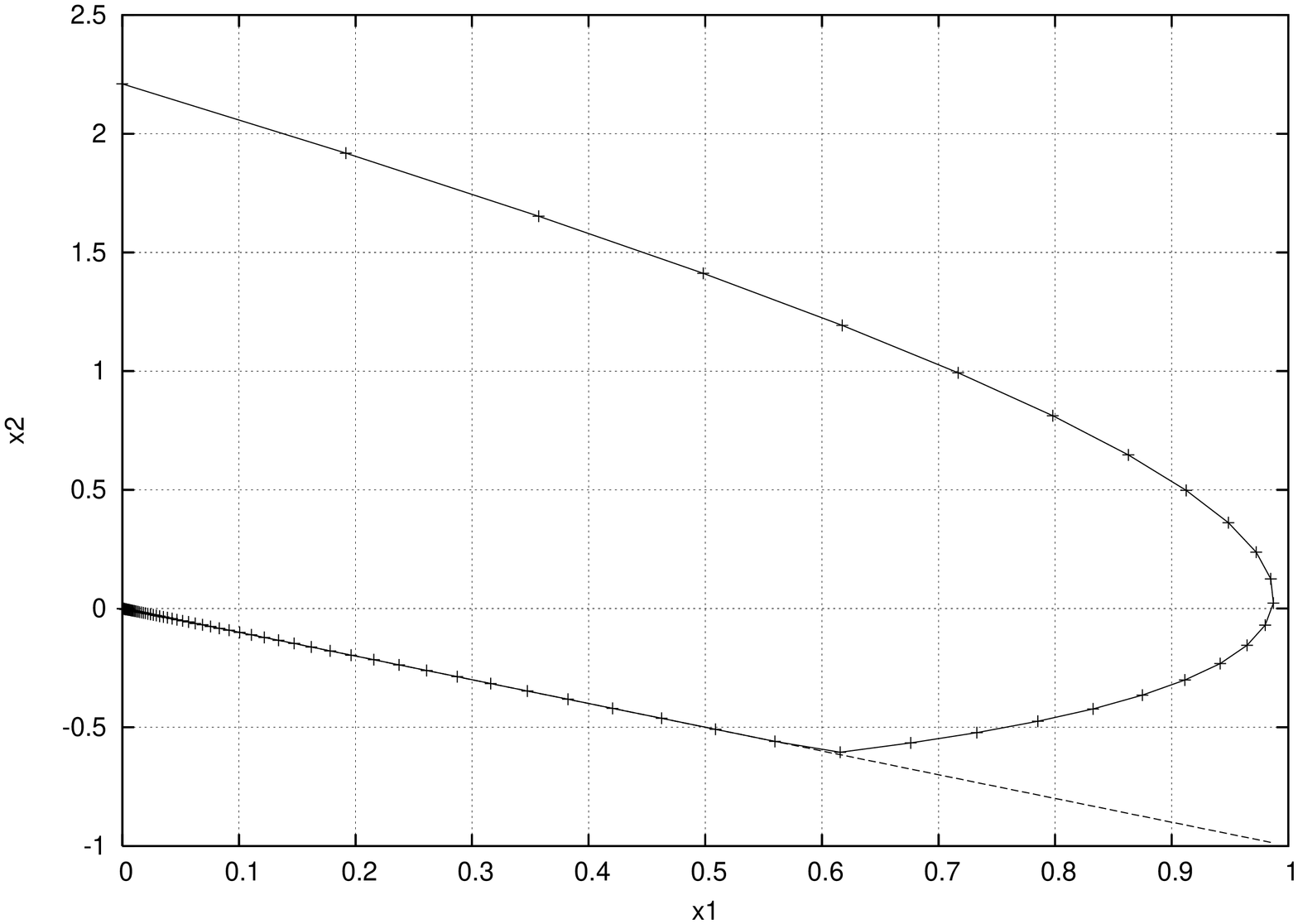}}
    \subfigure[$h=0.05$. Implicit Euler]
   {\label{Fig:GaliasYu2007d}\includegraphics[angle=-90,width=0.40\textwidth]{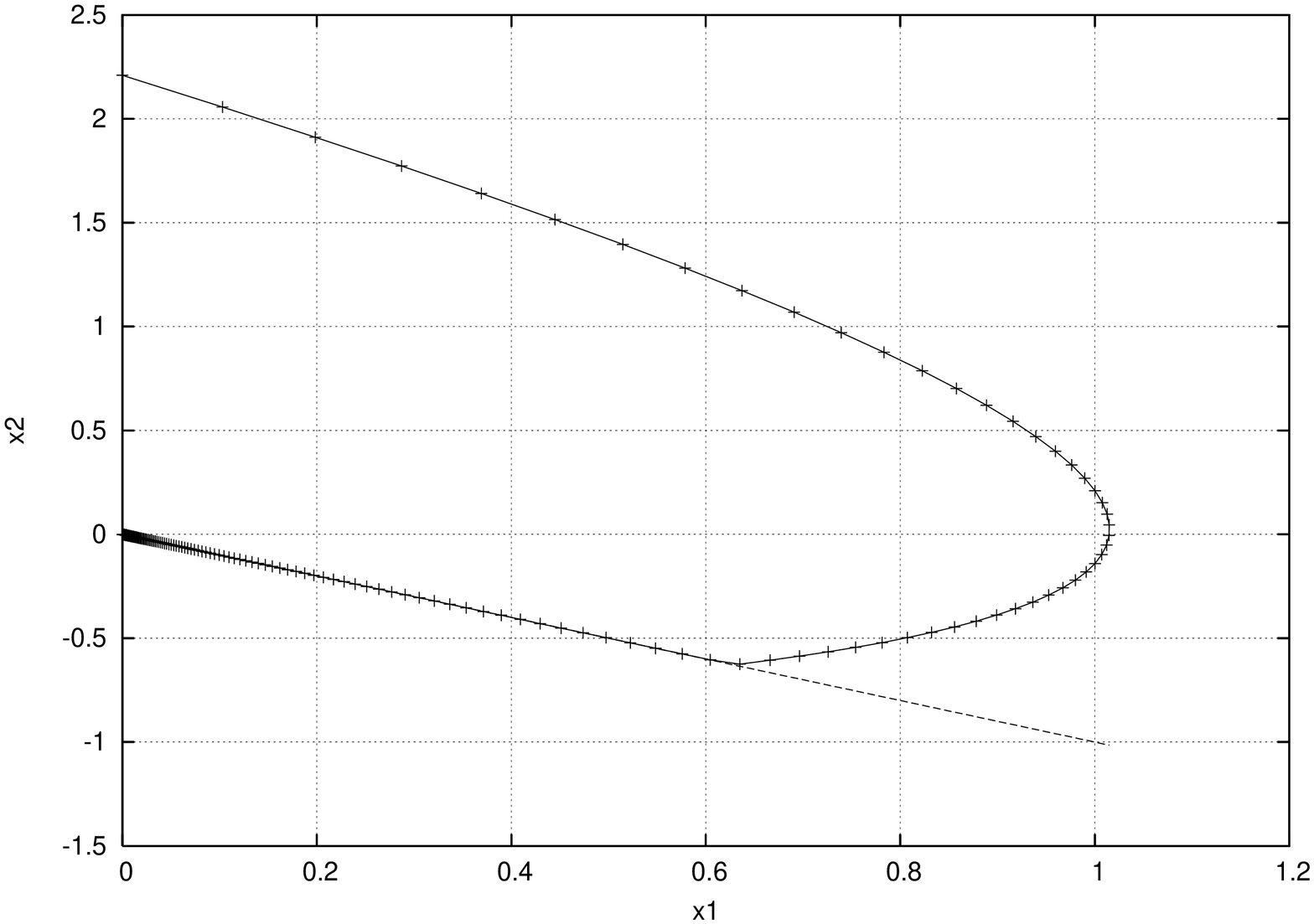}}
 \caption{Equivalent control based SMC, $c_1=1, \alpha =1$ and $x_0 = [0, 2.21]^T$. State $x_1(t)$ versus $x_2(t)$.}  \label{Fig:GaliasYu2007}
\end{figure}

Let us consider the following continuous--time closed loop system from~\cite{galias2007} given by
\begin{equation}
 \label{eq:GaliasYu2007-1}
 \dot x  = \left[
   \begin{array}{cc}
     0 & 1 \\
     0 & -c_1
   \end{array}
\right] x -  \left[\begin{array}{c}
      0 \\
      \alpha
   \end{array}\right] \mbox{sgn} ( \left[
   \begin{array}{cc}
     c_1 &  1 
   \end{array}
\right] x).
\end{equation}
As it is shown in~\cite{galias2007} the trajectories obtained by an explicit Euler discretization exhibit spurious oscillations  which are described by period-2 cycle around the sliding manifold. On Figure~\ref{Fig:GaliasYu2007}, the trajectories obtained by implicit discretization are shown using the implicit Euler discretization with $h=1$, $h=0.3$, $h=0.1$ and  $h=0.01$ and with $c_1=1$ and $\alpha=1$. As it has been predicted by theoretical discussions of Section \ref{section3}, the sliding manifold is reached in finite time and without any chattering. Indeed, the matrix $CB = \alpha  = 1$ satisfies the assumptions of Lemma~\ref{lemmaoups}. Note that the algorithm is also very robust in the sense that the simulation can be performed with relatively large time--steps.

\subsection{Example~\ref{Ex:Multi}: Multiple sliding surfaces}

Let us consider the example~\ref{Ex:Multi}. The system can be defined in the form~(\ref{eq:deux}) with
\begin{equation}
 \label{eq:example2-1}
 B= \left[
   \begin{array}{cc}
     1 & 2 \\
     2 & -1
   \end{array}
\right],\quad C= \left[
   \begin{array}{cc}
     1 & 2 \\
     2 & -1
   \end{array}
\right],\quad D = 0,\quad f(x(t),t) = 0
\end{equation}
\begin{figure}[htbp]
 \psfrag{time}[][]{ \Huge time $\sf t$}
 \psfrag{x1}[r][l]{\Huge $\sf x_1(t)$}
 \psfrag{x2}[r][l]{\Huge $\sf x_2(t)$}
 \psfrag{s1}[r][l]{\Huge $\sf s_1(t)$}
 \psfrag{s2}[r][l]{\Huge $\sf s_2(t)$}
 \psfrag{s values}[r][l]{\Huge \sf s values}
 \psfrag{xv}[][]{\Huge \sf state $\sf x_1(t)$ and $ \sf x_2(t)$}
 \centering
 \subfigure[state $x_1(t)$ and $x_2(t)$ versus time]
   {\label{Fig:Multia}\includegraphics[angle=-90,width=0.32\textwidth]{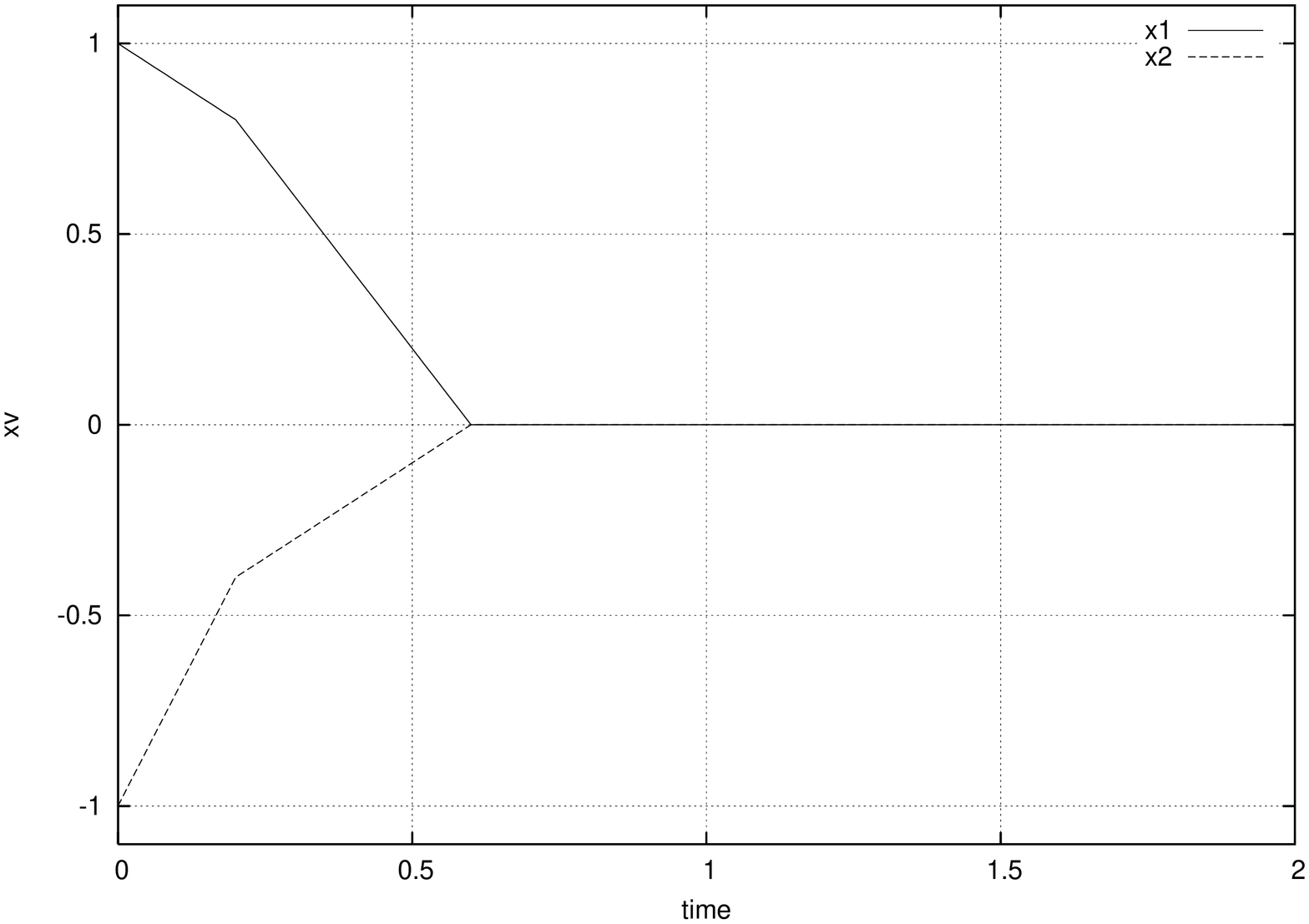}}\hspace{-0.5mm}
 \subfigure[phase portrait  $x_2(t)$ versus $x_1(t)$]
   {\label{Fig:Multib}\includegraphics[angle=-90,width=0.32\textwidth]{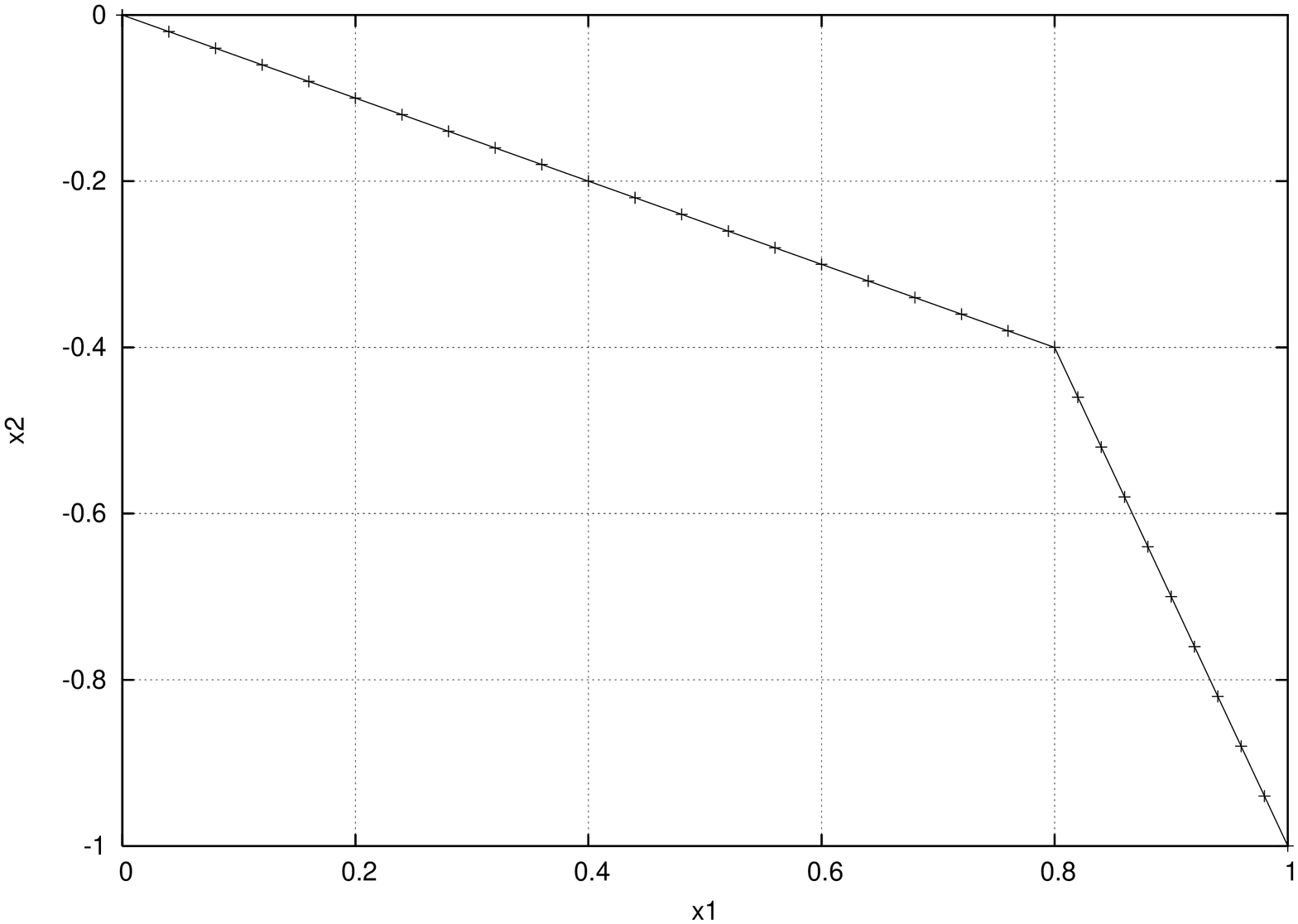}}\hspace{-0.5mm}
   \subfigure[sgn function  $s_1(t)$ and $s_2(t)$]
   {\label{Fig:Multic}\includegraphics[angle=-90,width=0.32\textwidth]{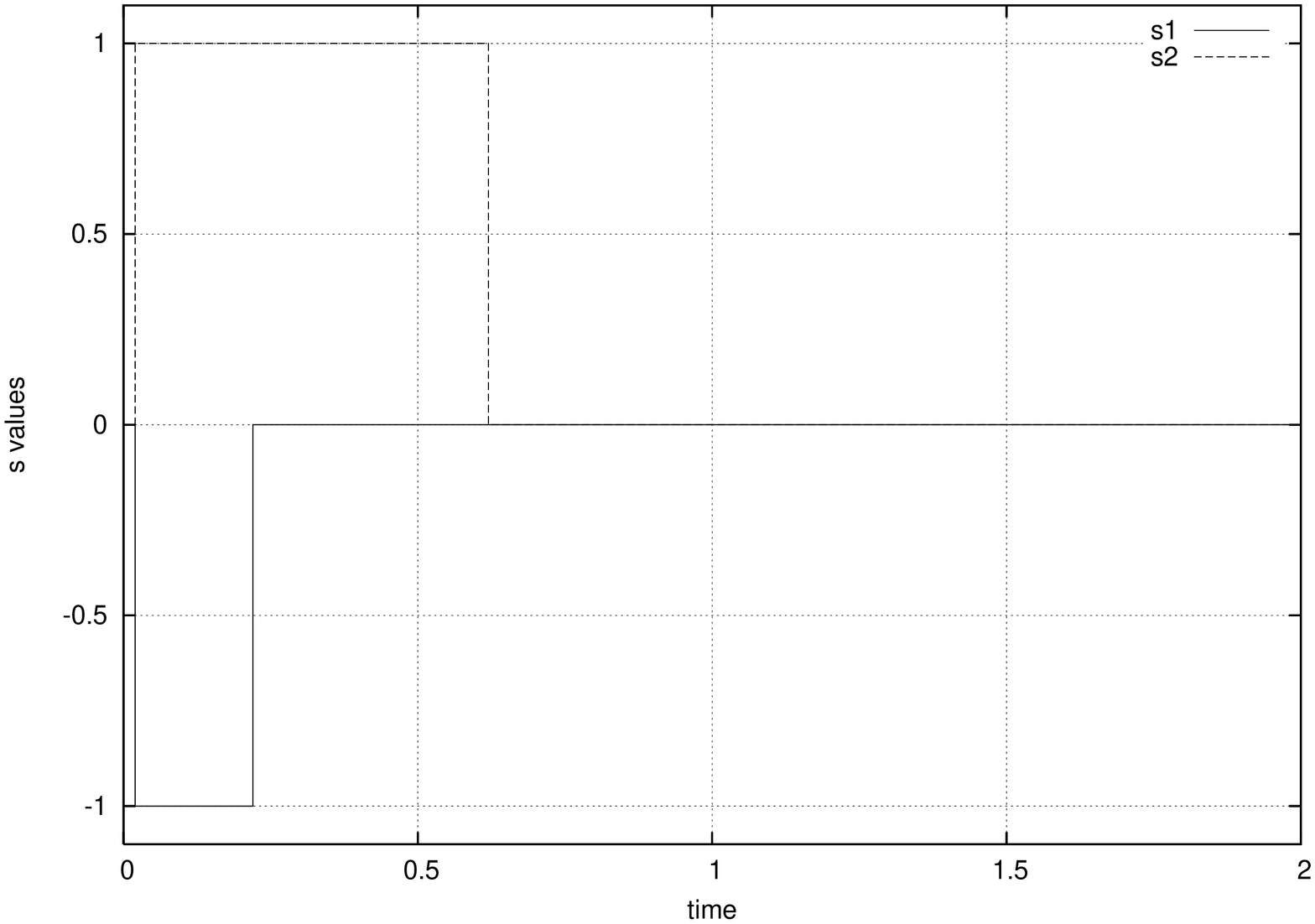}}
 \caption{Multiple Sliding surface. $h = 0.02$, $x(0) = [1.0 , -1.0]^T$} \label{Fig:Multi}
\end{figure}
This example illustrates Lemma~\ref{lemmaoups} since $CB =
\left[\begin{array}{cc}
 5 & 0 \\
 0 & 5 
\end{array}
\right]$. The results displayed on Figure~\ref{Fig:Multi} show that  the system reaches firstly  the sliding surface $2x_2+x_1 = 0$ without any chattering, The system then slides on the surface up to reaching the second sliding surface $2x_1-x_2=0$ and comes to rest at the origin.

\subsection{Extensions to ZOH discretized systems}

\begin{figure}[htbp]
 \psfrag{time}[][]{\Huge time $\sf t$}
 \psfrag{x1}[][]{\Huge $\sf x_1(t)$}
 \psfrag{x2}[][]{\Huge $\sf x_2(t)$}
 \psfrag{s1}[][]{\Huge $\sf s_1(t)$}
 \psfrag{s2}[][]{\Huge $\sf s_2(t)$}
 \psfrag{xv}[][]{\Huge  state}
 \centering
 \subfigure[$h=0.3$. Explicit ZOH]
   {\label{Fig:GaliasYuSISO2008a}\includegraphics[angle=-90,width=0.46\textwidth]{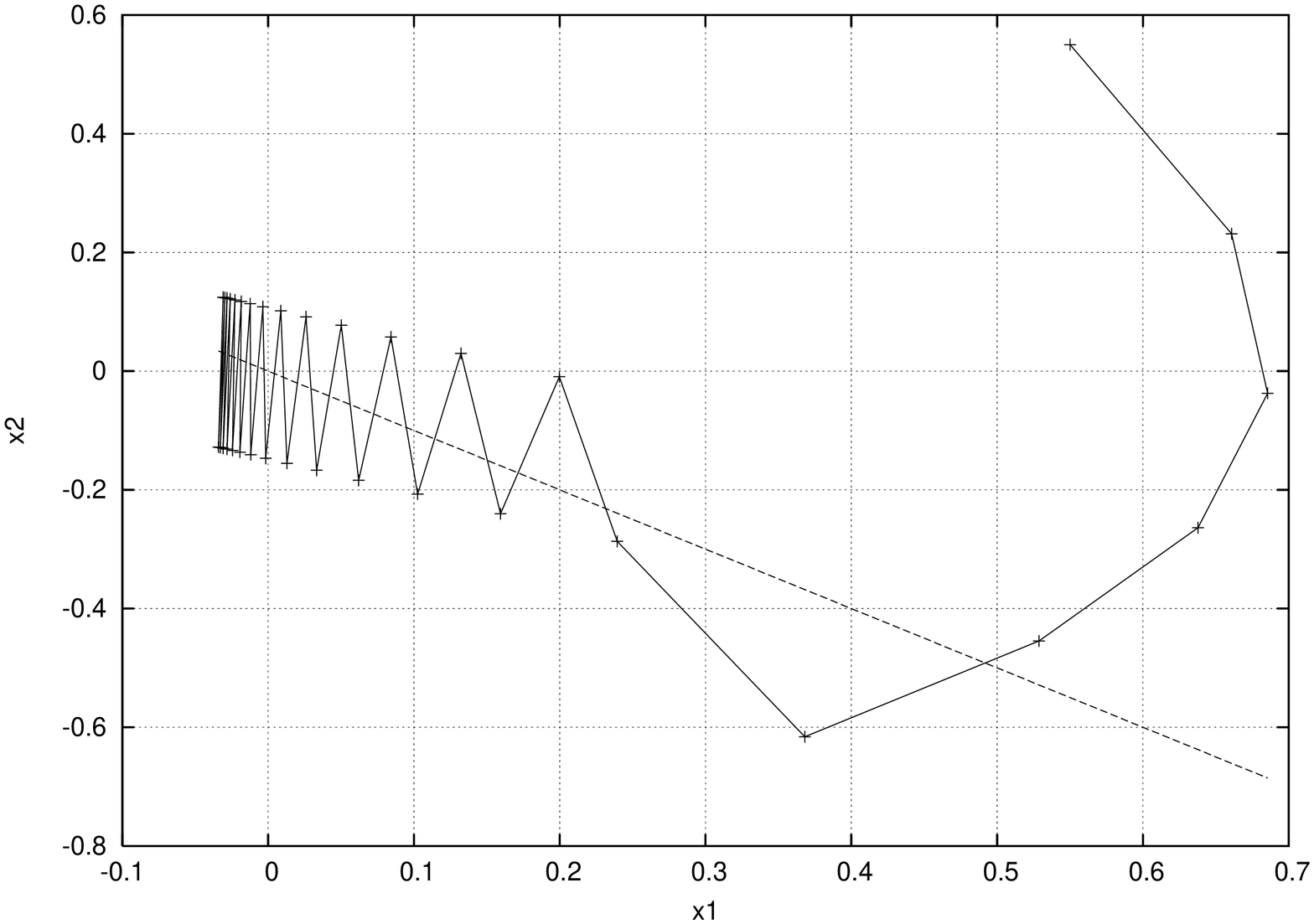}}
 \subfigure[$h=0.3$. Implicit ZOH]
   {\label{Fig:GaliasYuSISO2008b}\includegraphics[angle=-90,width=0.46\textwidth]{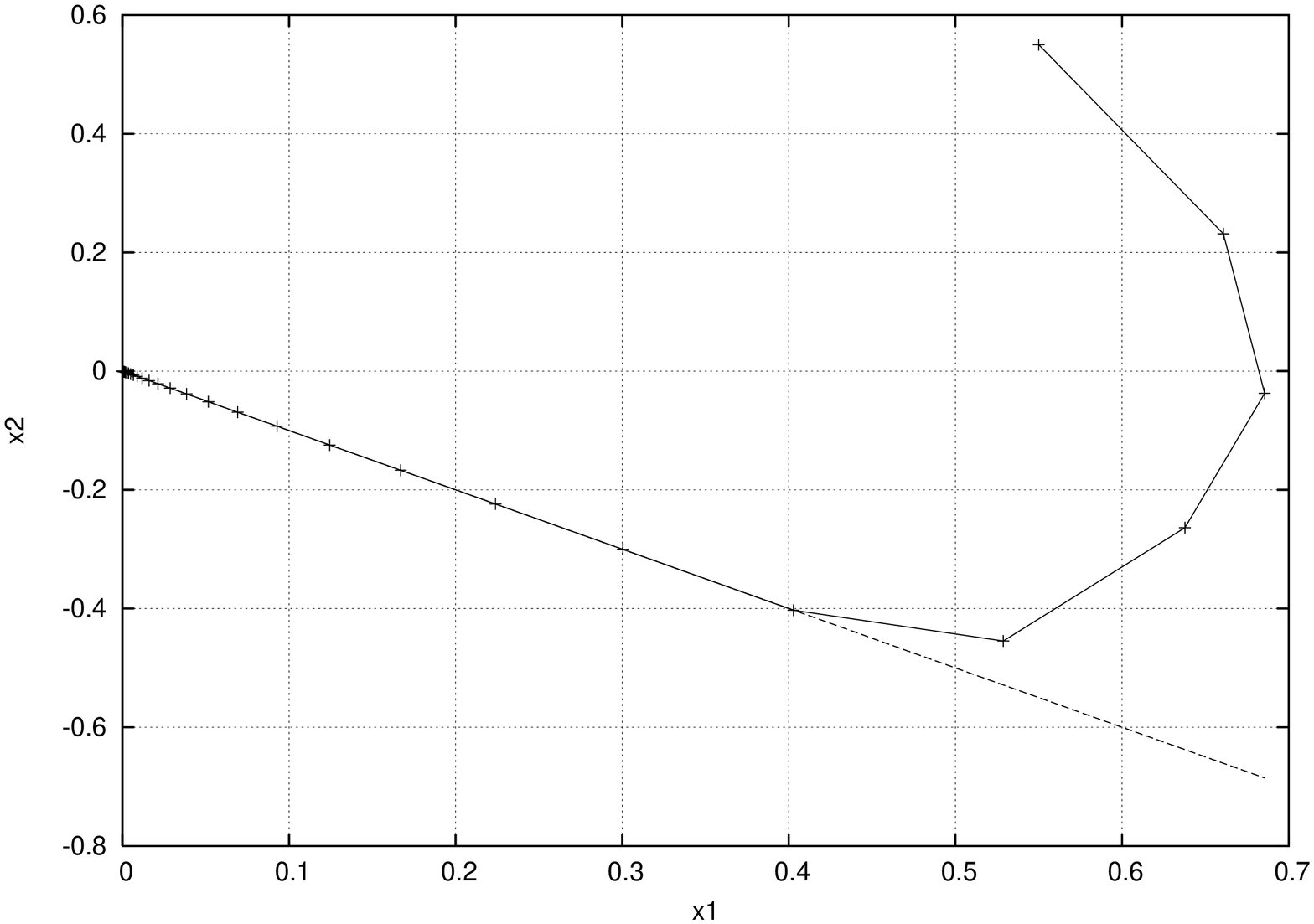}}
 \caption{Equivalent control based SMC,  $a_1=-2$, $a_2=2$, $c_1=1$ and $h=0.3$. $x_0=[0.55,0,55]^T$ State $x_1(t)$ versus $x_2(t)$.}  \label{Fig:GaliasYuSISO2008}
\end{figure}

\begin{figure}[htbp]
 \psfrag{time}[][]{\Huge  time $\sf t$}
 \psfrag{x1}[][]{\Huge $\sf x_1(t)$}
 \psfrag{x2}[][]{\Huge $\sf x_2(t)$}
 \psfrag{s1}[][]{\Huge $\sf s_1(t)$}
 \psfrag{s2}[][]{\Huge $\sf s_2(t)$}
 \psfrag{xv}[][]{state}
 \centering
  \subfigure[Explicit implementation]
 {\includegraphics[angle=-90,width=0.46\textwidth]{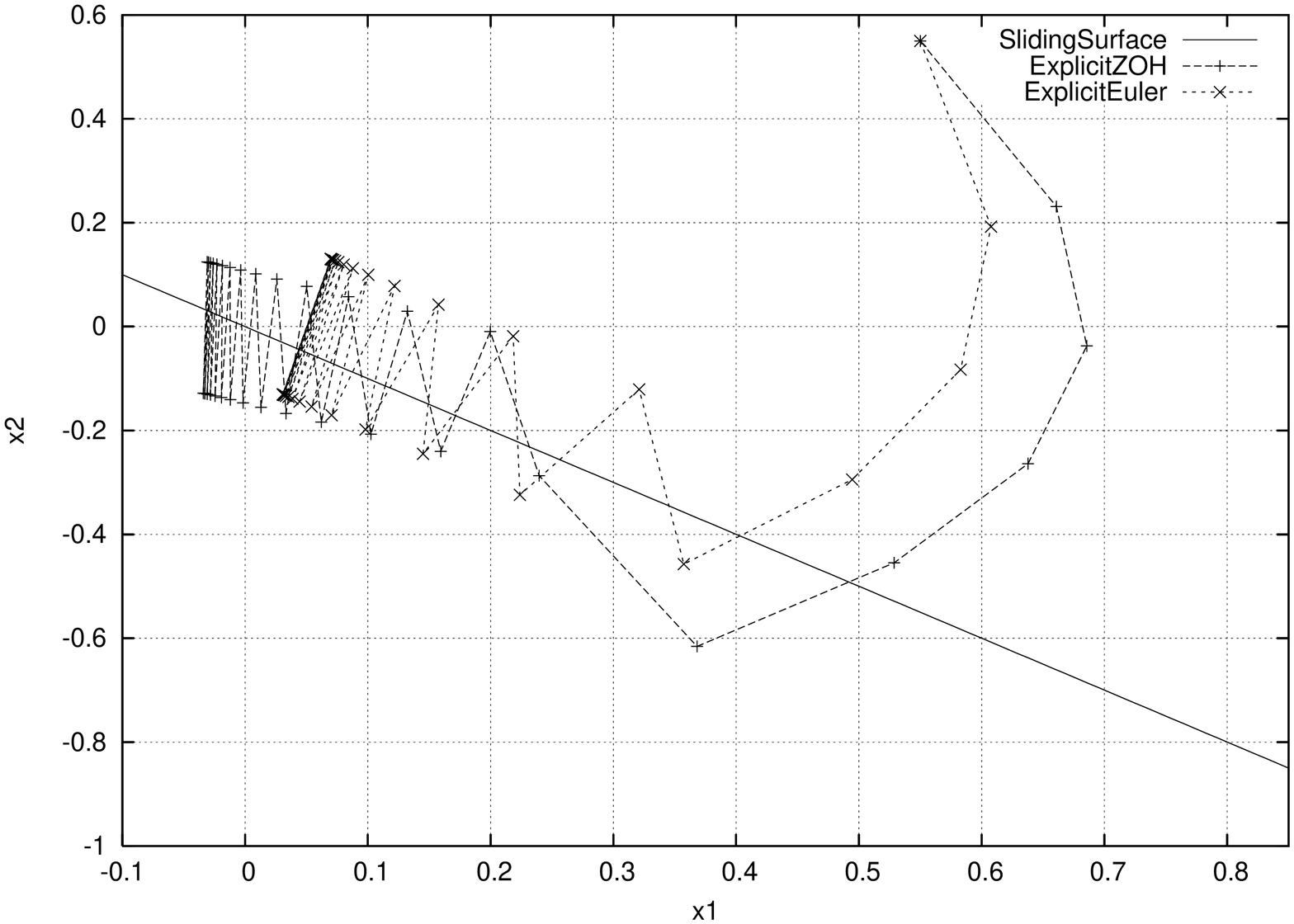}} 
 \subfigure[Implicit implementation]
 {\includegraphics[angle=-90,width=0.46\textwidth]{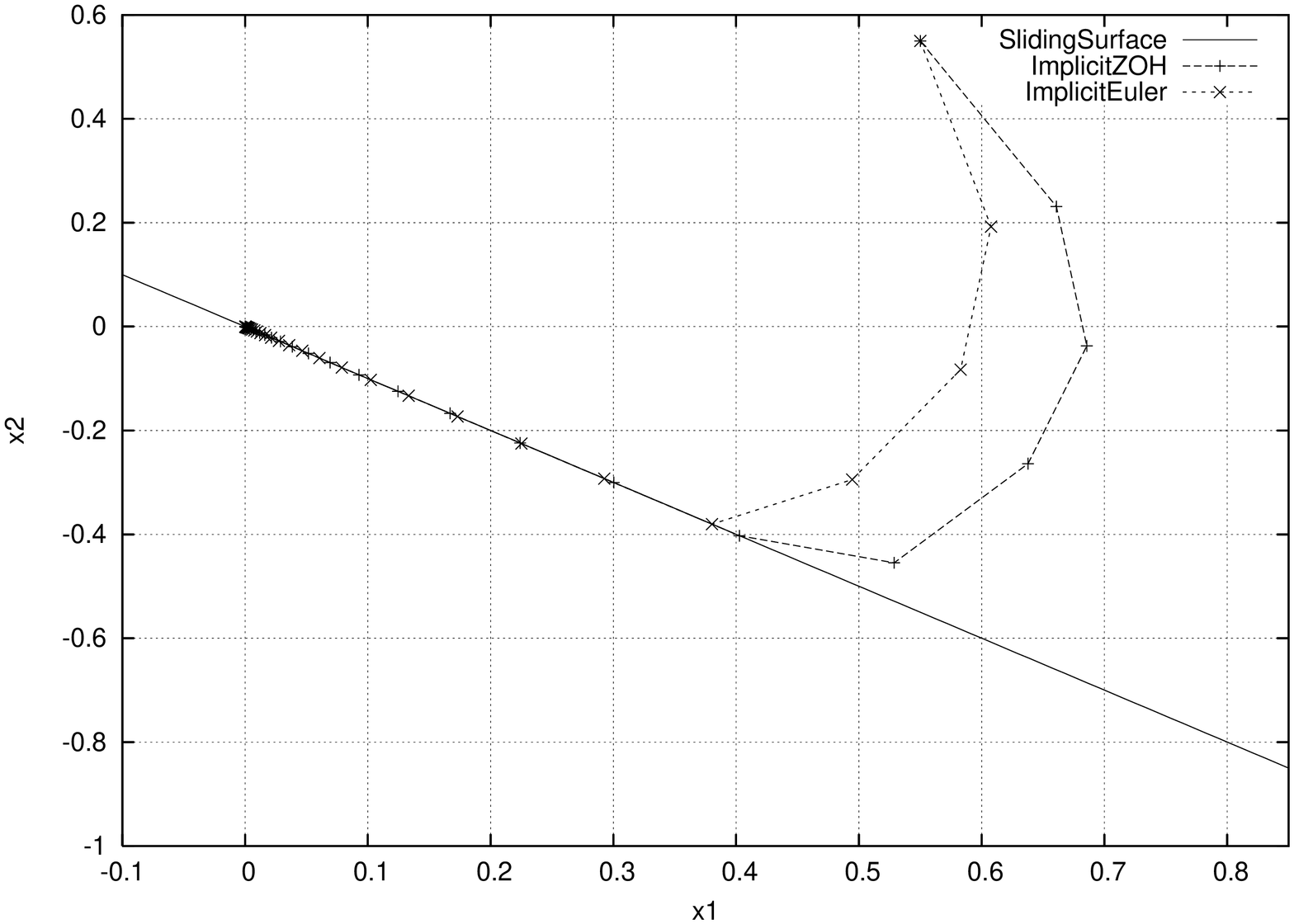}} 
\caption{Comparison of Euler and ZOH discretizations of ECB-SMC system,  $a_1=-2$, $a_2=2$, $c_1=1$ and $h=0.3$. $x_0=[0.55,0,55]^T$ State $x_1(t)$ versus $x_2(t)$.}
  \label{Fig:GaliasYuSISO2008-Compa}
\end{figure}

\begin{figure}[htbp]
 \psfrag{x1}[][]{\Huge $\sf x_1(t)$}
 \psfrag{x2}[][]{\Huge $\sf x_2(t)$}
  \psfrag{x3}[][]{\Huge $\sf x_3(t)$}
 \centering
 \subfigure[$h=0.3$. Explicit ZOH]
   {\label{Fig:GaliasYuMIMO2008a}\includegraphics[angle=-90,width=0.46\textwidth]{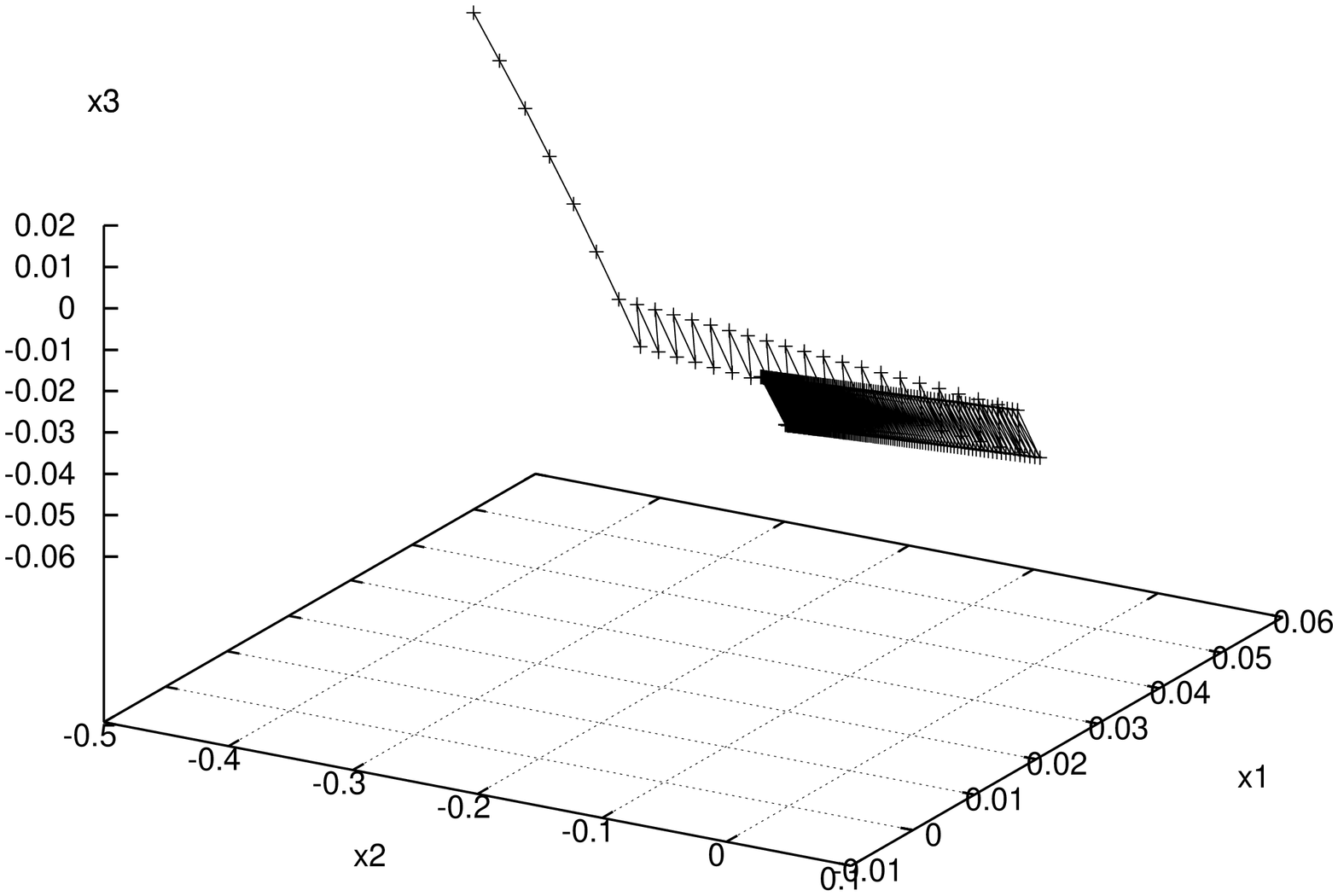}}
 \subfigure[$h=0.3$. Implicit ZOH]
   {\label{Fig:GaliasYuMIMO2008b}\includegraphics[angle=-90,width=0.46\textwidth]{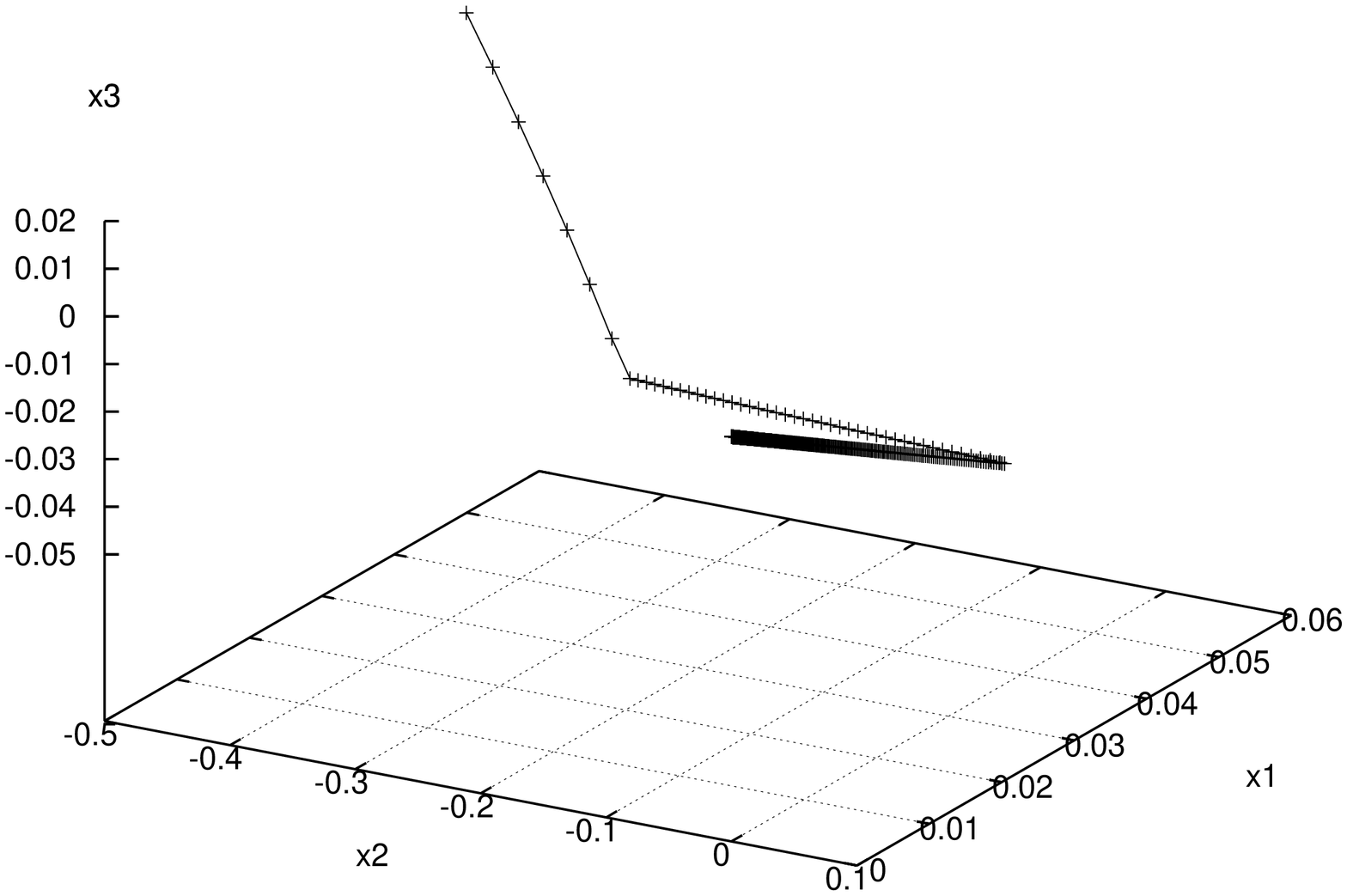}}
 \caption{Equivalent control based SMC,  Trajectory with initial point $x_0=[0.05,-0.5,0.02]$}  \label{Fig:GaliasYuMIMO2008}
\end{figure}

The extension to ZOH discretized systems is illustrated on a first example taken from \cite{Galias2008a}. In the notation of Section~\ref{subsecZOH}, the LTI system with an ECB-SMC controller is defined by the following data,
\begin{equation}
  \label{eq:GaliasYuSISO2008}
  F = \left[
    \begin{array}{cc}
      0 & 1 \\
      -a_1 & -a_2
    \end{array}
\right],\quad G =\left[
    \begin{array}{c}
      0  \\
      1
    \end{array}
\right],\quad  C =\left[
    \begin{array}{cc}
      c_1 & 1 \\
    \end{array}
\right].
\end{equation}
Starting from the initial data, $x_0=[0.55,0,55]^T$, Galias and Yu \cite{Galias2008a} have shown that the Explicit ZOH discretization of the system with $a_1=-2$, $a_2=2$, $c_1=1$ and $h=0.3$ exhibits a period--2 orbit. The results are reproduced on Figure~\ref{Fig:GaliasYuSISO2008a}. On Figure~\ref{Fig:GaliasYuSISO2008b}, the Implicit ZOH discretization as proposed in Section~\ref{subsecZOH} is free of chattering. On Figure~\ref{Fig:GaliasYuSISO2008-Compa}, a comparison is given between the ZOH and the Euler discretization of the vector field $f(\cdot,\cdot)$.
Another example taken from \cite{yu2008} in the MIMO case is given by the following parameters,
\begin{equation}
  \label{eq:GaliasYuMIMO2008}
  F = \left[
    \begin{array}{ccc}
      0  & 0 & 1 \\
      1 & 1 & 1 \\
      -1 & -3 & 1 \\
    \end{array}
\right],\quad G =\left[
    \begin{array}{cc}
      0 & 0 \\
      1 & 0 \\
      0 & 1
    \end{array}
\right],\quad  C =\left[
    \begin{array}{ccc}
      1 & 0 & 1 \\
      0 & 1&1 \\
    \end{array}
\right].
\end{equation}
Similar results are depicted on Figure~\ref{Fig:GaliasYuMIMO2008}.

\subsection{Lyapunov-based robust control}

\begin{figure}[htbp]
 \psfrag{t}[][]{\Huge\sf  time $\sf t$}
 \psfrag{x1}[][]{\Huge \sf state $\sf x_1(t)$}
 \psfrag{tau}[][]{\Huge \sf control $\sf u(t)$}
 \psfrag{xv}[][]{state}
 \centering
 \subfigure[State $x_1(t)$ vs. time. $h=0.1$. Implicit Euler]
   {\label{Fig:LyapunovBased1a}\includegraphics[angle=-90,width=0.40\textwidth]{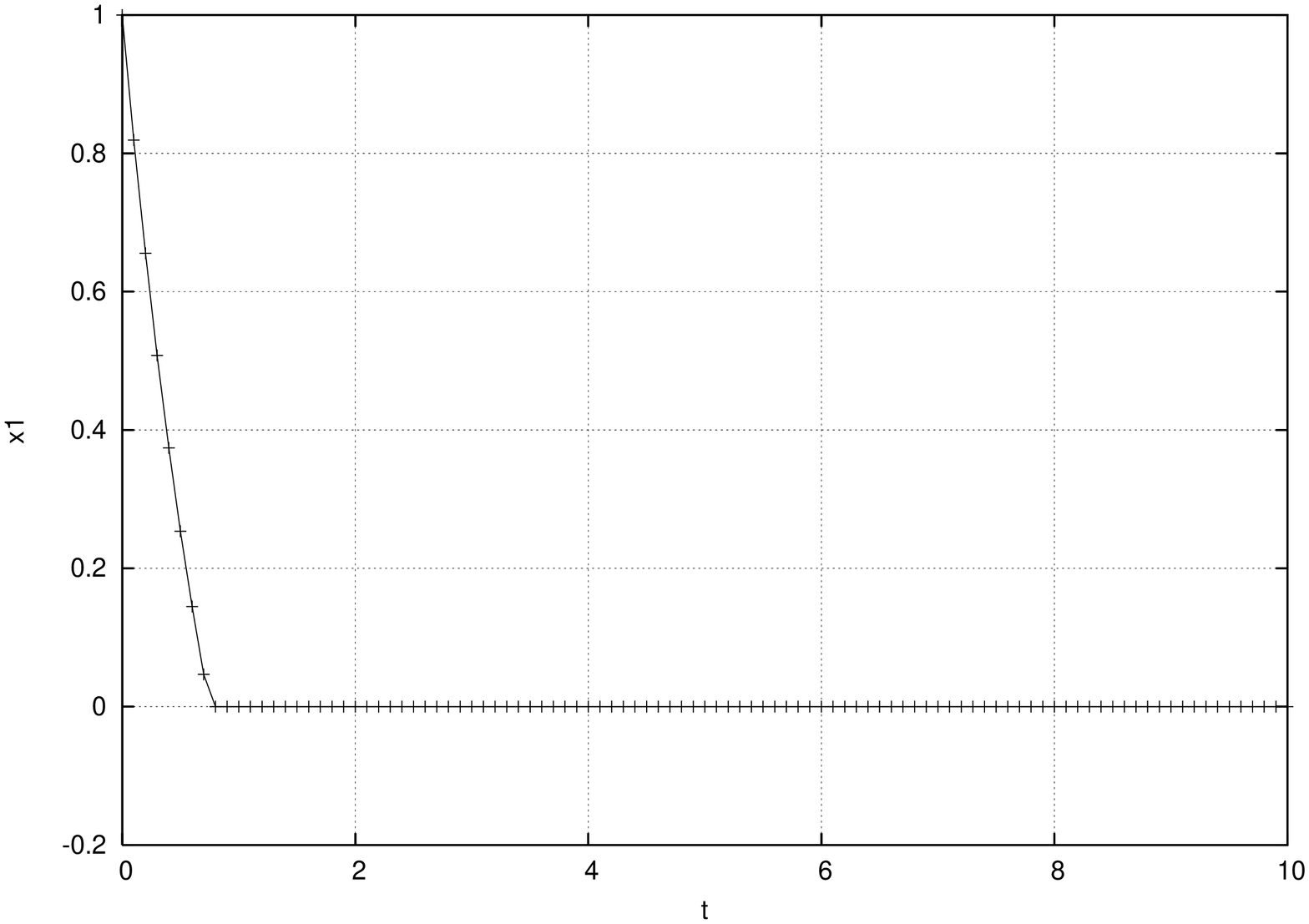}}
 \subfigure[State $x_1(t)$ vs. time. $h=0.1$. Explicit Euler]
   {\label{Fig:LyapunovBased1b}\includegraphics[angle=-90,width=0.40\textwidth]{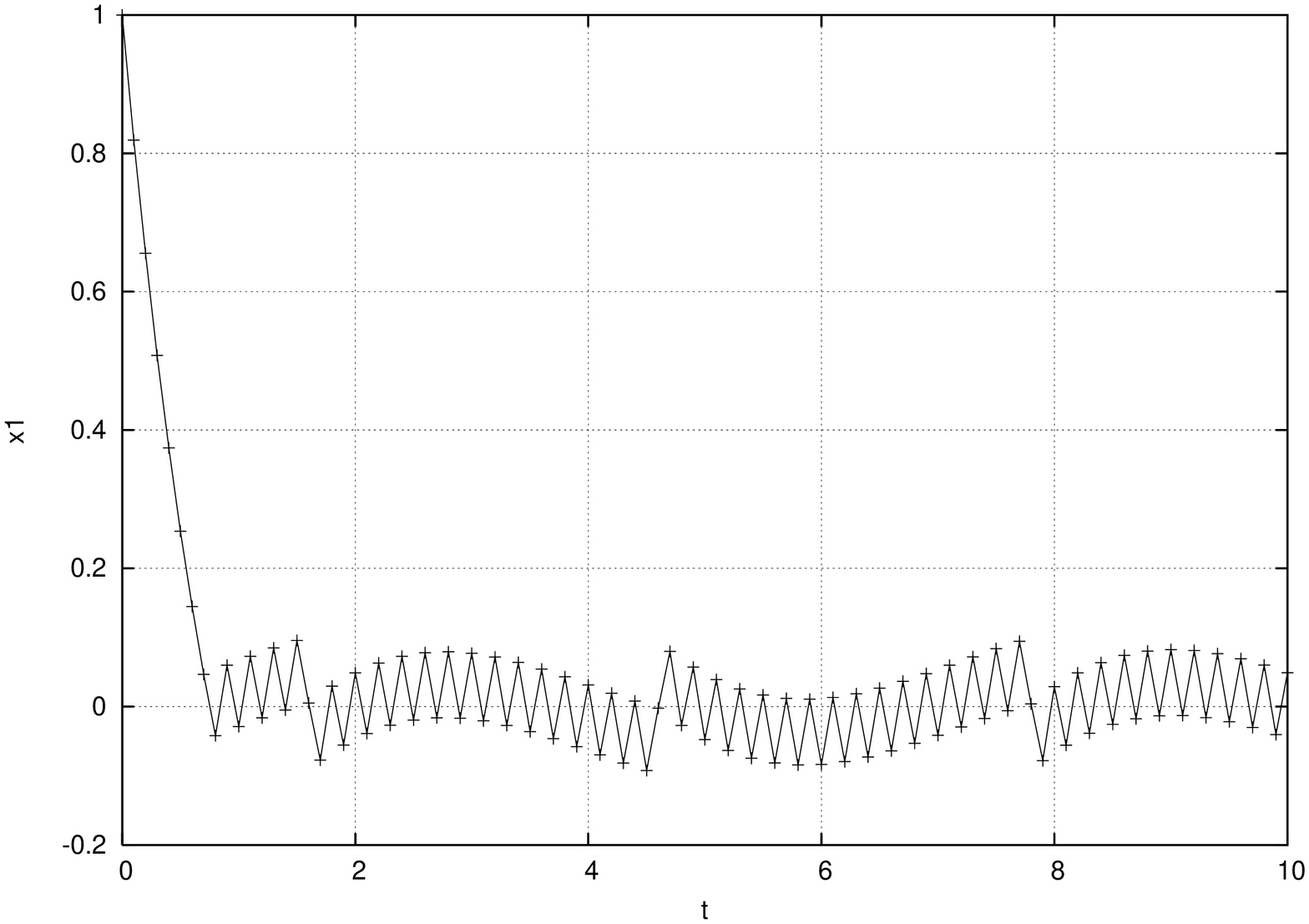}}
 \subfigure[Control $u(t)$ vs. time. $h=0.1$. Implicit Euler]
   {\label{Fig:LyapunovBased1c}\includegraphics[angle=-90,width=0.40\textwidth]{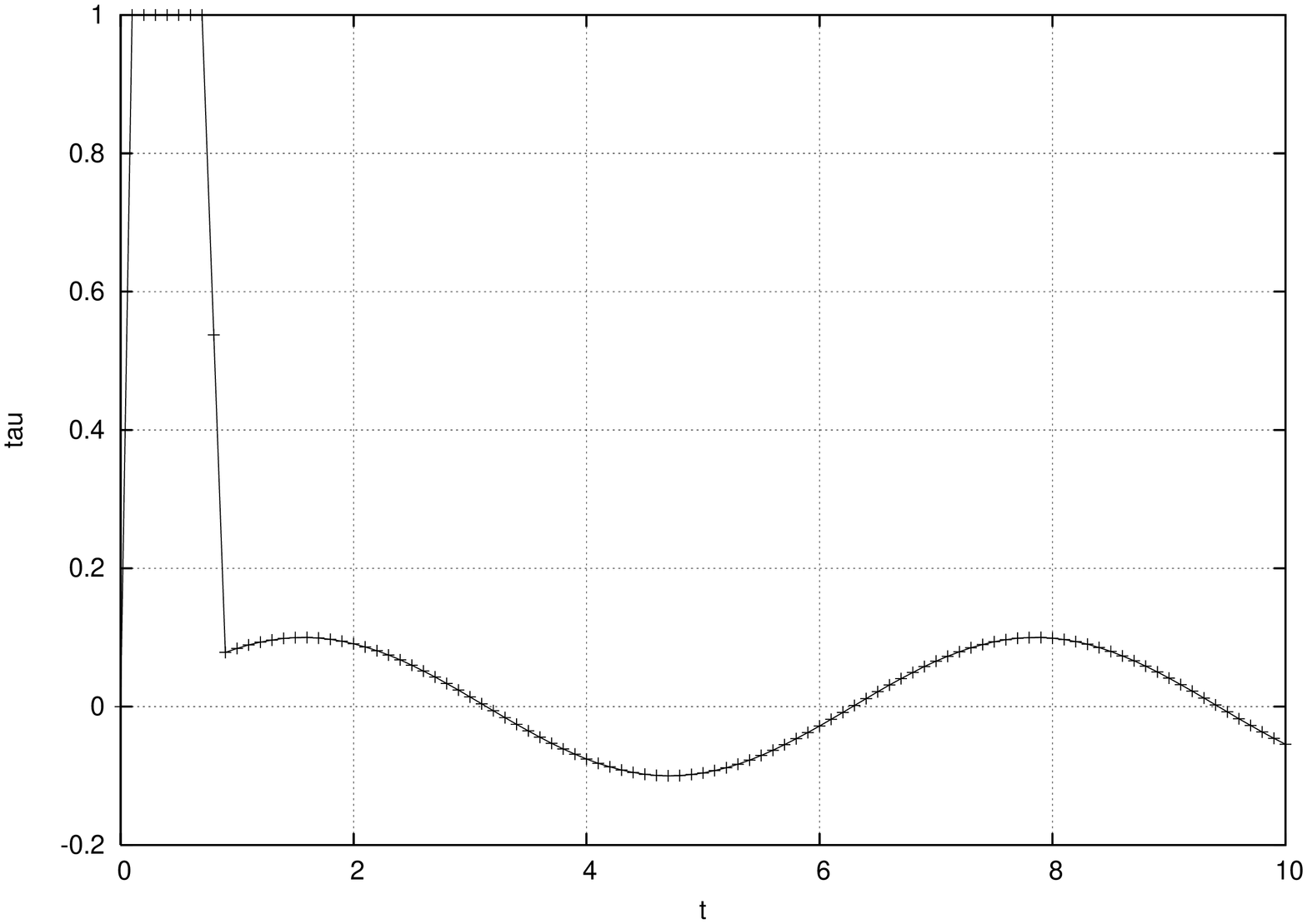}}
 \subfigure[Control $u(t)$ vs. time. $h=0.1$. Explicit Euler]
   {\label{Fig:LyapunovBased1d}\includegraphics[angle=-90,width=0.40\textwidth]{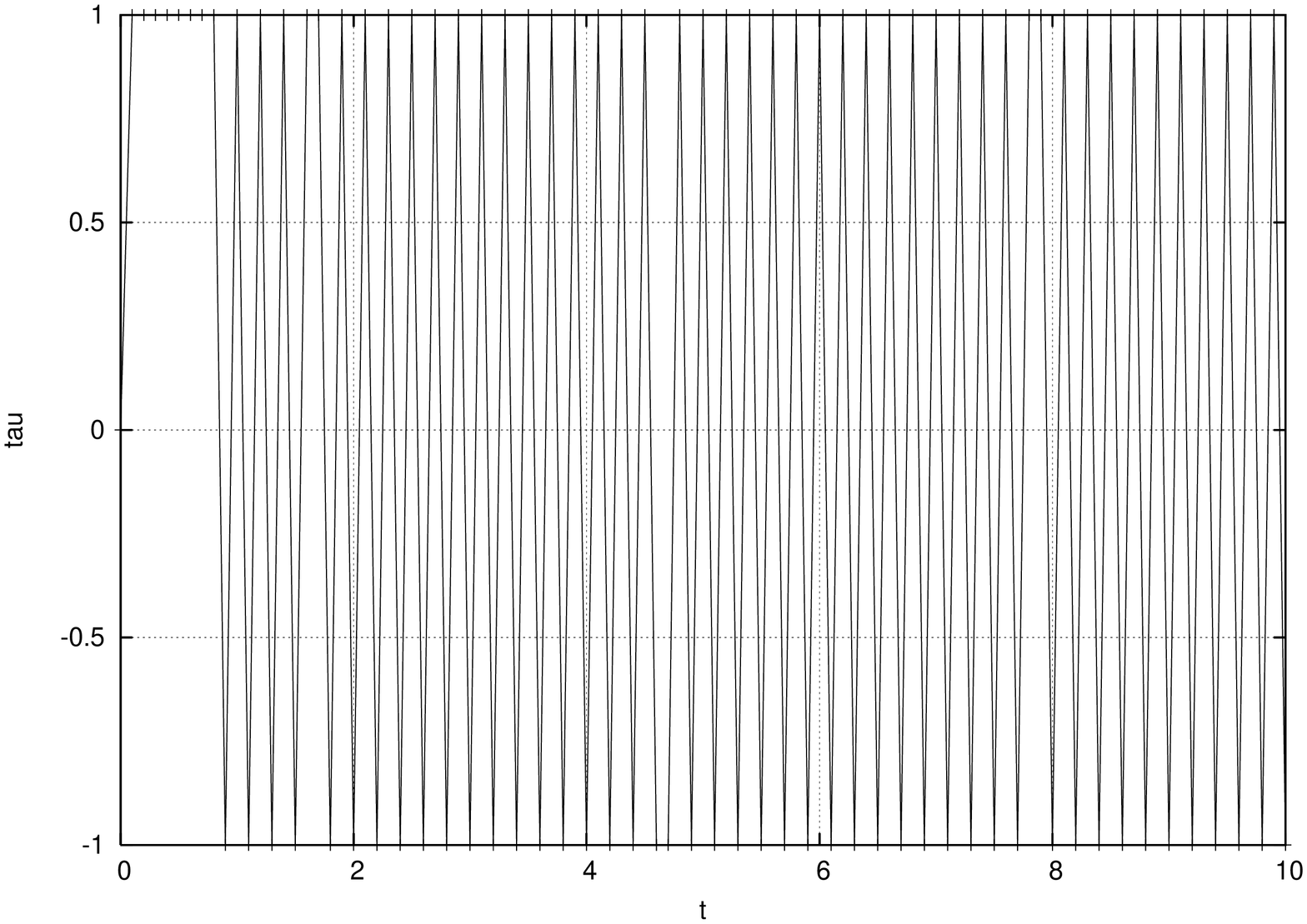}}
 \caption{Lyapunov-based discontinuous robust control. $h=0.1$ $\alpha=0.1$}  \label{Fig:LyapunovBased1}
\end{figure}

We propose in this section to give an numerical example which fits with the example~(\ref{cor1}) of a Lyapunov-based discontinuous robust control. Let us consider the following system
\begin{equation}
  \label{eq:lyapu1}
  \dot x(t) = - x(t) - u(t) + \gamma(t)
\end{equation}
with $\gamma(t) = \alpha \sin(t)$ and $u(t) = \mbox{sgn}(x(t))$. Is is obvious that, as expected, the implicit method yields a smooth stabilization at $x=0$ whereas the explicit Euler has significant chattering. Figure~\ref{Fig:LyapunovBased1c} illustrates the fact that the controller varies inside the multivalued part of the $\mbox{sgn}$ function in order to assure the existence of an equilibrium point.

\subsection{The  Filippov example}

\begin{example}Let us consider now the well known Filippov example which  can be defined in the form~(\ref{eq:deux}) with
\begin{equation}
  \label{eq:example2-1bis}
  B= \left[
    \begin{array}{cc}
      1 & -2 \\
      2 & 1
    \end{array}
\right],\quad C= \left[
    \begin{array}{cc}
      1 & 0 \\
      0 & 1
    \end{array}
\right],\quad D = 0,\quad f(x(t),t) = 0.
\end{equation} The  trajectories may slide on the codimension 2 surface given by $Cx=0$, that $x = 0$. 
\label{Ex:Filippov}
\end{example}

\begin{figure}[htbp]
  \psfrag{time}[][]{ time $\sf t$}
  \psfrag{x1}[][]{$\sf x_1(t)$}
  \psfrag{x2}[][]{$\sf x_2(t)$}
  \psfrag{s1}[][]{$\sf s_1(t)$}
  \psfrag{s2}[][]{$\sf s_2(t)$}
  \psfrag{xv}[][]{state}
  \centering
  \subfigure[state $x_1(t)$ and $x_2(t)$ versus time]
    {\label{Fig:Multi2a}\includegraphics[angle=-90,width=0.49\textwidth]{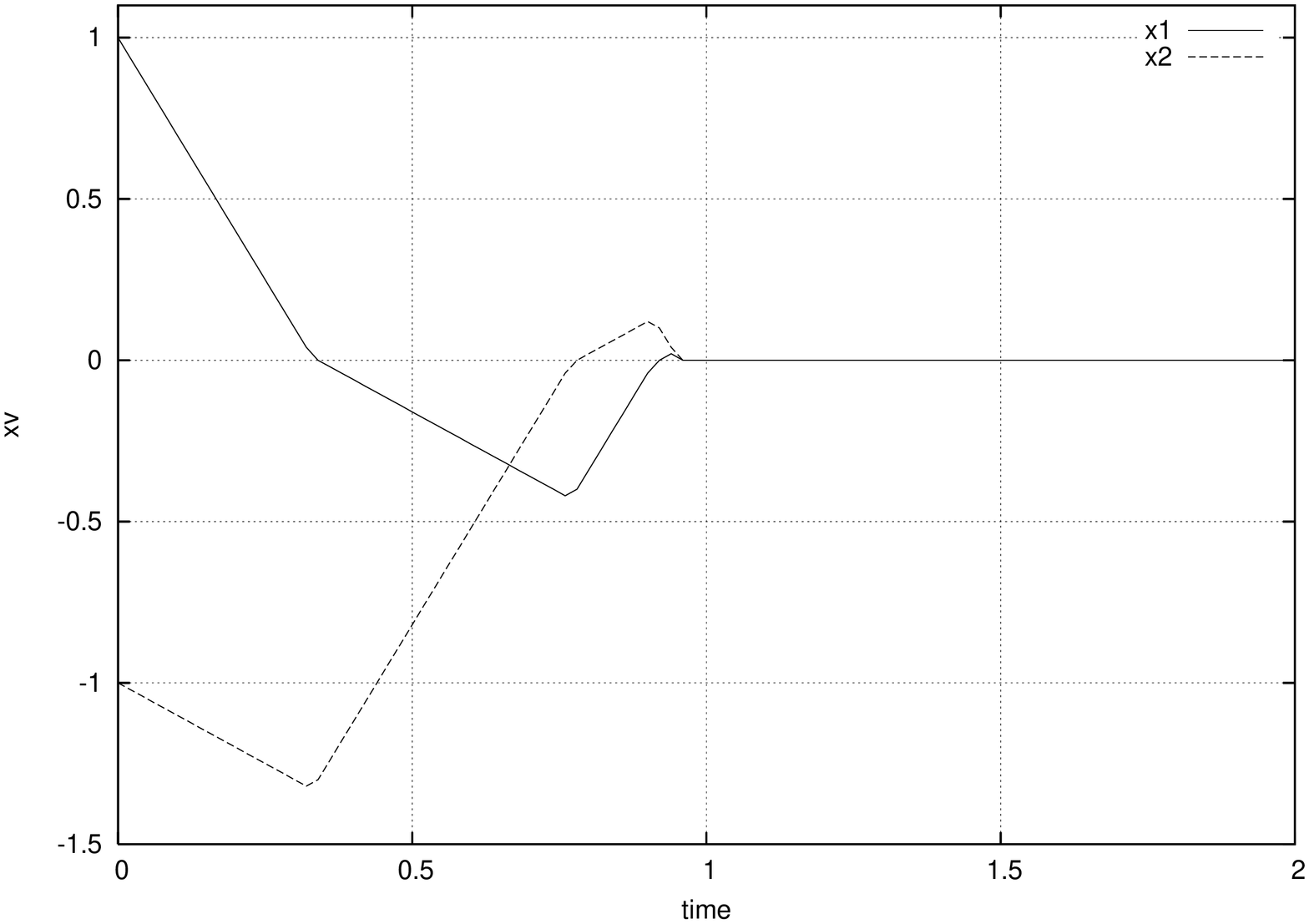}}
  \subfigure[phase portrait  $x_2(t)$ versus $x_1(t)$]
    {\label{Fig:Multi2b}\includegraphics[angle=-90,width=0.49\textwidth]{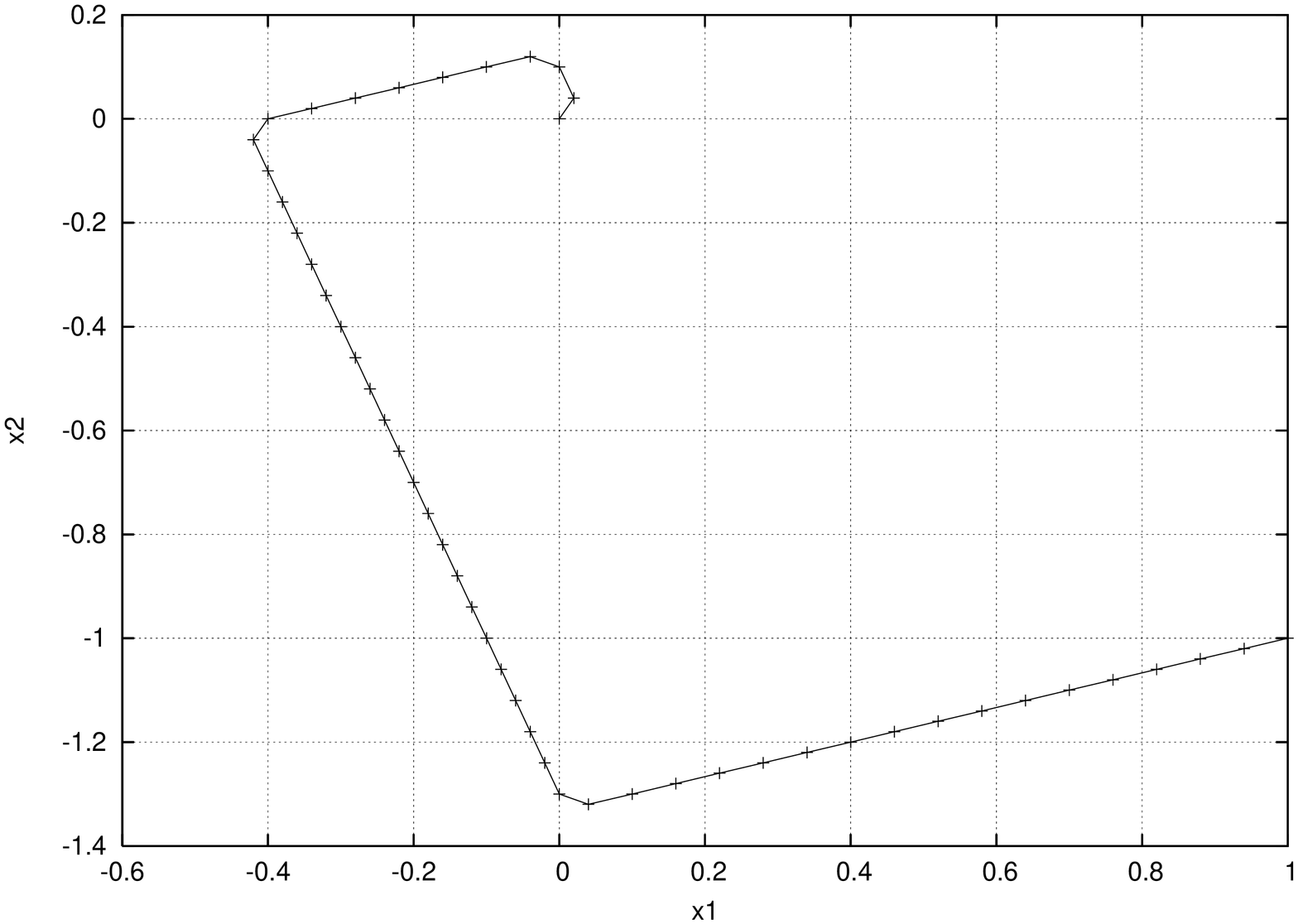}}
    \subfigure[sgn function  $s_1(t)$ and $s_2(t)$]
    {\label{Fig:Multi2c}\includegraphics[angle=-90,width=0.49\textwidth]{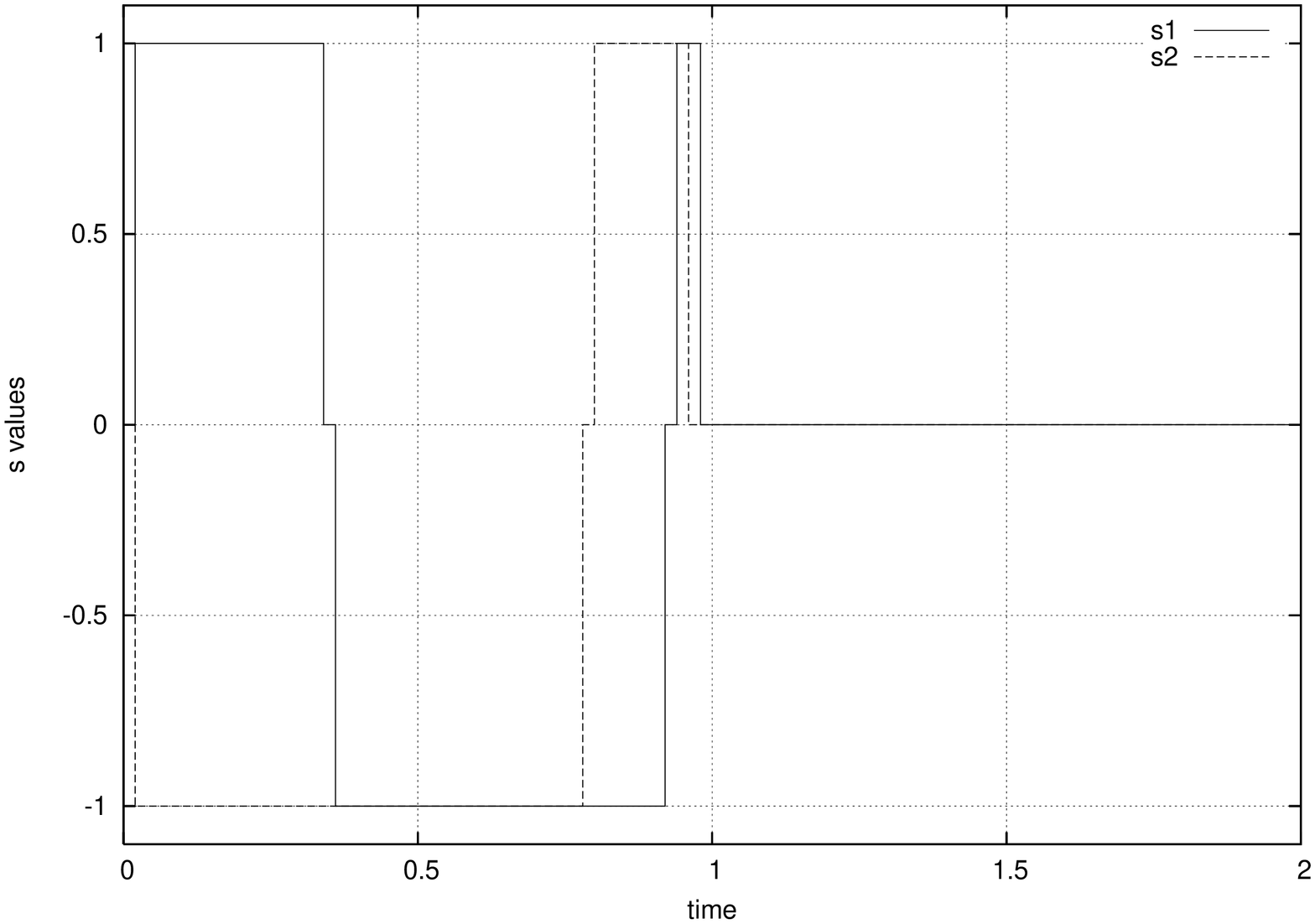}}
  \label{Fig:Multi2}
  \caption{Multiple Sliding surface. Filippov Example. $h = 0.002$, $x(0) = [1.0 , -1.0]^T$}
\end{figure}
The results displayed on Figure~\ref{Fig:Multi2} show that the system reaches the origin without any chattering. The sufficient conditions of the Lemma~\ref{lemmaoups} are not satisfied when seems to indicated that these conditions has to be improved.

\subsection{Example~\ref{Ex:observer}: Observer based SMC}

\begin{figure}[htbp]
  \centering
  \includegraphics[angle=-90,width=0.6\textwidth]{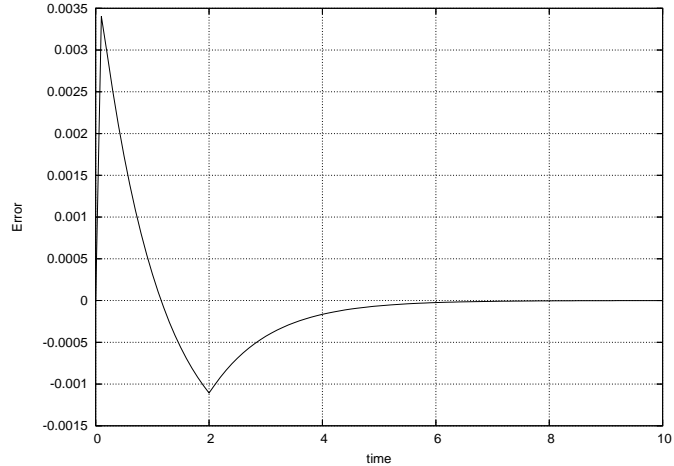}
  \caption{Observer based SMC: Error $e(t)$.  $k=1$ and $\tau = 0.001$. $h = 0.1$  $\theta=1, \gamma=1$  }\label{Fig:SMC-Error}
\end{figure}

\begin{figure}[htbp]
  \centering
  \subfigure[control $u(t) = \mbox{Sgn}(Cx(t))$]
    {\label{Fig:SMC-Controla}\includegraphics[angle=-90,width=0.49\textwidth]{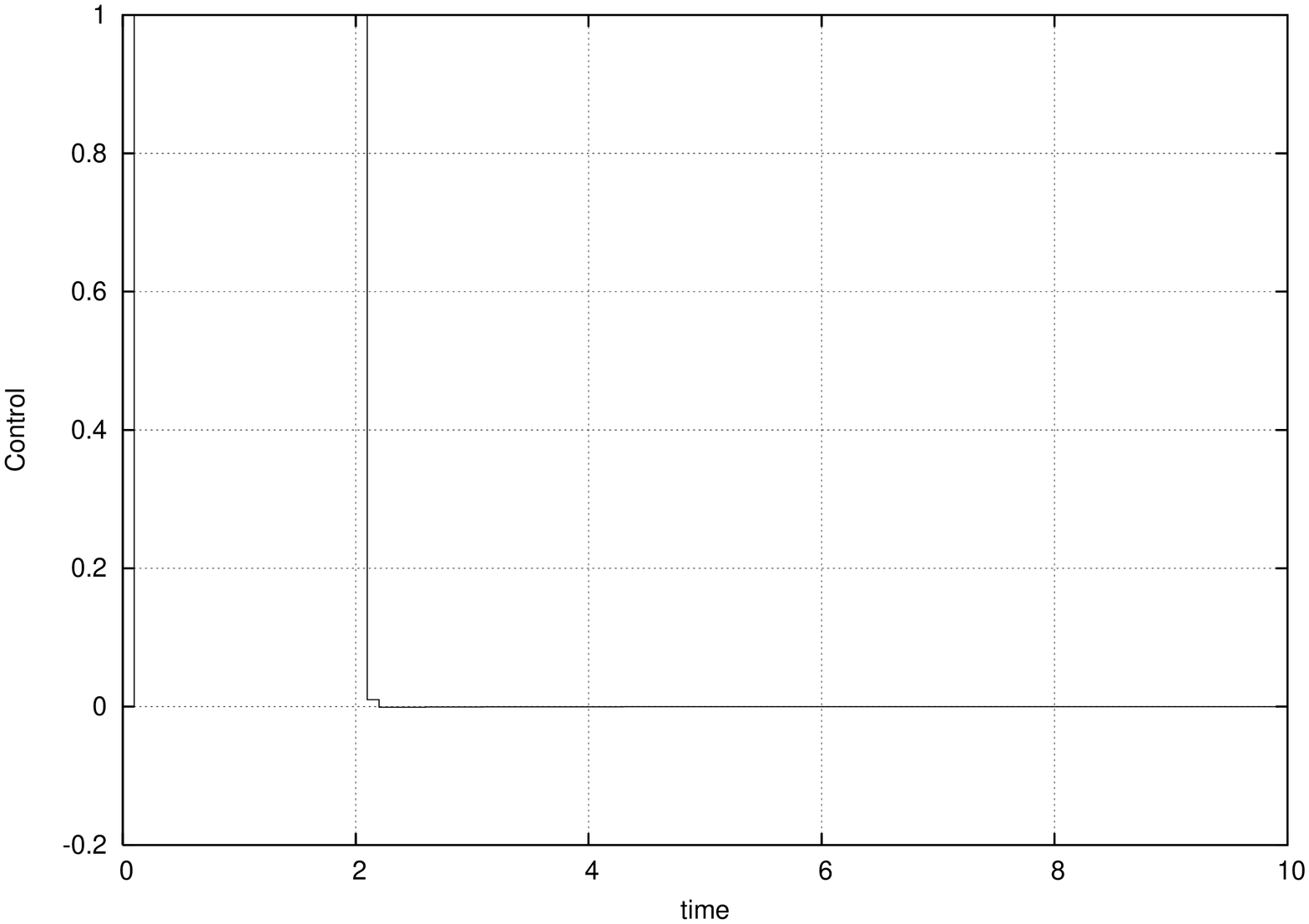}}
    \subfigure[Zoom on control $u(t) = \mbox{Sgn}(Cx(t))$]{\label{Fig:SMC-Controlb}\includegraphics[angle=-90,width=0.49\textwidth]{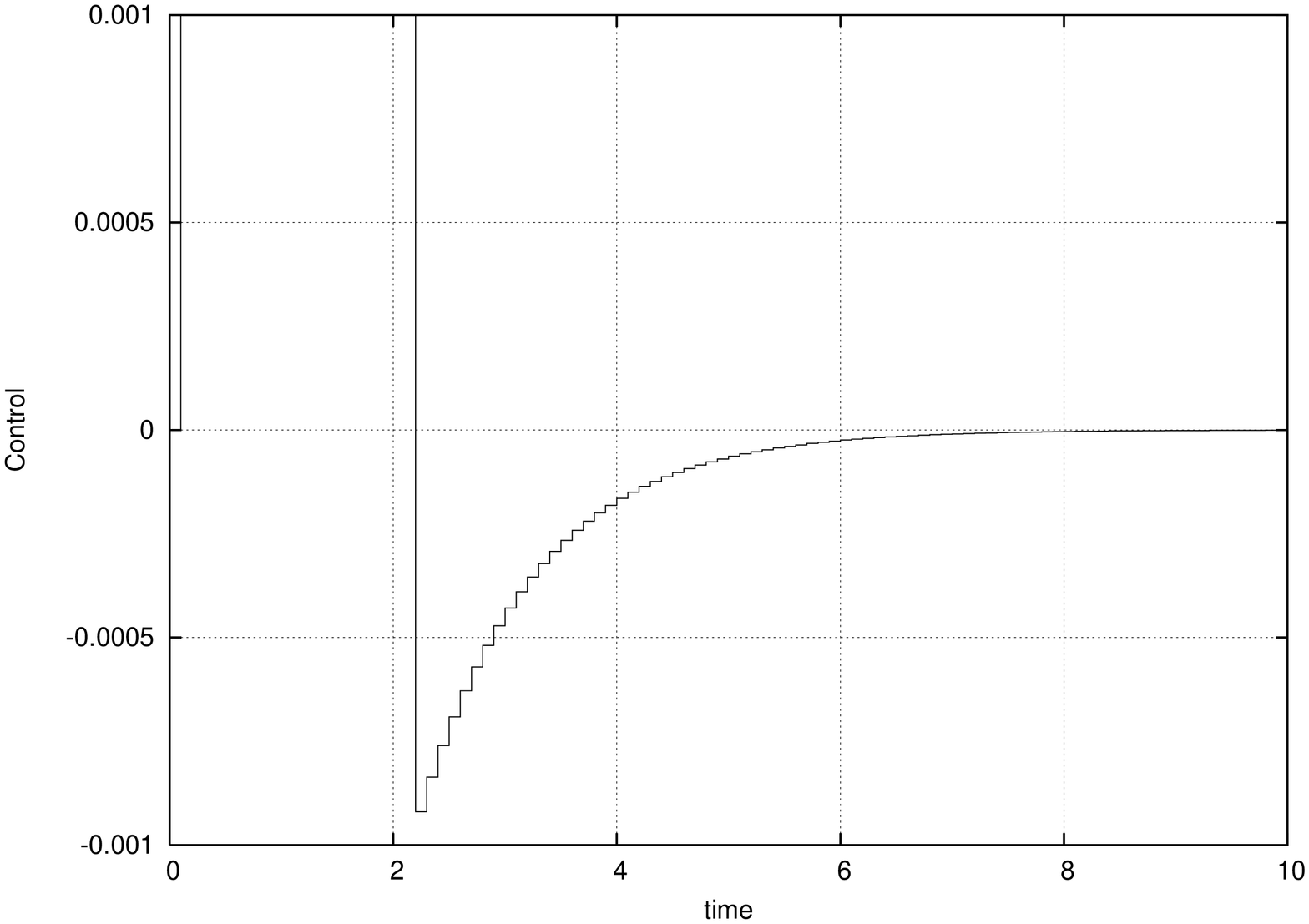}}
  \caption{Observer based SMC: Control.  $k=1$ and $\tau = 0.001$. $h = 0.1$  $\theta=1, \gamma=1$} \label{Fig:SMC-Control}
\end{figure}

Let us illustrate the performance of our implementation on the observer based SMC described by the Example \ref{Ex:observer}. The dynamics is given by \eqref{observer} with $k=1$ and $\tau = 0.001$. The initial conditions are chosen as $[2.0, 0, 0 , 0]^T$.  The numerical parameters are given by $h = 0.1$ that is a sampling of $10Hz$ and $\theta=1, \gamma=1$. On Figure~\ref{Fig:SMC-Error}, the error between the reference command and the observer state is given. On Figure~\ref{Fig:SMC-Control}, we can observe the behavior of the command without any chattering.

On the three previous examples, one sees a very accurate and smooth stabilization on the sliding surface, even for values of $h$ not so small.

\begin{figure}[htbp]
  \centering
  \subfigure[state $x_1(t)$]{\label{Fig:SMC-EulerExplicita}\includegraphics[angle=-90,width=0.49\textwidth]{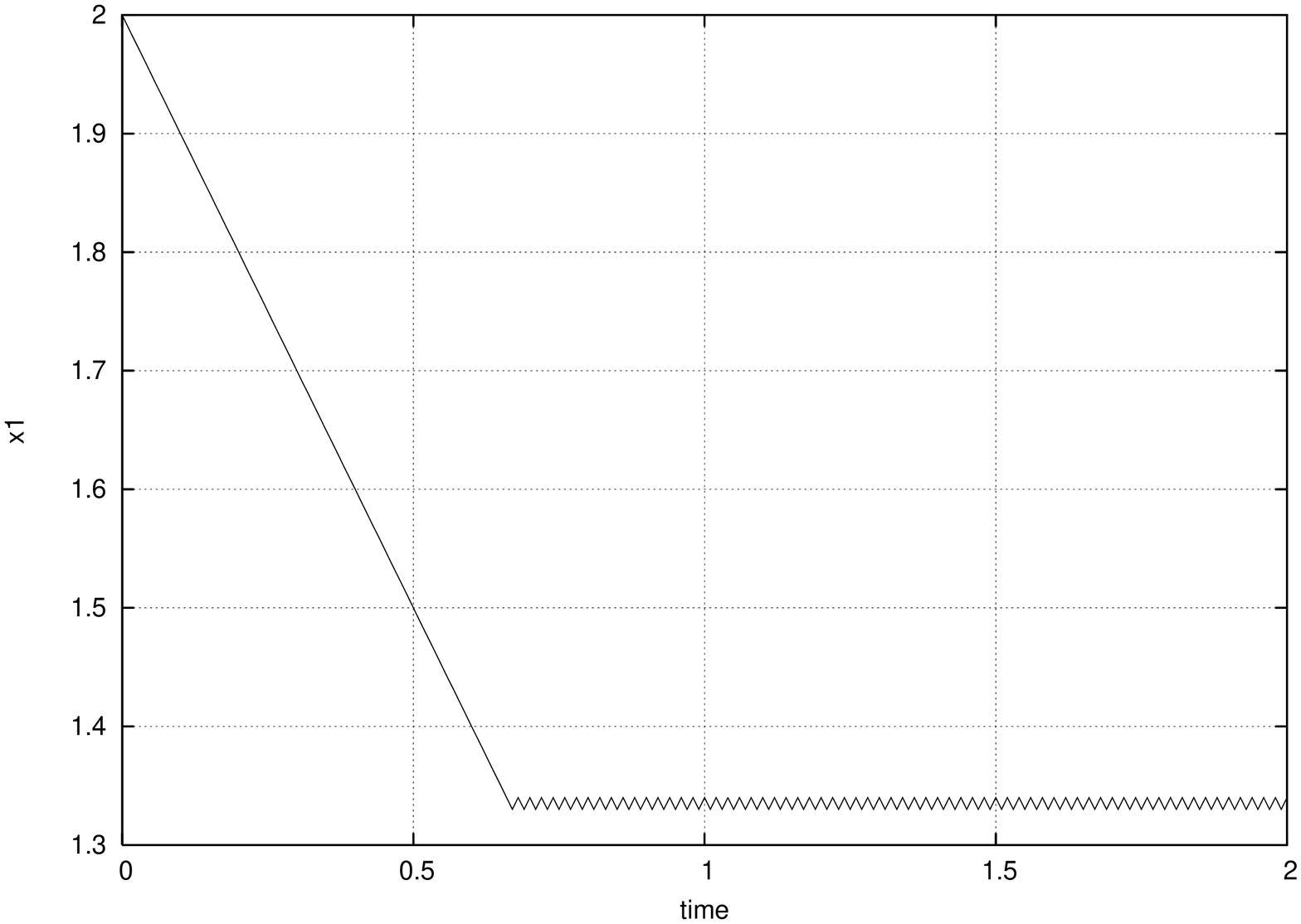}}
  \subfigure[control $u(t) = \mbox{Sgn}(Cx(t))$]{\label{Fig:SMC-EulerExplicitb}\includegraphics[angle=-90,width=0.49\textwidth]{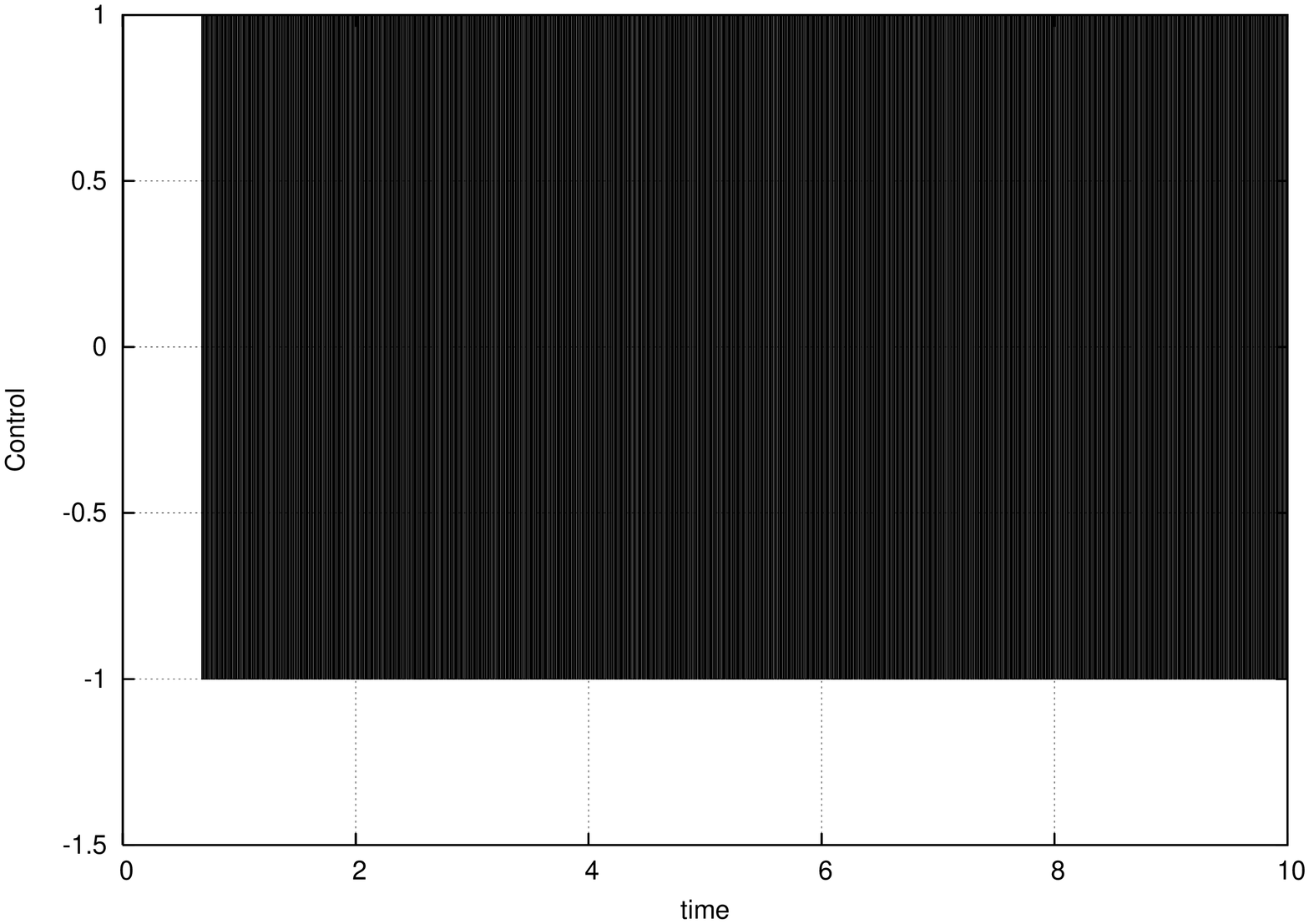}}
    \subfigure[error $e(t)$]{\label{Fig:SMC-EulerExplicitc}\includegraphics[angle=-90,width=0.49\textwidth]{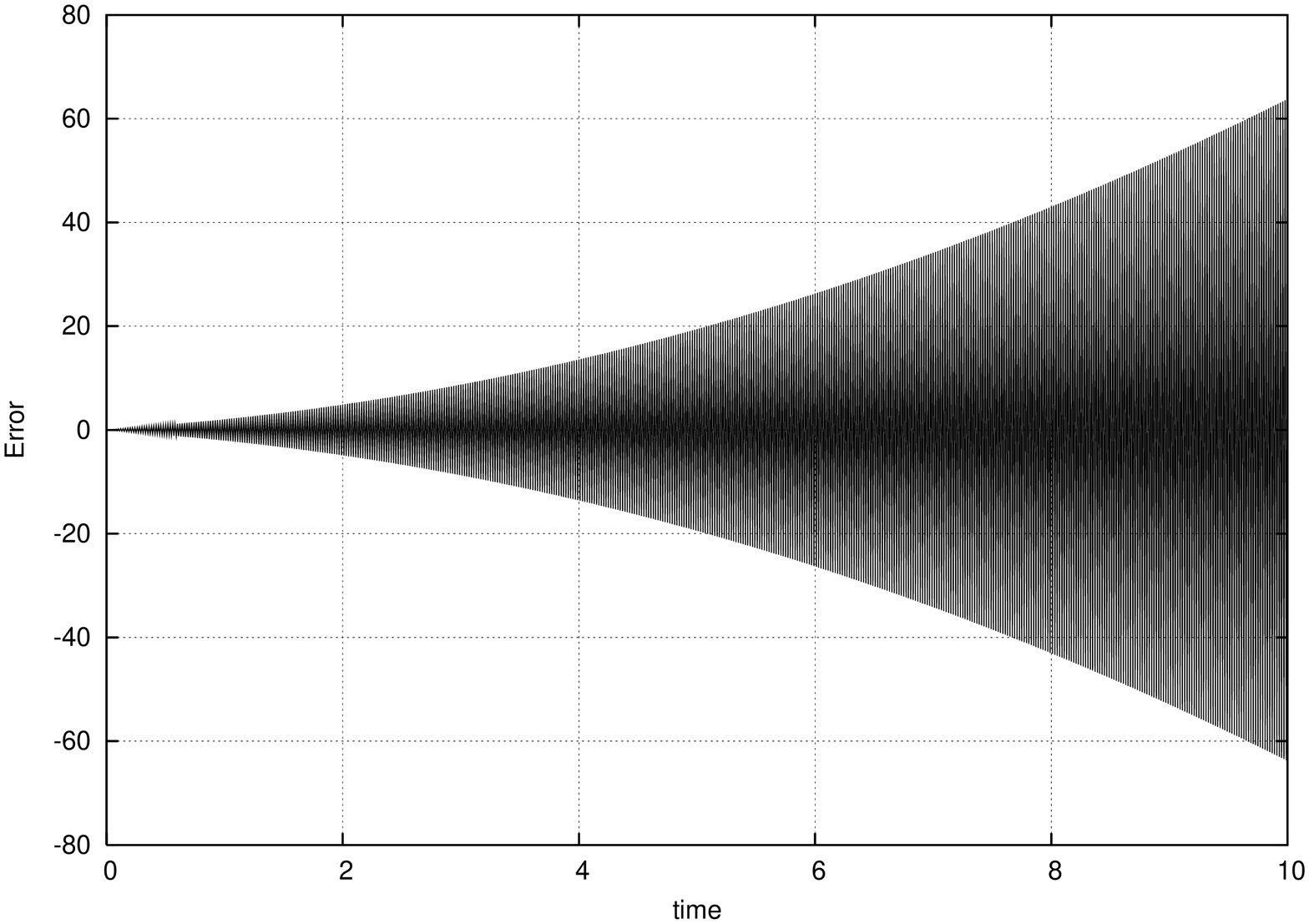}}
  \caption{Observer based SMC. $k=1$ and $\tau = 0.001$. $h = 0.01$  $\theta=0, \gamma=1$ }  \label{Fig:SMC-Control-EulerExplicit1}
\end{figure}

\paragraph{Influence of the integration parameters $\theta$ and $\gamma$} In the following numerical experiments, we discuss the role of the numerical parameters $\theta$ and $\gamma$. Due to the fact that the function $g(\cdot)$ is linear and reduced to a matrix $B = [1 ,0 ,0 , 0]^T$, the parameter $\gamma$ has no  influence on the numerical time--integration. On the contrary, the parameter $\theta$ has a huge influence on the stability of the integration. Indeed, the implicit Euler integration ($\theta =1)$ of the smooth term $f(\cdot)$ is unconditionally stable. This is not the case for the explicit Euler $\theta=0$ and for the chosen parameters $k$ and  $\tau$, the instability of the scheme for $h = 0.1$ does not allow to proceed to integration. On Figure~\ref{Fig:SMC-Control-EulerExplicit1}, the instability of the scheme is illustrated and appears as a chattering on the state $x$. The stability is retrieved for $h < 0.005$. 

\begin{figure}[htbp]
  \centering
  \psfrag{theta1}[][]{$\sf \theta=1$}
  \psfrag{thetamid}[][]{$\sf \theta=1/2$}
  \subfigure[state $x_1(t)$]{\label{Fig:SMC-Comparisona}\includegraphics[angle=-90,width=0.49\textwidth]{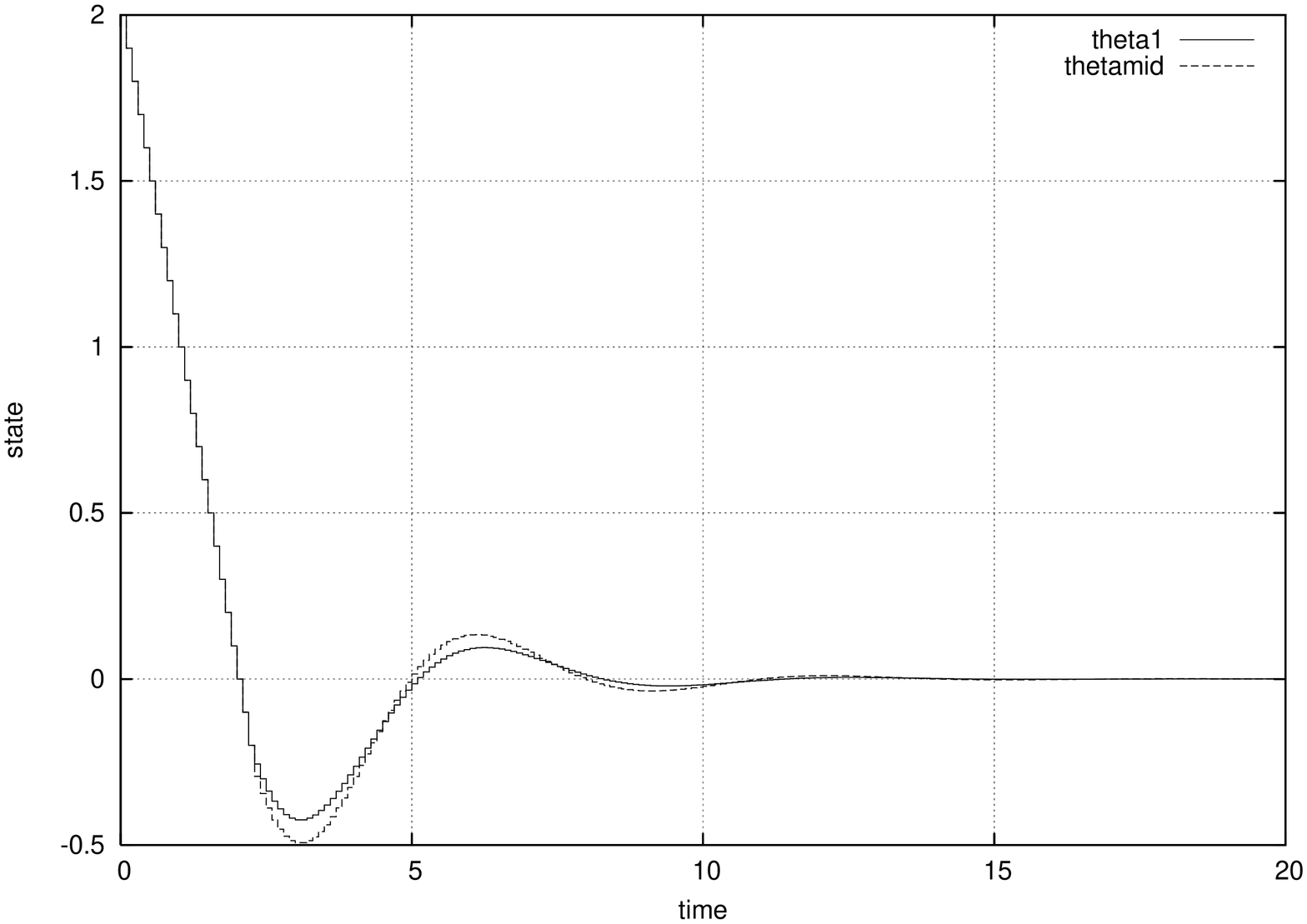}}
  \subfigure[control $u(t) = \mbox{Sgn}(Cx(t))$]{\label{Fig:SMC-Comparisonb}\includegraphics[angle=-90,width=0.49\textwidth]{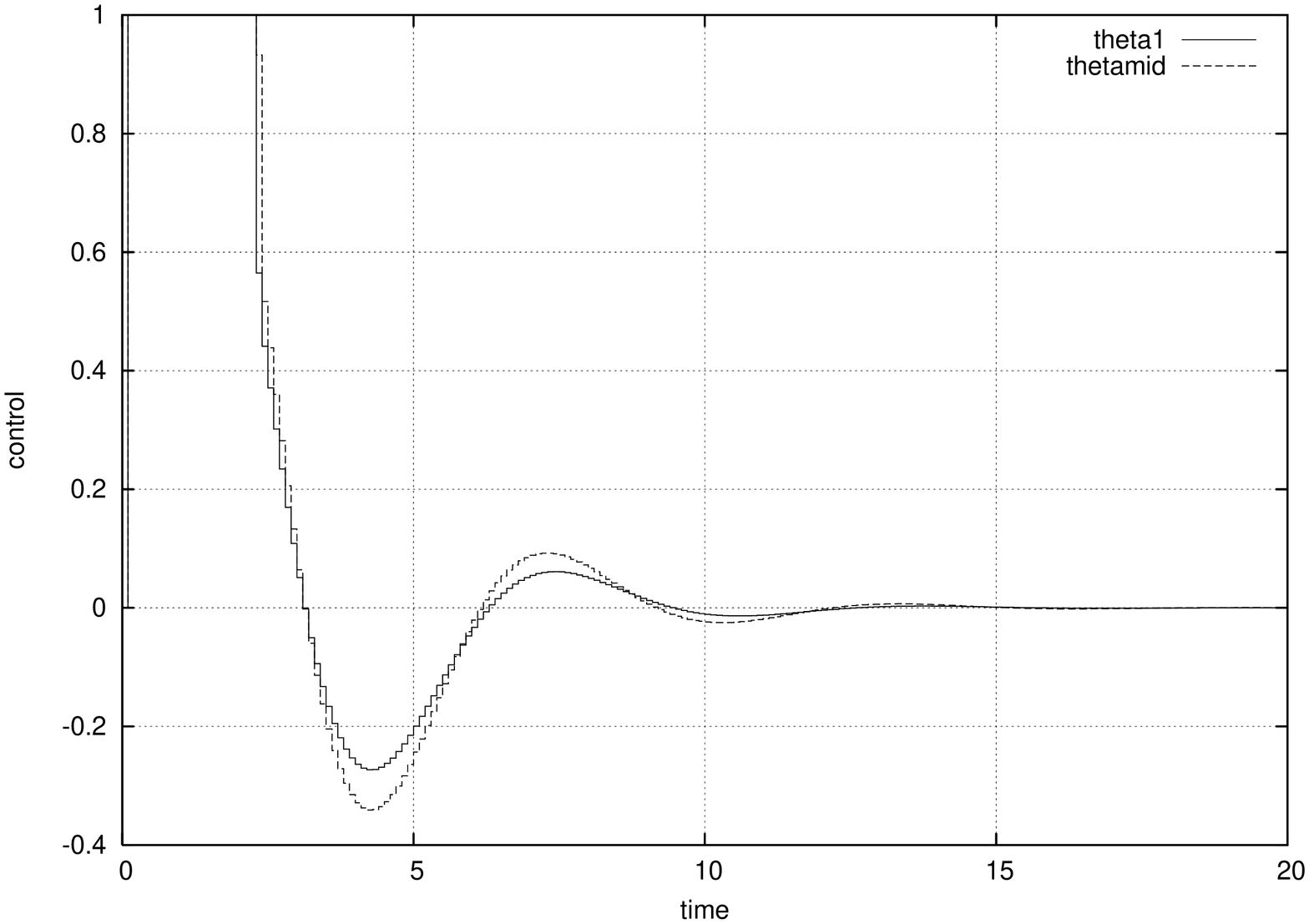}}
  \caption{Observer based SMC with $\theta=1$ and $\theta=1/2$. $k=1$ and $\tau = 0.5$. $h = 0.01,\gamma=1$ }  \label{Fig:SMC-Control-Comparison}
\end{figure}

Other choices of $\theta$ can be made to improve the numerical time--integration of the smooth dynamics. For instance, $\theta = 1/2$ yields a second order scheme for the integration of $f$. Unfortunately, the scheme is not of second order due to the fully implicit integration of the multi-valued part. Nevertheless, it can be interested to use of $\theta=1/2$ to decrease the numerical damping of the scheme on the smooth term. On the Figure~\ref{Fig:SMC-Control-Comparison}, two simulations of the observer based SMC are presented for $\theta =1$ and $\theta = 1/2$. The parameter $\tau$ has been modified to $0.5$ to correctly integrate  the parasitic dynamics with the same time--step $h=0.1$. 



\section{Conclusions}
\label{sectionConc}

In this paper the backward Euler method is studied on specific classes of Filippov's systems that encompass sliding mode control systems. It is shown that such implict schemes allow a smooth accurate stabilization on the sliding surface, even in case of codimension larger than one. Despite the backward Euler method has been studied and used for a long time in other fields like contact mechanics and electric circuits simulation \cite{acary2008}, it seems it has not yet been used in the sliding mode control community. This work therefore constitutes the introduction of a new discretization  method for EBC-SMC systems. The novelty compared to numerical simulation is that this time one has to consider not only the numerical simulation, but also the implementation on real processes. Perhaps one obstacle to the dissemination of the method is that at first sight, the controller designed from a backward philosophy looks like a non causal controller. However as shown in this paper this is not the case. This paper paves the way towards the study of a new family of discrete-time sliding mode controllers.

\bibliographystyle{plain}
\bibliography{Switch}

\tableofcontents

\end{document}